\documentclass[12pt,twoside]{article}

\usepackage[hmargin=0.56in,vmargin=0.56in]{geometry}
\usepackage{fancyhdr}
\usepackage{graphicx}
\usepackage{amssymb}
\usepackage{amsmath}
\usepackage{amsfonts}
\usepackage{amsthm}
\usepackage{mathrsfs}
\usepackage{mathtools}
\usepackage{bm}
\usepackage{color}
\usepackage{setspace}
\usepackage{exscale}
\usepackage{relsize}
\usepackage{epstopdf}
\usepackage{float}
\usepackage{cite}
\usepackage{nicefrac}
\usepackage{setspace}
\usepackage{algorithm,algorithmic}
\usepackage{tikz}
\usetikzlibrary{positioning}
\usepackage{xcolor}
\usepackage{rotating}
\usepackage{enumerate}
\usepackage{hyperref}
\usepackage{authblk}

\usepackage{textcomp}
\usepackage[title]{appendix}

\usepackage{cleveref}

\usepackage{booktabs,multirow} 
\usepackage{array} 
\usepackage{paralist} 
\usepackage{verbatim} 
\usepackage{subfigure} 

\usepackage[normalem]{ulem}

\newcommand{\ve}{\varepsilon}
\newcommand{\bmF}{\bm{\mathcal{F}}}
\newcommand{\bmG}{\bm{\mathcal{G}}}
\newcommand{\bmK}{\bm{\mathcal{K}}}
\newcommand{\bmL}{\bm{\mathcal{L}}}
\newcommand{\dx}{\Delta x} 
\newcommand{\dy}{\Delta y} 
\newcommand{\dt}{\Delta t} 
\newcommand{\U}{\bm U} 
\newcommand{\V}{\bm V} 
\newcommand{\F}{\bm F} 
\newcommand{\G}{\bm G} 
\newcommand{\K}{\bm K} 
\newcommand{\bL}{\bm L} 
\newcommand{\bS}{\bm S} 
\newcommand{\R}{\bm R} 

\newcommand{\B}{\bm B} 

\newcommand{\C}{\bm{C}} 

\newcommand{\hf}{\frac{1}{2}}

\newcommand{\jph}{{j+\frac{1}{2}}}
\newcommand{\jmh}{{j-\frac{1}{2}}}
\newcommand{\kph}{{k+\frac{1}{2}}}
\newcommand{\kmh}{{k-\frac{1}{2}}}

\newcommand{\p}{\mathcal{P}}

\newcommand{\zero}{\bm0}

\def\softd{{\leavevmode\setbox1=\hbox{d}%
          \hbox to 1.05\wd1{d\kern-0.4ex{\char039}\hss}}}

\newcommand\eref[1]{(\ref{#1})}

\newcommand*\xbar[1]{%
  \hbox{%
    \vbox{%
      \hrule height 0.5pt 
      \kern0.4ex
      \hbox{%
        \kern-0.05em
        \ensuremath{#1}%
        \kern-0.00em
      }%
    }%
  }%
}

\newtheorem{remark}{Remark}

\numberwithin{theorem}{section}
\numberwithin{equation}{section}
\numberwithin{remark}{section}
\numberwithin{figure}{section}
\numberwithin{table}{section}
\numberwithin{lemma}{section}

\graphicspath{{Figures/}}

\title{New Adaptive Numerical Methods Based on Dual Formulation of Hyperbolic Conservation Laws}

\author{Alina Chertock\thanks{Department of Mathematics, North Carolina State University, Raleigh, NC 27695, USA;
{\tt chertock@math.ncsu.edu}}, Qingcheng Fu\thanks{Department of Mathematics, Southern University of Science and Technology, Shenzhen,
518055, China; {\tt 12431005@mail.sustech.edu.cn}}, Alexander Kurganov\thanks{Department of Mathematics and Shenzhen International Center
for Mathematics, Southern University of Science and Technology, Shenzhen, 518055, China; {\tt alexander@sustech.edu.cn}}, and Lorenzo
Micalizzi\thanks{Department of Mathematics, North Carolina State University, Raleigh, NC 27695, USA; {\tt lmicali@ncsu.edu}}}

\date{}

\begin{document}

\maketitle

\begin{abstract}
In this paper, we propose an adaptive high-order method for hyperbolic systems of conservation laws. The proposed method is based on a dual
formulation approach: Two numerical solutions, corresponding to conservative and nonconservative formulations of the same system, are
evolved simultaneously. Since nonconservative schemes are known to produce nonphysical weak solutions near discontinuities, we exploit the
difference between these two solutions to construct a smoothness indicator (SI). In smooth regions, the difference between the conservative
and nonconservative solutions is of the same order as the truncation error of the underlying discretization, whereas in nonsmooth regions,
it is ${\cal O}(1)$. We apply this idea to the Euler equations of gas dynamics and define the SI using differences in the momentum and
pressure variables. This choice allows us to further distinguish neighborhoods of contact discontinuities from other nonsmooth parts of the
computed solution. The resulting classification is used to adaptively select numerical discretizations. In the vicinities of contact
discontinuities, we employ the low-dissipation central-upwind numerical flux and a second-order piecewise linear reconstruction with the
slopes computed using an overcompressive SBM limiter. Elsewhere, we use an alternative weighted essentially non-oscillatory (A-WENO)
framework with the central-upwind finite-volume numerical fluxes and either unlimited (in smooth regions) or Ai-WENO-Z (in the nonsmooth
regions away from contact discontinuities) fifth-order interpolation. Numerical results for the one- and two-dimensional compressible Euler
equations show that the proposed adaptive method improves both the computational efficiency and resolution of complex flow features compared
with the non-adaptive fifth-order A-WENO scheme.
\end{abstract}



\section{Introduction}
In this paper, we introduce new adaptive high-order numerical methods for hyperbolic systems of conservation laws within a dual formulation 
(DF) framework. In the two-dimensional (2-D) case, the governing equations can be expressed as 
\begin{equation}
\U_t+\F(\U)_x+\G(\U)_y=\zero,
\label{1.1}
\end{equation}
respectively. Here, $x$ and $y$ denote spatial variables, $t$ denotes time, $\U\in\mathbb R^d$ is the vector of conserved variables, and 
$\F:\mathbb R^d\to\mathbb R^d$ and $\G:\mathbb R^d\to\mathbb R^d$ are nonlinear flux functions.

Hyperbolic systems of conservation laws often generate complex wave structures, including shock, rarefaction, and contact waves, even when 
the initial data are smooth. The numerical approximation of such systems is particularly challenging, as the schemes should not only 
maintain high accuracy in smooth regions but also capture discontinuous structures without producing spurious oscillations. Over the past 
decades, a wide range of numerical schemes for such systems have been proposed; see, e.g., monographs and review papers 
\cite{GR3,Hes,LeV02,shu2009high,shu2020essentially,toro2009riemann}. 

It is well-known that in order to achieve high resolution of discontinuous solution structures on practically affordable meshes, one has to 
use numerical methods that are at least second-order accurate. Moreover, when the solution contains nontrivial smooth structures, one may
need to use higher-order methods to accurately resolve them. Among high-order numerical methods, weighted essentially non-oscillatory 
(WENO) schemes, originally introduced in \cite{JIANG1996202}, have become one of the most widely used approaches. However, a well-known 
drawback of high-order WENO schemes is the high computational cost of WENO reconstructions/interpolations, which employ nonlinear weights 
needed to ensure stability.

One way to reduce the computational cost is to construct adaptive schemes that employ WENO reconstructions/interpolations only in regions 
where the solution is ``rough'', while adopting some unlimited discretization elsewhere.  Such an adaptive approach was inroduced, for 
example, in \cite{CCK23}; see also \cite{MR2330302,MR3780797,WDGH2020,YDGW2020}. A key component of adaptive schemes is the design of an 
accurate, robust, and efficient smoothness indicator (SI) that automatically detects ``rough'' parts of the computed solution. A variety of 
SIs have been developed. Some of the SIs are based on data analyses that help measure the smoothness of the computed solution; see, e.g., 
\cite{ABD2008Vec,FS2017,GT2002,GT2006,LOHNER1987323,QS2005995,WDGWW2018,WSYK2015,WDGH2020}. There are also SIs, which measure the residual
and thus check how well the computed solution satisfies the studied system of PDEs; see, e.g., 
\cite{DKL2015,GUERMOND20114248,KK2005,KKP2002,KL2012,PS2011}. Another SI, introduced in \cite{CCK23}, is based on a time-Taylor expansion 
applied to the computed solutions at each stage of the Runge-Kutta time integrator. In this paper, we propose an alternative approach to designing SIs based on two numerical solutions evolved simultaneously within the DF framework recently introduced in~\cite{ACKM}; see also 
\cite{A2023,abgrall2026bound,AL2024,gaburro2024discontinuous,pidatella2019semi}.

We focus on numerical schemes that employ two formulations of the same governing equations, namely the conservative and nonconservative 
(primitive) forms. To this end, we rewrite \eref{1.1} as an equivalent (for smooth solutions) nonconservative system
\begin{equation}
\V_t+\widetilde{\F}(\V)_x+\widetilde{\G}(\V)_y=B(\V)\bm V_x+C(\V)\bm V_y.
\label{1.2}
\end{equation}
where  $\V\in\mathbb R^d$ is the vector of primitive variables, $\widetilde\F(\V):\mathbb R^d\to\mathbb R^d$ and
$\widetilde\G(\V):\mathbb R^d\to\mathbb R^d$ are nonlinear flux functions, and $B(\V),C(\V)\in\mathbb R^{d\times d}$ are matrices 
characterizing the nonconservative products. As a representative example, one may consider the Euler equations of gas dynamics, which can 
be expressed in the form \eref{1.1} with 
\begin{equation}
\begin{aligned}
&\U=(\rho,\rho u,\rho v,E)^\top\!\!,\\
&\F(\U)=\left(\rho u,\rho u^2+p,\rho uv,(E+p)u\right)^\top\!\!,~\G(\U)=\left(\rho v,\rho uv,\rho v^2+p,(E+p)v\right)^\top\!\!,
\end{aligned}
\label{1.3}
\end{equation}
where $\rho$ denotes the density, $u$ and $v$ are the velocities in the $x$- and $y$-directions, respectively, $p$ is the pressure, and $E$
stands for the total energy. The system \eref{1.1}, \eref{1.3} is completed using the equation of state, which in the ideal gas case reads
as
\begin{equation}
E=\frac{p}{\gamma-1}+\frac{\rho}{2}(u^2+v^2),
\label{1.4f}	
\end{equation}
where $\gamma$ is the specific heat ratio. This system can be reformulated as \eref{1.2} with 
\begin{equation}
\begin{aligned}			
&\V=(\rho,u,v,p)^\top,\quad\widetilde\F(\V)=\Big(\rho u,\frac{u^2}{2},0,pu\Big)^\top,\quad
\widetilde\G(\V)=\Big(\rho v,0,\frac{v^2}{2},pv\Big)^\top,\\[0.5ex]
&B(\V)=\begin{pmatrix}0&0&0&0\\0&0&0&-\frac{1}{\rho}\\0&0&-u&0\\0&-(\gamma-1)p&0&0\end{pmatrix},~
C(\V)=\begin{pmatrix}0&0&0&0\\0&-v&0&0\\0&0&0&-\frac{1}{\rho}\\0&0&-(\gamma-1)p&0\end{pmatrix}.
\end{aligned}
\label{1.4}	
\end{equation}

In the proposed adaptive DF methods, the solution is advanced in time as follows. At a given time level, we assume that we have both the
computed solution and the corresponding SI values, which enable us to classify different regions of the solution as ``rough'' or smooth. We
then advance the solution to the next time level by solving both the $\U$- and $\V$-systems, \eref{1.1} and \eref{1.2}, respectively. This
yields two approximations of the solution that are expected to differ only by truncation errors in smooth regions, but exhibit significant
discrepancies in nonsmooth regions, since nonconservative schemes generally converge to incorrect solutions in the presence of
discontinuities; see, e.g., \cite{MR2595792}. Motivated by this observation, we design a new SI based on the difference between these two
solutions.

We implement this SI for the Euler equations of gas dynamics. Specifically, we first consider differences in the momentum field to identify 
the ``rough'' regions, and then further examine differences in the pressure field to distinguish contact waves from the remaining ``rough'' 
areas, exploiting the fact that the pressure remains continuous across contact discontinuities.

Our adaptive strategy for solving the $\U$-system employs different numerical methods depending on the local flow features. In 
neighborhoods of contact waves, we follow \cite{MR4830367} and use the second-order low-dissipation central-upwind (LDCU) scheme proposed 
in \cite{chu2024new}, which incorporates the overcompressive limiter of \cite{SBM2003}. In the remaining nonsmooth regions, we apply the 
fifth-order alternative weighted essentially non-oscillatory (A-WENO) scheme introduced in \cite{MR4923663}. A-WENO schemes, originally 
proposed in \cite{Jiang13} (see also \cite{MR4923663,Gao20,Liu17}), can be combined with arbitrary finite-volume numerical fluxes and 
WENO-type interpolations. Specifically, we use the central-upwind (CU) numerical flux from \cite{kurganov2007reduction} together with the 
Ai-WENO-Z interpolation from \cite{WD22}. In smooth regions, we again employ the fifth-order A-WENO scheme with the CU numerical flux, but 
replace the nonlinear interpolation by an unlimited one.

The $\V$-system is integrated using a simplified version of the fifth-order A-WENO scheme from \cite{CHU2022111508}. It is important to
emphasize that the $\V$-system solution is nonphysical and is used solely for the construction of the new SI. Therefore, once the SI values 
are computed, we overwrite the nonconservative $\V$-solution with the physically relevant conservative $\U$-solution before proceeding to 
the next time step.

The rest of the paper is organized as follows. In \S\ref{sec2}, we introduce the general idea behind the proposed approach. Our specific 
applications to the one-dimensional (1-D) and 2-D compressible Euler equations are described in \S\ref{sec3} and \S\ref{sec4}, respectively.
The thorough numerical validation of the approach is presented in \S\ref{sec5}. Finally, \S\ref{sec6} outlines conclusions and further 
perspectives.

\section{New Adaptive Schemes within the DF Framework}\label{sec2}
We begin with a description of the general idea by considering a 2-D hyperbolic system of conservation laws \eref{1.1}, which admits a
nonconservative formulation \eref{1.2}.

We introduce a uniform Cartesian mesh consisting of nonoverlapping cells $I_{j,k}=[x_\jmh,x_\jph]\times[y_\kmh,y_\kph]$ of size $\dx\dy$
with $x_\jph-x_\jmh\equiv\dx$ and $y_\kph-y_\kmh\equiv\dy$, centered at $(x_j,y_k)$ with $x_j=\hf(x_\jmh+x_\jph)$ and
$y_k=\hf(y_\kmh+y_\kph)$, for $j=1,\dots,N$ and $k=1,\dots,M$. We denote by $\U_{j,k}(t)\approx\U(x_j,y_k,t)$ the point values of $\U$
computed at time $t$. These point values will be evolved in time using a one-step method, for which the solution update from the time level
$t^n$ to the next time level $t^{n+1}:=t^n+\dt$ can be, in general, written as
\begin{equation}
\underline\U^{n+1}={\cal C}(\underline\U^n,\underline{\bm{{\cal E}}}^n).
\label{2.1}
\end{equation}
Here, $\underline\U^n:=\{\U_{j,k}^n\}$ and ${\cal C}(\cdot,\cdot)$ is a high-order nonlinear conservative and consistent operator with an
adaptive nature, which depends on the second argument $\underline{\bm{{\cal E}}}^n:=\{\bm\ve_{j,k}^n\}$ with $\bm\ve_{j,k}^n$ representing
an SI for the solution in the cell $I_{j,k}$.

According to $\underline{\bm{{\cal E}}}^n$, we distinguish between smooth (S) and ``rough'' (R) areas, and within the latter one, we make a 
further distinction between the neighborhoods of contact discontinuities (RC) and the remaining ``rough'' non-contact areas (RNC). In each
of these areas, we adopt a different suitable discretization as described below:

\smallskip
\noindent
$\bullet$ In {\bf Region S}, we achieve high accuracy by evolving the solution using the fifth-order A-WENO scheme based on the CU 
numerical fluxes and unlimited interpolation;

\smallskip
\noindent
$\bullet$ In {\bf Region RNC}, we use the same fifth-order A-WENO scheme but, in order to guarantee a robust handling of discontinuities
avoiding spurious oscillations, we replace the unlimited interpolation with a local characteristic decomposition (LCD)-based fifth-order
Ai-WENO-Z interpolation;

\smallskip
\noindent
$\bullet$ In {\bf Region RC}, we switch to the second-order LDCU scheme and employ an overcompressive limiter applied to the local
characteristic variables in order to sharply capture linearly degenerate contact waves.

\begin{remark}
It is well-known that the numerical treatment of contact discontinuities is particularly challenging. In contrast to shocks, which are
supported by a self-stabilizing mechanism due to the convergence of characteristic curves and the resulting compression, linearly degenerate
contact discontinuities are inevitably smeared by the numerical dissipation present in any stable numerical scheme. Despite being formally
second-order accurate, the discretization employed in Region RC is well-suited for handling contact discontinuities. In fact, the LDCU
numerical flux is designed to accurately capture contact discontinuities and, when combined with the overcompressive limiter, the resulting
scheme is able to provide a very sharp resolution of discontinuous contact profiles. We stress, however, that overcompressive limiters must
be employed in a careful way, as their use in smooth regions typically leads to the generation of artificial structures that contain kinks
or spurious discontinuities.
\end{remark}

The SI $\underline{\bm{{\cal E}}}^n$ plays a crucial role in the proposed algorithm, and it is designed following ideas introduced in
\cite{chertock2025new}. For this purpose, we define the point values as approximations of the primitive variables
$\V_{j,k}(t):=\V(\U_{j,k}(t))$, with $\V(\cdot)$ representing the transformation from the conserved to primitive variables, the vector
$\underline\V^n:=\{\V_{j,k}^n\}$ with $\V^n_{j,k}:=\V(\U^n_{j,k})$, and we complement \eref{2.1} with the following
nonconservative update:
\begin{equation}
\underline\V^*={\cal N}(\underline\U^n,\underline\V^n),
\label{2.2}
\end{equation}
where ${\cal N}(\cdot,\cdot)$ is a nonlinear high-order consistent operator.

It is well-known that nonconservative schemes tend to converge to wrong weak solutions in the presence of discontinuities; see
\cite{MR2595792,hou1994nonconservative}. Therefore, the difference between $\U^{n+1}_{j,k}$ and $\U(\V^*_{j,k})$, where $\U(\cdot)$
represents the transformation from the primitive to conserved variables, is expected to be consistent with the order of the scheme in the
smooth areas only, while being ${\cal O}(1)$ in the ``rough'' areas. It is therefore natural to define the components of
$\underline{\bm{{\cal E}}}^{n+1}$ as functions of the differences $\bm\alpha(\U_{j,k}^{n+1})-\bm\alpha(\U(\V_{j,k}^*))$ with
$\bm\alpha(\cdot)$ being a vector of quantities of interest obtained from the conserved variables. We emphasize that the updated values
$\underline\V^*$ are, in general, unreliable and thus they are only used in the computation of the SI. Therefore, any linearly stable
high-order scheme can be used in \eref{2.2}. Our particular choice is a simplified version of the A-WENO scheme from \cite{CHU2022111508},
which will be described below.

\section{Application to the One-Dimensional Euler Equations}\label{sec3}
In this section, we apply the new adaptive scheme to the one-dimensional (1-D) Euler equations of gas dynamics, which read as
\begin{equation}
\begin{aligned}
&\U_t+\F(\U)_x=\zero,\\
&\U=(\rho,\rho u,E)^\top,\quad\F(\U)=\left(\rho u,\rho u^2+p,u(E+p)\right)^\top,\quad
E=\frac{p}{\gamma-1}+\frac{\rho u^2}{2}.
\end{aligned}
\label{3.1f}
\end{equation}
The primitive formulation of \eref{3.1f} reads as
\begin{equation}
\begin{aligned}
&\V_t+\widetilde\F(\V)_x=B(\V)\V_x,\quad\V=\begin{pmatrix}\rho\\u\\p\end{pmatrix},\quad\widetilde\F(\V)=
\begin{pmatrix}\rho u\\u^2/2\\pu\end{pmatrix},\\
&B(\V)=\begin{pmatrix}0&0&0\\0&0&-1/\rho\\0&-(\gamma-1)p&0\end{pmatrix}.
\end{aligned}
\label{3.3f}
\end{equation}
This system can be rewritten in the following quasi-conservative form:
\begin{equation}
\V_t+\K_x=\zero,\quad\K:=\widetilde{\F}(\V)-\R,\quad\R:=\int\limits^x_{\hat x}B(\V)\V_\xi(\xi,t)\,{\rm d}\xi,
\label{3.4f}
\end{equation}
where $\K$ is a global flux and $\hat x$ is an arbitrary number.

We introduce a uniform mesh consisting of the cells $I_j=[x_\jmh,x_\jph]$ of size $x_\jph-x_\jmh\equiv\dx$ with centers
$x_j=\hf(x_\jmh+x_\jph)$, for $j=1,\ldots,N$. We denote by $\U_j(t)\approx\U(x_j,t)$ and $\V_j(t):=\V(\U_j(t))$, and assume that at time
level $t^n$ the solution, realized in terms of $\U_j^n:=\U_j(t^n)$ and $\V_j^n:=\V_j(t^n)$, is available. Within the time interval
$[t^n,t^{n+1}]$, the solution is evolved according to the following semi-discretizations of the conservative $\U$-system \eref{3.1f} and
nonconservative $\V$-system \eref{3.3f}:
\begin{align}
\frac{{\rm d}}{{\rm d}t}\,\U_j(t)&=-\frac{\bmF_\jph(t)-\bmF_\jmh(t)}{\dx},\label{3.5f}\\[0.5ex]
\frac{{\rm d}}{{\rm d}t}\V_j(t)&=-\frac{\bmK_\jph(t)-\bmK_\jmh(t)}{\dx},\label{3.5ff}
\end{align}
where $\bmF_\jph$ and $\bmK_\jph$ are the corresponding numerical fluxes, which will be described below.

Notice that most of the indexed quantities are time-dependent, but from now on, we will omit this dependence for the sake of brevity.

\subsection{Adaptive Criteria}\label{sec31}
For the 1-D Euler equations, we denote $\bm\ve:=(\ve^{\rho u},\ve^p)^\top$ and define the grid values of $\ve^{\rho u}$ and $\ve^p$ as
\begin{equation*}
\ve^{\rho u}_j:=\left[\alpha^{\rho u}(\U_j)-\alpha^{\rho u}(\U(\V^*_j))\right]^2,\quad
\ve^p_j:=\left[\alpha^p(\U_j)-\alpha^p(\U(\V^*_j))\right]^2,
\end{equation*}
where we select $\bm\alpha=(\alpha^{\rho u},\alpha^p)^\top:=(\rho u,p)^\top$. We then smooth out these quantities by introducing
\begin{equation*}
\begin{aligned}
(\xbar{\ve^{\rho u}})_j&=\frac{\ve^{\rho u}_{j-2}+4\ve^{\rho u}_{j-1}+8\ve^{\rho u}_j+4\ve^{\rho u}_{j+1}+\ve^{\rho u}_{j+2}}{18},\\
(\xbar{\ve^p})_j&=\frac{\ve^p_{j-2}+4\ve^p_{j-1}+8\ve^p_j+4\ve^p_{j+1}+\ve^p_{j+2}}{18},
\end{aligned}
\end{equation*}
take their averages over the entire computational domain:
\begin{equation*}
(\xbar{\ve^{\rho u}})_{\rm ave}:=\frac{1}{N}\sum\limits^N_{j=1}(\xbar{\ve^{\rho u}})_j,\quad
(\xbar{\ve^p})_{\rm ave}:=\frac{1}{N}\sum\limits^N_{j=1}(\xbar{\ve^p})_j, 	
\end{equation*}
and define the following quantities at the cell interfaces:
\begin{equation*}
\ve^{\rho u}_\jph:=\max\left\{(\xbar{\ve^{\rho u}})_j,(\xbar{\ve^{\rho u}})_{j+1}\right\},\quad
\ve^p_\jph:=\max\left\{(\xbar{\ve^p})_j,(\xbar{\ve^p})_{j+1}\right\}.
\end{equation*}

We now select Regions S, RNC, and RC based on the consideration that while in the smooth regions both $\ve^{\rho u}_\jph$ and $\ve^p_\jph$
are supposed to be very small, $\ve^{\rho u}_\jph$ is expected to be ${\cal O}(1)$ near both shocks and contact discontinuities, while
$\ve^p_\jph$ is expected to be ${\cal O}(1)$ near the shocks only as $p$ stays continuous across the contact discontinuities. In particular,
the selection is carried out as follows:
\begin{equation*}
x_\jph\in
\begin{cases}
\mbox{Region S},&\mbox{if }\,\ve^{\rho u}_\jph<\kappa_{\rho u}\,(\xbar{\ve^{\rho u}})_{\rm ave},\\
\mbox{Region RC},&\mbox{if }\,\ve^{\rho u}_\jph>\kappa_{\rho u}\,(\xbar{\ve^{\rho u}})_{\rm ave}~~\mbox{and}~~
\ve^p_\jph<\kappa_p\,(\xbar{\ve^p})_{\rm ave},\\
\text{Region RNC},&\mbox{if }\,\ve^{\rho u}_\jph>\kappa_{\rho u}\,(\xbar{\ve^{\rho u}})_{\rm ave}~~\mbox{and}~~
\ve^p_\jph>\kappa_p\,(\xbar{\ve^p})_{\rm ave},
\end{cases}
\end{equation*}
where $\kappa_{\rho u}$ and $\kappa_p$ are two test-dependent adaption coefficients to be tuned.

\subsection{Numerical Solvers for the Conservative System}\label{sec32}
We use the semi-discretization \eref{3.5f} with numerical fluxes $\bmF_\jph$, which are obtained out of the finite-volume fluxes 
$\bmF^{\rm FV}_\jph:=\bmF^{\rm FV}\big(\U_\jph^-,\U_\jph^+\big)$, where $\U_\jph^\pm$ are point values at the cell interface $x_\jph$, obtained adaptively using a proper interpolation/reconstruction.

In particular, given the reconstructed/interpolated values at each cell interface (we will describe later how to obtain them), we first
compute the finite-volume numerical fluxes $\bmF^{\rm FV}_\jph$ as follows: if $x_\jph$ is in
either Regions S or RNC, $\bmF^{\rm FV}_\jph$ is the CU numerical flux from \cite{kurganov2007reduction}, while, if $x_\jph$ is in Region
RC, $\bmF^{\rm FV}_\jph$ is the LDCU numerical flux from \cite{chu2024new}.

We then proceed with computing the numerical fluxes $\bmF_\jph$. If $x_\jph$ is in either Regions S or RNC, we use the fifth-order A-WENO
numerical flux from \cite{MR4923663}, which reads as
\begin{equation}
\bmF_\jph=\bmF^{\rm FV}_\jph-\frac{1}{24}(\dx)^2(\F_{xx})_\jph+\frac{7}{5760}(\dx)^4(\F_{xxxx})_\jph,
\label{3.6ff}
\end{equation}
where the high-order correction terms are computed through central differences
\begin{equation}
\begin{aligned}
&(\F_{xx})_\jph=\frac{-\bmF^{\rm FV}_{j-\frac{3}{2}}+16\bmF^{\rm FV}_\jmh-30\bmF^{\rm FV}_\jph+16\bmF^{\rm FV}_{j+\frac{3}{2}}-
\bmF^{\rm FV}_{j+\frac{5}{2}}}{12(\dx)^2},\\
&(\F_{xxxx})_\jph=\frac{\bmF^{\rm FV}_{j-\frac{3}{2}}-4\bmF^{\rm FV}_\jmh+6\bmF^{\rm FV}_\jph-4\bmF^{\rm FV}_{j+\frac{3}{2}}+
\bmF^{\rm FV}_{j+\frac{5}{2}}}{(\dx)^4},
\end{aligned}
\label{3.7ff}
\end{equation}
applied to the already computed finite-volume numerical fluxes. If $x_\jph$ is in Region RC, we set $\bmF_\jph=\bmF^{\rm FV}_\jph$ to fully
exploit the ability of the second-order LDCU scheme to capture contact discontinuities.
\begin{remark}
We do not provide details on the CU and LDCU numerical fluxes for the sake of brevity.
\end{remark}

We compute the point values $\U_\jph^\pm$ as follows. If $x_\jph\in$ Region S, we use the fourth-degree polynomials and compute the 
fifth-order interpolants over the stencils $\{x_{j-2},x_{j-1},x_j,x_{j+1},x_{j+2}\}$ and $\{x_{j-1},x_j,x_{j+1},x_{j+2},x_{j+3}\}$ to 
obtain
\begin{equation}
\bm\U_\jph^-=\frac{3\bm\U_{j-2}-20\bm\U_{j-1}+90\bm\U_j+60\bm\U_{j+1}-5\bm\U_{j+2}}{128},
\label{3.8f}
\end{equation}
and
\begin{equation}
\bm\U_\jph^+=\frac{-5\bm\U_{j-1}+60\bm\U_j+90\bm\U_{j+1}-20\bm\U_{j+2}+3\bm\U_{j+3}}{128},
\label{3.9f}
\end{equation}
respectively.

In ``rough'' regions, in order to minimize spurious oscillations, which can be triggered by the use of either a high-order interpolation or
an overcompressive limiter, we perform the interpolation/reconstruction in the local characteristic variables $\bm\Gamma$, which are
obtained using the LCD. To this end, we first compute averaged values $\widehat{\U}_\jph$ (in the numerical examples presented in
\S\ref{sec5}, we have used the Roe averages, and define the local linearized Jacobians
$\widehat A_\jph:=\frac{\partial\F}{\partial\U}(\widehat{\U}_\jph)$. These matrices can be diagonalized using the matrices of their
eigenvectors, $R_\jph$, so that $R_\jph^{-1}\widehat A_\jph R_\jph$ are diagonal matrices. We then introduce the local characteristic
variables in the neighborhood of $x=x_\jph$,
\begin{equation}
\bm\Gamma_\ell:=R_\jph^{-1}\U_\ell,\quad\ell=j-2,\dots,j+3,
\label{3.10ff}
\end{equation}
and distinguish between the two different ``rough'' regions.

\smallskip
\noindent
$\bullet$ If $x_\jph\in$ Region RNC, we use the values $\Gamma_\ell$ in \eref{3.10ff} to compute the point values $\bm\Gamma_\jph^\pm$ with
the help of the fifth-order Ai-WENO-Z interpolation, whose detailed description can be found in \cite[\S2.2]{MR4938975}.

\smallskip
\noindent
$\bullet$ If $x_\jph\in$ Region RC, we employ the overcompressive limiter to evaluate the slopes
\begin{equation}
\begin{aligned}
&(\bm\Gamma_x)_j=\phi^{\rm SBM}_{\theta,\tau}\left(\frac{\bm\Gamma_{j+1}-\bm\Gamma_j}{\bm\Gamma_j-\bm\Gamma_{j-1}}\right)
\frac{\bm\Gamma_j-\bm\Gamma_{j-1}}{\dx},\\[1ex]
&(\bm\Gamma_x)_{j+1}=\phi^{\rm SBM}_{\theta,\tau}
\left(\frac{\bm\Gamma_{j+2}-\bm\Gamma_{j+1}}{\bm\Gamma_{j+1}-\bm\Gamma_j}\right)\frac{\bm\Gamma_{j+1}-\bm\Gamma_j}{\dx},
\end{aligned}
\label{3.6f}
\end{equation}
where the two-parameter SBM function
\begin{equation}
\phi^{\rm SBM}_{\theta,\tau}(r):=\begin{cases}0&\mbox{if $r<0$},\\\min\{r\theta,1+\tau(r-1)\}&\mbox{if}~0<r\le1,\\
r\phi^{\rm SBM}_{\theta,\tau}(1/r)&\text{otherwise,}\end{cases}
\label{3.7f}
\end{equation}
is applied in the component-wise manner. We set $\theta=2$ and $\tau=-0.25$, which correspond to the overcompressive regime; see
\cite{SBM2003}. Equipped with \eref{3.6f}, we evaluate
$$
\bm\Gamma^-_\jph=\bm\Gamma_j+\frac{\dx}{2}(\bm\Gamma_x)_j\quad\mbox{and}\quad
\bm\Gamma^+_\jph=\bm\Gamma_{j+1}-\frac{\dx}{2}(\bm\Gamma_x)_{j+1}.
$$

Finally, in both Regions RNC and RC, we obtain the corresponding point values of $\U$ by
\begin{equation*}
\U^\pm_\jph=R_\jph\bm\Gamma^\pm_\jph.
\end{equation*}

\subsection{Numerical Solver for the Primitive System}\label{sec33}
We employ the semi-discretization \eref{3.5ff} with the fifth-order global numerical fluxes
\begin{equation}
\bmK_\jph=\bmK_\jph^{\rm FV}-\frac{1}{24}(\dx)^2(\K_{xx})_\jph+\frac{7}{5760}(\dx)^4(\K_{xxxx})_\jph,
\label{3.12f}
\end{equation}
where $\bmK_\jph^{\rm FV}$ is a simplified version of the finite-volume PCCU numerical global flux from \cite{CKX_Beijing} and the
high-order correction terms are computed the same way as in \eref{3.6ff}--\eref{3.7ff}, namely,
\begin{equation}
\begin{aligned}
&(\K_{xx})_\jph=\frac{-\bmK^{\rm FV}_{j-\frac{3}{2}}+16\bmK^{\rm FV}_\jmh-30\bmK^{\rm FV}_\jph+
16\bmK^{\rm FV}_{j+\frac{3}{2}}-\bmK^{\rm FV}_{j+\frac{5}{2}}}{12(\dx)^2},\\
&(\K_{xxxx})_\jph=\frac{\bmK^{\rm FV}_{j-\frac{3}{2}}-4\bmK^{\rm FV}_\jmh+6\bmK^{\rm FV}_\jph-
4\bmK^{\rm FV}_{j+\frac{3}{2}}+\bmK^{\rm FV}_{j+\frac{5}{2}}}{(\dx)^4}.
\end{aligned}
\label{3.13f}
\end{equation}
The finite-volume numerical global fluxes are given by
\begin{equation}
\bmK_\jph^{\rm FV}=\frac{a_\jph^+\K_\jph^--a_\jph^-\K_\jph^+}{a_\jph^+-a_\jph^-}+\frac{a_\jph^+a_\jph^-}{a_\jph^+-a_\jph^-}
\left(\V_\jph^+-\V_\jph^-\right),
\label{3.12}
\end{equation}
where $\V_\jph^\pm$ are one-sided cell interface point values. We stress that our goal is to enforce the linear stability of the primitive
system solver, and hence $\V_\jph^\pm$ are obtained using the same fifth-order unlimited interpolation as in \eref{3.8f}--\eref{3.9f}, 
which gives
$$
\begin{aligned}
&\bm\V_\jph^-=\frac{3\bm\V_{j-2}-20\bm\V_{j-1}+90\bm\V_j+60\bm\V_{j+1}-5\bm\V_{j+2}}{128},\\[0.5ex]
&\bm\V_\jph^+=\frac{-5\bm\V_{j-1}+60\bm\V_j+90\bm\V_{j+1}-20\bm\V_{j+2}+3\bm\V_{j+3}}{128}.
\end{aligned}
$$
In \eref{3.12}, $a_\jph^\pm$ are the one-sided local speeds of propagation estimated by
\begin{equation}
\begin{aligned}
&a_\jph^+=\max\left\{u_\jph^-+c_\jph^-,u_\jph^++c_\jph^+,\delta\right\},\\
&a_\jph^-=\min\left\{u_\jph^--c_\jph^-,u_\jph^+-c_\jph^+,-\delta\right\},
\end{aligned}
\label{3.15}
\end{equation}
where $c:=\sqrt{\gamma p/\rho}$ is the speed of sound, $u_\jph^\pm$ and $c_\jph^\pm$ are computed from $\bm\V_\jph^\pm$, and
$\delta=10^{-10}$ is a small constant needed to prevent division by zero.

A simplification compared with the numerical global flux from \cite{CKX_Beijing} is in the evaluation of $\K_\jph^\pm$, which is now
conducted by neglecting a special path-conservative treatment of the interface jump associated with the nonconservative product terms. 
Namely, we set
\begin{equation}
\K_\jph^\pm=\widetilde{\F}\big(\V_\jph^\pm\big)-\R_\jph,
\label{3.16}
\end{equation}
where the global part of the numerical flux is evaluated recursively as follows. First, we set $\hat x=x_\hf$ so that $\R_\hf=\zero$, and
then we compute
\begin{equation}
\R_\jph=\R_\jmh+\B_j,\quad j=1,\dots,N,
\label{3.17}
\end{equation}
where $\B_j$ is a fifth-order approximation of the integral of the nonconservative product terms over the cell $I_j$:
\begin{equation}
\B_j\approx\int\limits_{I_\jph}B(\V)\V_x\,{\rm d}x.
\label{3.10f}
\end{equation}
To evaluate these integrals, we use the fifth-order Newton-Cotes quadrature with the nodes $x_\jmh$, $x_{j-\frac{1}{4}}$, $x_j$,
$x_{j+\frac{1}{4}}$, $x_\jph$; see \cite[Appendix A]{MR4938975} for details. When this quadrature is applied, one needs to interpolate the
values $\V_{j\pm\frac{1}{4}}$ at the points $x=x_{j\pm\frac{1}{4}}$. This is done with the help of the fifth-order (fourth-degree)
polynomial interpolant over the stencil $\{x_{j-2},x_{j-1},x_j,x_{j+1},x_{j+2}\}$, which gives
\begin{equation*}
\begin{aligned}
\bm\V_{j-\frac{1}{4}}&=\frac{-45\bm\V_{j-2}+420\bm\V_{j-1}+1890\bm\V_j-252\bm\V_{j+1}+35\bm\V_{j+2}}{2048},\\[0.5ex]
\bm\V_{j+\frac{1}{4}}&=\frac{35\bm\V_{j-2}-252\bm\V_{j-1}+1890\bm\V_j+420\bm\V_{j+1}-45\bm\V_{j+2}}{2048}.
\end{aligned}
\end{equation*}
\begin{remark}
We stress again that the solution $\underline\V^*$, obtained upon completion of one time evolution step according to the simplified and only
linearly stable semi-discrete scheme \eref{3.5ff}, \eref{3.12f}--\eref{3.10f}, will be accurate in the smooth regions only, while it may
contain severe oscillations near the discontinuities. This behavior is not only acceptable, but, in fact, advantageous since the solution
$\underline\V^*$ is used exclusively to compute the SI and thus to detect Regions S, RC, and RNC.
\end{remark}

\section{Application to the Two-Dimensional Euler Equations}\label{sec4}
In this section, we apply the new adaptive scheme to the 2-D Euler equations of gas dynamics, whose conservative and primitive formulations
are given by \eref{1.1}, \eref{1.3}--\eref{1.4f}  and \eref{1.2}, \eref{1.4}, respectively. As in the 1-D case, the primitive system can be
rewritten in the quasi-conservative form as
\begin{equation}
\begin{aligned}
&\V_t+\K_x+\bL_y=\zero,\quad\K=\widetilde\F(\V)-\R,\quad\bL=\widetilde\G(\V)-\bS,\\\
&\R:=\int\limits^x_{\hat x}B(\V)\V_\xi(\xi,y,t)\,{\rm d}\xi,\quad\bS:=\int\limits^y_{\hat y}C(\V)\V_\eta(x,\eta,t)\,{\rm d}\eta,
\end{aligned}
\label{3.6}
\end{equation}
where $\K$ and $\bL$ are global fluxes and $\hat x$ and $\hat y$ are arbitrary numbers.

Within the time interval $[t^n,t^{n+1}]$, the solution is evolved according to the following semi-discretizations of the conservative
$\U$-system \eref{1.1}, \eref{1.3}--\eref{1.4f} and nonconservative $\V$-system \eref{1.2}, \eref{1.4}:
\begin{align}
\frac{{\rm d}}{{\rm d}t}\,\U_{j,k}&=-\frac{\bmF_{\jph,k}-\bmF_{\jmh,k}}{\dx}-\frac{\bmG_{j,\kph}-\bmG_{j,\kmh}}{\dy},\label{4.2}\\
\frac{{\rm d}}{{\rm d}t}\,\V_{j,k}&=-\frac{\bmK_{\jph,k}-\bmK_{\jmh,k}}{\dx}-\frac{\bmL_{j,\kph}-\bmL_{j,\kmh}}{\dy},\label{4.3}
\end{align}
where $\bmF_{\jph,k},\bmG_{j,\kph}$ and $\bmK_{\jph,k},\bmL_{j,\kph}$ are the corresponding numerical fluxes, which will be described below.

\subsection{Adaptive Criteria}
For the 2-D Euler equations, we define the grid values of $\ve^{\rho u}$, $\ve^{\rho v}$, and $\ve^p$ as
\begin{equation*}
\begin{aligned}
&\ve^{\rho u}_{j,k}:=\left[\alpha^{\rho u}(\U_{j,k})-\alpha^{\rho u}(\U(\V^*_{j,k}))\right]^2,\quad
\ve^{\rho v}_{j,k}:=\left[\alpha^{\rho v}(\U_{j,k})-\alpha^{\rho v}(\U(\V^*_{j,k}))\right]^2,\\
&\ve^p_{j,k}:=\left[\alpha^p(\U_{j,k})-\alpha^p(\U(\V^*_{j,k}))\right]^2,
\end{aligned}
\end{equation*}
and select $\bm\alpha=(\alpha^{\rho u},\alpha^{\rho v},\alpha^p)^\top:=(\rho u,\rho v,p)^\top$, denote
$\bm\ve:=(\ve^{\rho u},\ve^{\rho v},\ve^p)^\top$, and define the grid values of $\ve^{\rho u}$, $\ve^{\rho v}$, and $\ve^p$ as

In the 2-D case, we extend the adaptive criteria introduced in \S\ref{sec31} in a dimension-by-dimension manner. We therefore check
smoothness of the computed solution at the midpoints of the cell interfaces in each of the directions separately.

We begin with the cell interfaces in the $ x$-direction. We first smooth out $\ve^{\rho u}_{j,k}$ and $\ve^p_{j,k}$ by introducing
\begin{equation*}
\begin{aligned}
&(\xbar{\ve^{\rho u}})_{j,k}^x=
\frac{\ve^{\rho u}_{j-2,k}+4\ve^{\rho u}_{j-1,k}+8\ve^{\rho u}_{j,k}+4\ve^{\rho u}_{j+1,k}+\ve^{\rho u}_{j+2,k}}{18},\\
&(\xbar{\ve^p})_{j,k}^x=\frac{\ve^p_{j-2,k}+4\ve^p_{j-1,k}+8\ve^p_{j,k}+4\ve^p_{j+1,k}+\ve^p_{j+2,k}}{18},
\end{aligned}
\end{equation*}
take their averages over the entire computational domain:
\begin{equation*}
(\xbar{\ve^{\rho u}})_{\rm ave}^x:=\frac{1}{NM}\sum\limits^N_{j=1}\sum\limits^M_{k=1}(\xbar{\ve^{\rho u}})_{j,k}^x,\quad
(\xbar{\ve^p})_{\rm ave}^x:=\frac{1}{NM}\sum\limits^N_{j=1}\sum\limits^M_{k=1}(\xbar{\ve^p})_{j,k}^x,       
\end{equation*}
and define the following quantities at $(x_\jph,y_k)$:
\begin{equation*}
\ve^{\rho u}_{\jph,k}:=\max\left\{(\xbar{\ve^{\rho u}})_{j,k}^x,(\xbar{\ve^{\rho u}})_{j+1,k}^x\right\},\quad
\ve^p_{\jph,k}:=\max\left\{(\xbar{\ve^p})_{j,k}^x,(\xbar{\ve^p})_{j+1,k}^x\right\}.
\end{equation*}
We then use the following simple strategy for automatically performing the classification of the $x$-directional interfaces:
\begin{equation*}
(x_\jph,y_k)\in
\begin{cases}
\mbox{Region S},&\mbox{if }\,\ve^{\rho u}_{\jph,k}<\kappa_{\rho u}\,(\xbar{\ve^{\rho u}})_{\rm ave}^x,\\
\mbox{Region RC},&\mbox{if }\,\ve^{\rho u}_{\jph,k}>\kappa_{\rho u}\,(\xbar{\ve^{\rho u}})_{\rm ave}^x~~\mbox{and}~~
\ve^p_{\jph,k}<\kappa_p\,(\xbar{\ve^p})_{\rm ave}^x,\\
\text{Region RNC},&\mbox{if }\,\ve^{\rho u}_{\jph,k}>\kappa_{\rho u}\,(\xbar{\ve^{\rho u}})_{\rm ave}^x~~\mbox{and}~~
\ve^p_{\jph,k}>\kappa_p\,(\xbar{\ve^p})_{\rm ave}^x,
\end{cases}
\end{equation*}
where $\kappa_{\rho u}$ and $\kappa_p$, as before, are two test-dependent adaption coefficients to be tuned.

Next, we proceed with the cell interfaces in the $y$-direction in a similar manner. We first compute
\begin{equation*}
\begin{aligned}
&(\xbar{\ve^{\rho v}})_{j,k}^y=
\frac{\ve^{\rho v}_{j,k-2}+4\ve^{\rho v}_{j,k-1}+8\ve^{\rho v}_{j,k}+4\ve^{\rho v}_{j,k+1}+\ve^{\rho v}_{j,k+2}}{18},\\
&(\xbar{\ve^p})_{j,k}^y=\frac{\ve^p_{j,k-2}+4\ve^p_{j,k-1}+8\ve^p_{j,k}+4\ve^p_{j,k+1}+\ve^p_{j,k+2}}{18},
\end{aligned}
\end{equation*}
\begin{equation*}
(\xbar{\ve^{\rho v}})_{\rm ave}^y:=\frac{1}{NM}\sum\limits^N_{j=1}\sum\limits^M_{k=1}(\xbar{\ve^{\rho v}})_{j,k}^y,\quad
(\xbar{\ve^p})_{\rm ave}^y:=\frac{1}{NM}\sum\limits^N_{j=1}\sum\limits^M_{k=1}(\xbar{\ve^p})_{j,k}^y,       
\end{equation*}
and
\begin{equation*}
\ve^{\rho v}_{j,\kph}:=\max\left\{(\xbar{\ve^{\rho v}})_{j,k}^y,(\xbar{\ve^{\rho v}})_{j,k+1}^y\right\},\quad
\ve^p_{j,\kph}:=\max\left\{(\xbar{\ve^p})_{j,k}^y,(\xbar{\ve^p})_{j,k+1}^y\right\},
\end{equation*}
and then perform the classification of the $y$-directional interfaces:
\begin{equation*}
(x_j,y_\kph)\in
\begin{cases}
\mbox{Region S},&\mbox{if }\,\ve^{\rho v}_{j,\kph}<\kappa_{\rho v}\,(\xbar{\ve^{\rho v}})_{\rm ave}^y,\\
\mbox{Region RC},&\mbox{if }\,\ve^{\rho v}_{j,\kph}>\kappa_{\rho v}\,(\xbar{\ve^{\rho v}})_{\rm ave}^y~~\mbox{and}~~
\ve^p_{j,\kph}<\kappa_p\,(\xbar{\ve^p})_{\rm ave}^y,\\
\text{Region RNC},&\mbox{if }\,\ve^{\rho u}_{j,\kph}>\kappa_{\rho v}\,(\xbar{\ve^{\rho v}})_{\rm ave}^y~~\mbox{and}~~
\ve^p_{j,\kph}>\kappa_p\,(\xbar{\ve^p})_{\rm ave}^y,
\end{cases}
\end{equation*}
where $\kappa_{\rho v}$ and $\kappa_p$ are two test-dependent adaption coefficients to be tuned.

\subsection{Numerical Solvers for the Conservative System}
We employ the semi-discretization \eref{4.2} with numerical fluxes $\bmF_{\jph,k}$ and $\bmG_{j,\kph}$, obtained out of pre-computed 
second-order finite-volume ones $\bmF^{\rm FV}_{\jph,k}:=\bmF^{\rm FV}\big(\U_{\jph,k}^-,\U_{\jph,k}^+\big)$ and
$\bmG^{\rm FV}_{j,\kph}:=\bmG^{\rm FV}\big(\U_{j,\kph}^-,\U_{j,\kph}^+\big)$, where the point values $\U_{\jph,k}^\pm$ and $\U_{j,\kph}^\pm$
are reconstructed/interpolated on both sides at $(x_\jph,y_k)$ and $(x_j,y_\kph)$, respectively. These values are obtained adaptively using
proper reconstruction/interpolation, as detailed later.

We first compute the finite-volume numerical fluxes $\bmF^{\rm FV}_{\jph,k}$ in the $x$-direction as follows: if $(x_\jph,y_k)$ is in either
Regions S or RNC, $\bmF^{\rm FV}_{\jph,k}$ is the 2-D CU numerical flux from \cite{kurganov2007reduction}, while, if $(x_\jph,y_k)$ is in
Region RC, $\bmF^{\rm FV}_{\jph,k}$ is the 2-D LDCU numerical flux from \cite{chu2024new}.

We then compute the numerical fluxes $\bmF_{\jph,k}$. If $(x_\jph,y_k)$ is in either Regions S or RNC, we use the fifth-order A-WENO
numerical flux from \cite{MR4923663}, which reads as
\begin{equation}
\bmF_{\jph,k}=\bmF^{\rm FV}_{\jph,k}-\frac{1}{24}(\dx)^2(\F_{xx})_{\jph,k}+\frac{7}{5760}(\dx)^4(\F_{xxxx})_{\jph,k},
\label{4.4}
\end{equation}
where the high-order correction terms are computed through central differences:
\begin{equation}
\begin{aligned}
&(\F_{xx})_{\jph,k}=\frac{-\bmF^{\rm FV}_{j-\frac{3}{2},k}+16\bmF^{\rm FV}_{\jmh,k}-30\bmF^{\rm FV}_{\jph,k}+
16\bmF^{\rm FV}_{j+\frac{3}{2},k}-\bmF^{\rm FV}_{j+\frac{5}{2},k}}{12(\dx)^2},\\
&(\F_{xxxx})_{\jph,k}=\frac{\bmF^{\rm FV}_{j-\frac{3}{2},k}-4\bmF^{\rm FV}_{\jmh,k}+6\bmF^{\rm FV}_{\jph,k}-
4\bmF^{\rm FV}_{j+\frac{3}{2},k}+\bmF^{\rm FV}_{j+\frac{5}{2},k}}{(\dx)^4}.
\end{aligned}
\label{4.5}
\end{equation}
applied to the already computed finite-volume numerical fluxes. If $(x_\jph,y_k)$ is in Region RC, we set
$\bmF_{\jph,k}=\bmF^{\rm FV}_{\jph,k}$, coinciding with the 2-D second-order LDCU numerical flux from \cite{chu2024new}.

The numerical fluxes in the $y$-direction are obtained similarly. We start by computing the finite-volume fluxes $\bmG^{\rm FV}_{j,\kph}$ as
follows: if $(x_j,y_\kph)$ is in either Regions S or RNC, $\bmG^{\rm FV}_{j,\kph}$ is the 2-D CU numerical flux from
\cite{kurganov2007reduction}, while, if $(x_j,y_\kph)$ is in Region RC, $\bmG^{\rm FV}_{j,\kph}$ is the 2-D LDCU numerical flux from
\cite{chu2024new}.

We then compute the numerical fluxes $\bmG_{j,\kph}$. If $(x_j,y_\kph)$ is in either Regions S or RNC, we use the fifth-order A-WENO
numerical flux from \cite{MR4923663}, which reads as
\begin{equation}
\bmG_{j,\kph}=\bmG^{\rm FV}_{j,\kph}-\frac{1}{24}(\dy)^2(\G_{yy})_{j,\kph}+\frac{7}{5760}(\dy)^4(\G_{yyyy})_{j,\kph},
\label{4.6}
\end{equation}
where the high-order correction terms are
\begin{equation}
\begin{aligned}
&(\G_{yy})_{j,\kph}=\frac{-\bmG^{\rm FV}_{j,k-\frac{3}{2}}+16\bmG^{\rm FV}_{j,\kmh}-30\bmG^{\rm FV}_{j,\kph}+
16\bmG^{\rm FV}_{j,k+\frac{3}{2}}-\bmG^{\rm FV}_{j,k+\frac{5}{2}}}{12(\dy)^2},\\
&(\G_{yyyy})_{j,\kph}=\frac{\bmG^{\rm FV}_{j,k-\frac{3}{2}}-4\bmG^{\rm FV}_{j,\kmh}+6\bmG^{\rm FV}_{j,\kph}-
4\bmG^{\rm FV}_{j,k+\frac{3}{2}}+\bmG^{\rm FV}_{j,k+\frac{5}{2}}}{(\dy)^4}.
\end{aligned}
\label{4.7}
\end{equation}
If $(x_j,y_\kph)$ is in Region RC, we set $\bmG_{j,\kph}=\bmG^{\rm FV}_{j,\kph}$, coinciding with the 2-D second-order LDCU numerical flux
from \cite{chu2024new}.
\begin{remark}
We do not provide details on the 2-D CU and LDCU numerical fluxes for the sake brevity.
\end{remark}

Let us now address the computation of the required point values at cell interfaces. We obtain the point values $\U_{\jph,k}^\pm$ and 
$\U_{j,\kph}^\pm$as follows. If $(x_\jph,y_k)\in$ Region S, we use the 1-D fourth-degree polynomials and compute fifth-order interpolants 
over the stencils $\{(x_{j-2},y_k),(x_{j-1},y_k),(x_j,y_k),(x_{j+1},y_k),(x_{j+2},y_k)\}$ and 
$\{(x_{j-1},y_k),(x_j,y_k),(x_{j+1},y_k),(x_{j+2},y_k),(x_{j+3},y_k)\}$ in the $x$-direction to obtain
\begin{equation}
\bm\U_{\jph,k}^-=\frac{3\bm\U_{j-2,k}-20\bm\U_{j-1,k}+90\bm\U_{j,k}+60\bm\U_{j+1,k}-5\bm\U_{j+2,k}}{128},
\label{4.8}
\end{equation}
and
\begin{equation}
\bm\U_{\jph,k}^+=\frac{-5\bm\U_{j-1,k}+60\bm\U_{j,k}+90\bm\U_{j+1,k}-20\bm\U_{j+2,k}+3\bm\U_{j+3,k}}{128},
\label{4.9}
\end{equation}
respectively. Acting similarly in the $y$-direction, we obtain
\begin{equation}
\begin{aligned}
&\bm\U_{j,\kph}^-=\frac{3\bm\U_{j,k-2}-20\bm\U_{j,k-1}+90\bm\U_{j,k}+60\bm\U_{j,k+1}-5\bm\U_{j,k+2}}{128},\\[0.5ex]
&\bm\U_{j,\kph}^+=\frac{-5\bm\U_{j,k-1}+60\bm\U_{j,k}+90\bm\U_{j,k+1}-20\bm\U_{j,k+2}+3\bm\U_{j,k+3}}{128}.
\end{aligned}
\label{4.10}
\end{equation}

In ``rough'' regions, we proceed as in the 1-D case and perform the interpolation/reconstruction in local characteristic variables
$\bm\Gamma$, which are obtained using the LCD in a dimension-by-dimension manner. In the $x$-direction, we first compute for all $k$
averaged values $\widehat{\U}_{\jph,k}$ (in the numerical examples presented in \S\ref{sec5}, we have used the Roe averages) and define the
local linearized Jacobians $\widehat A_{\jph,k}:=\frac{\partial\F}{\partial\U}(\widehat{\U}_{\jph,k})$, whose eigenvectors form the matrices
$R_{\jph,k}$. We then introduce the local characteristic variables in the neighborhood of $(x_\jph,y_k)$,
$$
\bm\Gamma_\ell:=R_{\jph,k}^{-1}\U_{\ell,k},\quad\ell=j-2,\dots,j+3,
$$
and obtain the values $\bm\Gamma_\jph^\pm$ using the 1-D adaptive interpolation/reconstruction described in \S\ref{sec32}. After that, we
obtain the corresponding point values of $\U$ by
\begin{equation*}
\U^\pm_{\jph,k}=R_{\jph,k}\bm\Gamma^\pm_\jph.
\end{equation*}

Finally, we proceed in the $y$-direction in a similar way to compute $\widehat{\U}_{j,\kph}$ in ``rough'' regions.

\subsection{Numerical Solver for the Primitive System}
We employ the semi-discretization \eref{4.3} with the 2-D fifth-order global numerical fluxes
\begin{equation}
\begin{aligned}
&\bmK_{\jph,k}=\bmK_{\jph,k}^{\rm FV}-\frac{1}{24}(\dx)^2(\K_{xx})_{\jph,k}+\frac{7}{5760}(\dx)^4(\K_{xxxx})_{\jph,k},\\
&\bmL_{j,\kph}=\bmL_{j,\kph}^{\rm FV}-\frac{1}{24}(\dy)^2(\bL_{yy})_{j,\kph}+\frac{7}{5760}(\dy)^4(\bL_{yyyy})_{j,\kph},
\end{aligned}
\label{4.11}
\end{equation}
where $\bmK_{\jph,k}^{\rm FV}$ and $\bmL_{j,\kph}^{\rm FV}$ are dimension-by-dimension extensions of the 1-D numerical global fluxes
introduced in \S\ref{sec33} and the high-order correction terms are
\begin{equation}
\begin{aligned}
&(\K_{xx})_{\jph,k}=\frac{-\bmK^{\rm FV}_{j-\frac{3}{2},k}+16\bmK^{\rm FV}_{\jmh,k}-30\bmK^{\rm FV}_{\jph,k}+
16\bmK^{\rm FV}_{j+\frac{3}{2},k}-\bmK^{\rm FV}_{j+\frac{5}{2},k}}{12(\dx)^2},\\
&(\K_{xxxx})_{\jph,k}=\frac{\bmK^{\rm FV}_{j-\frac{3}{2},k}-4\bmK^{\rm FV}_{\jmh,k}+6\bmK^{\rm FV}_{\jph,k}-
4\bmK^{\rm FV}_{j+\frac{3}{2},k}+\bmK^{\rm FV}_{j+\frac{5}{2},k}}{(\dx)^4},\\
&(\bL_{yy})_{j,\kph}=\frac{-\bmL^{\rm FV}_{j,k-\frac{3}{2}}+16\bmL^{\rm FV}_{j,\kmh}-30\bmL^{\rm FV}_{j,\kph}+
16\bmL^{\rm FV}_{j,k+\frac{3}{2}}-\bmL^{\rm FV}_{j,k+\frac{5}{2}}}{12(\dy)^2},\\
&(\bL_{xxxx})_{j,\kph}=\frac{\bmL^{\rm FV}_{j,k-\frac{3}{2}}-4\bmL^{\rm FV}_{j,\kmh}+6\bmL^{\rm FV}_{j,\kph}-
4\bmL^{\rm FV}_{j,k+\frac{3}{2}}+\bmL^{\rm FV}_{j,k+\frac{5}{2}}}{(\dy)^4}.
\end{aligned}
\label{4.12}
\end{equation}
The finite-volume numerical global fluxes are given by
\begin{equation}
\begin{aligned}
\bmK_{\jph,k}^{\rm FV}&=\frac{a_{\jph,k}^+\K_{\jph,k}^--a_{\jph,k}^-\K_{\jph,k}^+}{a_{\jph,k}^+-a_{\jph,k}^-}\\
&+\frac{a_{\jph,k}^+a_{\jph,k}^-}{a_{\jph,k}^+-a_{\jph,k}^-}\left(\V_{\jph,k}^+-\V_{\jph,k}^-\right),\\[0.5ex]
\bmL_{j,\kph}^{\rm FV}&=\frac{a_{j,\kph}^+\bL_{j,\kph}^--a_{j,\kph}^-\bL_{j,\kph}^+}{a_{j,\kph}^+-a_{j,\kph}^-}\\
&+\frac{a_{j,\kph}^+a_{j,\kph}^-}{a_{j,\kph}^+-a_{j,\kph}^-}\left(\V_{j,\kph}^+-\V_{j,\kph}^-\right),
\end{aligned}
\label{4.13}
\end{equation}
where $\V_{\jph,k}^\pm$ and $\V_{j,\kph}^\pm$ are one-sided point values at $(x_\jph,y_k)$ and $(x_j,y_\kph)$, respectively. As in the 1-D 
case, these values are computed using the unlimited fifth-order (fourth-degree) interpolations performed in the $x$- and $y$-directions
separately. This results in
\allowdisplaybreaks
\begin{align*}
&\bm\V_{\jph,k}^-=\frac{3\bm\V_{j-2,k}-20\bm\V_{j-1,k}+90\bm\V_{j,k}+60\bm\V_{j+1,k}-5\bm\V_{j+2,k}}{128},\\[0.5ex]
&\bm\V_{\jph,k}^+=\frac{-5\bm\V_{j-1,k}+60\bm\V_{j,k}+90\bm\V_{j+1,k}-20\bm\V_{j+2,k}+3\bm\V_{j+3,k}}{128},\\[0.5ex]
&\bm\V_{j,\kph}^-=\frac{3\bm\V_{j,k-2}-20\bm\V_{j,k-1}+90\bm\V_{j,k}+60\bm\V_{j,k+1}-5\bm\V_{j,k+2}}{128},\\[0.5ex]
&\bm\V_{j,\kph}^+=\frac{-5\bm\V_{j,k-1}+60\bm\V_{j,k}+90\bm\V_{j,k+1}-20\bm\V_{j,k+2}+3\bm\V_{j,k+3}}{128}.
\end{align*}
In \eref{4.13}, $a_{\jph,k}^\pm$ and $a_{j,\kph}^\pm$ are the one-sided local speeds of propagation in the $x$- and $y$-directions,
respectively. The speeds are estimated by
\begin{equation}
\begin{aligned}
&a_{\jph,k}^+=\max\left\{u_{\jph,k}^-+c_{\jph,k}^-,u_{\jph,k}^++c_{\jph,k}^+,\delta\right\},\\
&a_{\jph,k}^-=\min\left\{u_{\jph,k}^--c_{\jph,k}^-,u_{\jph,k}^+-c_{\jph,k}^+,-\delta\right\},\\
&a_{j,\kph}^+=\max\left\{v_{j,\kph}^-+c_{j,\kph}^-,v_{j,\kph}^++c_{j,\kph}^+,\delta\right\},\\
&a_{j,\kph}^-=\min\left\{v_{j,\kph}^--c_{j,\kph}^-,v_{j,\kph}^+-c_{j,\kph}^+,-\delta\right\},
\end{aligned}
\label{4.14}
\end{equation}
where $u_{\jph,k}^\pm$ and $c_{\jph,k}^\pm$ are computed from $\bm\V_{\jph,k}^\pm$, $v_{j,\kph}^\pm$ and $c_{j,\kph}^\pm$ are computed from
$\bm\V_{j,\kph}^\pm$, and $\delta=10^{-10}$.

Finally, the global fluxes $\K_{\jph,k}^\pm$ and $\bL_{j,\kph}^\pm$ used in \eref{4.13} are obtained by
\begin{equation}
\K_{\jph,k}^\pm=\widetilde\F\big(\V_{\jph,k}^\pm\big)-\R_{\jph,k},\quad\bL_{j,\kph}^\pm=\widetilde\G\big(\V_{j,\kph}^\pm\big)-\bS_{j,\kph},
\label{4.15}
\end{equation}
where the global part of the numerical flux is evaluated recursively by setting $\hat x=x_\hf$ and $\hat y=y_\hf$ so that
$\R_{\hf,k}=\bS_{j,\hf}=\zero$, and then computing
\begin{equation}
\begin{aligned}
&\R_{\jph,k}=\R_{\jmh,k}+\B_{j,k},&&j=1,\dots,N,\\
&\bS_{j,\kph}=\bS_{j,\kmh}+\C_{j,k},&&k=1,\dots,M,
\end{aligned}
\label{4.16}
\end{equation}
where $\B_j$ and $\C_{j,k}$ are fifth-order approximations of the integrals of the nonconservative product terms:
\begin{equation}
\B_{j,k}\approx\int\limits_{x_\jmh}^{x_\jph}B(\V)\V_x\,{\rm d}x,\quad\C_{j,k}\approx\int\limits_{y_\kmh}^{y_\kph}C(\V)\V_y\,{\rm d}y.
\label{4.17}
\end{equation}
We evaluate the first of these integrals using the fifth-order Newton-Cotes quadrature with the nodes
$(x_\jmh,y_k)$, $(x_{j-\frac{1}{4}},y_k)$, $(x_j,y_k)$, $(x_{j+\frac{1}{4}},y_k)$, $(x_\jph,y_k)$. When this quadrature is applied, one
needs to interpolate the values $\V_{j\pm\frac{1}{4},k}$ at the points $(x_{j\pm\frac{1}{4}},y_k)$. This is done using the unlimited 
fifth-order (fourth-degree) polynomial interpolant over the stencil
$\{(x_{j-2},y_k),(x_{j-1},y_k),(x_j,y_k),(x_{j+1},y_k),(x_{j+2},y_k)\}$, which gives
\begin{equation}
\begin{aligned}
\bm\V_{j-\frac{1}{4},k}&=\frac{-45\bm\V_{j-2,k}+420\bm\V_{j-1,k}+1890\bm\V_{j,k}-252\bm\V_{j+1,k}+35\bm\V_{j+2,k}}{2048},\\[0.5ex]
\bm\V_{j+\frac{1}{4},k}&=\frac{35\bm\V_{j-2,k}-252\bm\V_{j-1,k}+1890\bm\V_{j,k}+420\bm\V_{j+1,k}-45\bm\V_{j+2,k}}{2048}.
\end{aligned}
\label{4.18}
\end{equation}
The second integral in \eref{4.17} is evaluated using the same fifth-order Newton-Cotes quadrature, which would require the values
$\V_{j,k\pm\frac{1}{4}}$ at the points $(x_j,y_{k\pm\frac{1}{4}})$, which are computed similarly to \eref{4.18}:
\begin{equation*}
\begin{aligned}
\bm\V_{j,k-\frac{1}{4}}&=\frac{-45\bm\V_{j,k-2}+420\bm\V_{j,k-1}+1890\bm\V_{j,k}-252\bm\V_{j,k+1}+35\bm\V_{j,k+2}}{2048},\\[0.5ex]
\bm\V_{j,k+\frac{1}{4}}&=\frac{35\bm\V_{j,k-2}-252\bm\V_{j,k-1}+1890\bm\V_{j,k}+420\bm\V_{j,k+1}-45\bm\V_{j,k+2}}{2048}.
\end{aligned}
\end{equation*}

\section{Numerical Examples}\label{sec5}
In this section, we illustrate the performance of the proposed adaptive approach on benchmarks for the 1-D (\S\ref{sec51}) and 2-D
(\S\ref{sec52}) Euler equations of gas dynamics.

In all of the numerical examples:

\smallskip
\noindent
$\bullet$ The three-stages third-order strong stability-preserving (SSP) Runge-Kutta method \cite{GST,GKS} is employed to numerically
integrate the ODE systems \eref{3.5f}, \eref{3.5ff}, \eref{4.2}, and \eref{4.3};

\smallskip
\noindent
$\bullet$ The time step is selected with the CFL number $0.45$; 

\smallskip
\noindent
$\bullet$ The first time evolution step is performed by the non-adaptive A-WENO schemes, which coincides with the A-WENO schemes used in the
Regions RNC, but applied throughout the entire computational domain. These non-adaptive methods will be referred to as A-WENO schemes;

\smallskip
\noindent
$\bullet$ $\gamma=1.4$, except for Example 6, in which $\gamma=5/3$;

\smallskip
\noindent
$\bullet$ In addition, the symmetry-enforcing algorithm from \cite[Appendix B]{WANG2020SIAM} is implemented in 2-D examples.

\smallskip
In order to demonstrate the advantages of the proposed adaptive methods, we compare the obtained results with those computed by the 
non-adaptive A-WENO schemes.
\begin{remark}
While conducting the numerical experiments, we have observed that the distribution of smooth and ``rough'' areas over the computational
domain varies only slightly over a short time interval. Therefore, computational cost can be reduced by avoiding shock detection at every
time step. In fact, in all of the numerical experiments reported below, the shock detection is performed only once every three time steps
and the Regions S, RC, and RNC are kept unchanged until the next shock detection is carried out.
\end{remark}
\begin{remark}
In the case the solution is smooth, our adaptive scheme reduces to the fifth-order A-WENO scheme with the unlimited reconstruction. We have
verified that the experimental order of convergence (EOC) is the expected fifth one on several 1-D and 2-D accuracy tests, not reported here
for the sake of brevity. We stress that when the EOC was checked, the time step was adjusted to be proportional to $(\dx)^\frac{5}{3}$ in
order to guarantee a matching between the temporal and spatial accuracy.
\end{remark}

\subsection{One-Dimensional Examples}\label{sec51}
The 1-D tests discussed in this section include two shock-turbulence interaction problems characterized by several smooth structures and a
blast wave problem featuring strong discontinuities.

\smallskip
\paragraph{Example 1 (Shock-Density Wave Interaction Problem)}
In this example, we consider the shock–density wave interaction problem originally proposed in \cite{SHU1988439}. The initial conditions,
\begin{equation*}
(\rho(x,0),u(x,0),p(x,0))=
\begin{cases}\bigg(\dfrac{27}{7},\dfrac{4\sqrt{35}}{9},\dfrac{31}{3}\bigg),&x<-4,\\(1+0.2\sin(5x),0,1),&x>-4,\end{cases}
\end{equation*}
are prescribed in the computational domain $[-5,15]$ subject to the free boundary conditions. In this case, a Mach 3 shock interacts with
small-amplitude density fluctuations ahead of it, generating a sequence of intricate waves behind the shock.

\begin{figure}[ht!]
\centerline{\includegraphics[trim=0.1cm 0.0cm 0.1cm 0.3cm, clip, width=0.35\textwidth]{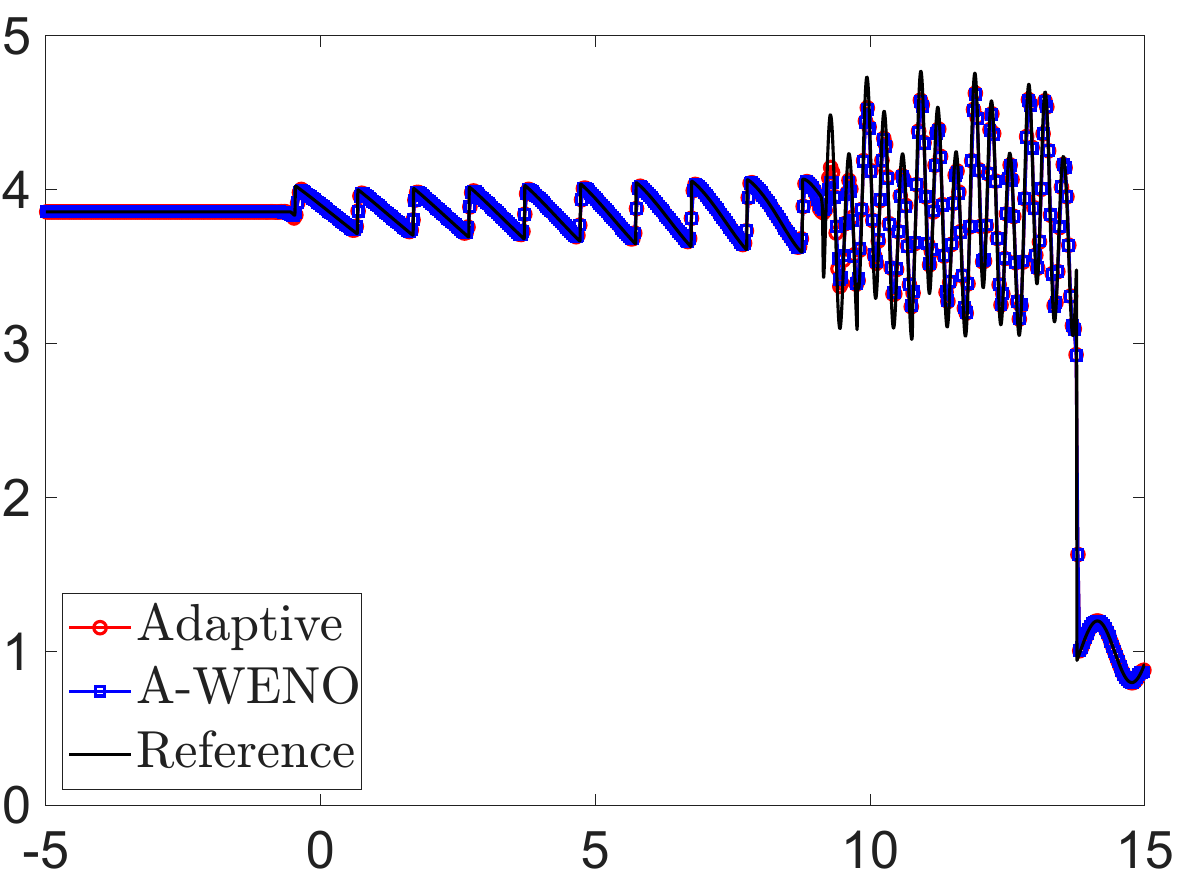}\hspace*{0.8cm}
            \includegraphics[trim=0.1cm 0.1cm 0.1cm 0.0cm, clip, width=0.36\textwidth]{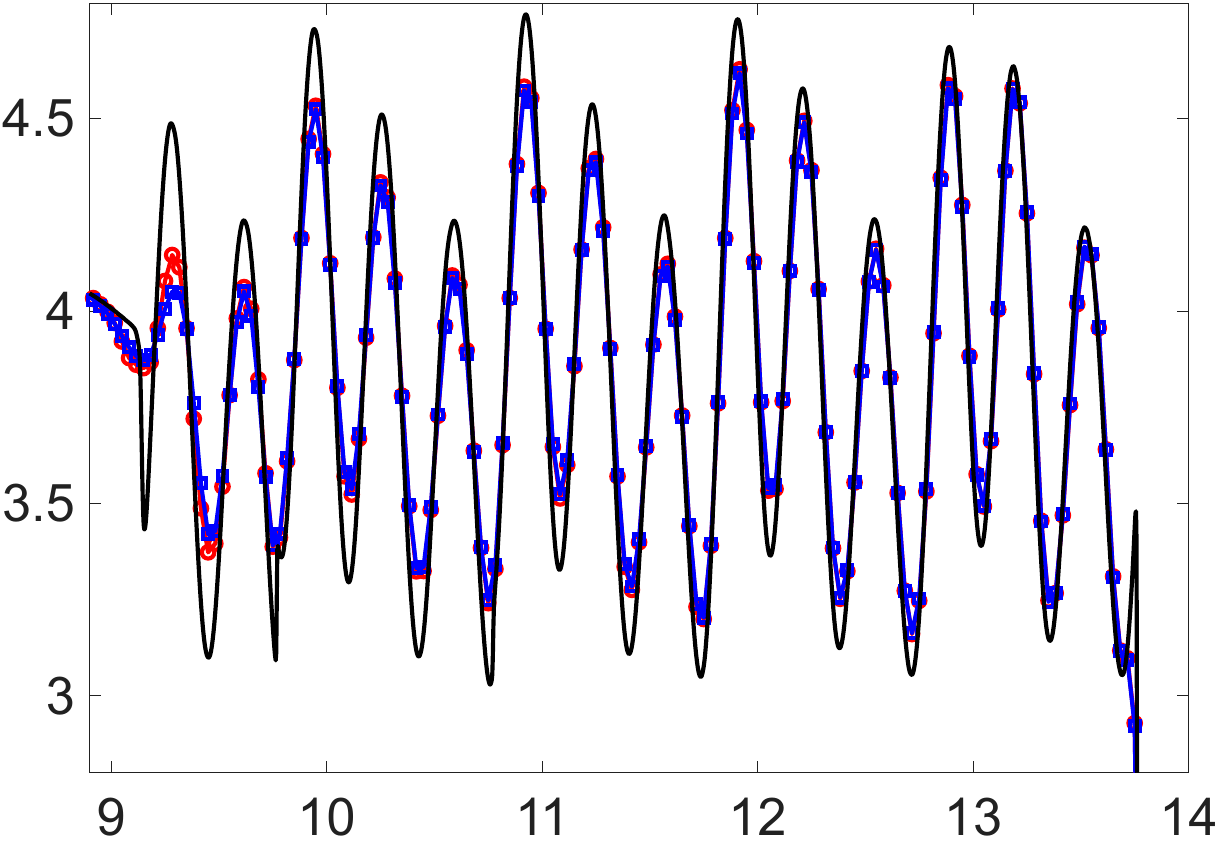}}
\caption{\sf Example 1: Density $\rho$ computed by the adaptive and A-WENO scheme with $\dx=1/30$ (left) and zoom at $x\in[8.9,14]$ (right).
\label{fig51}}
\end{figure}

\begin{figure}[ht!]
\centerline{\includegraphics[trim=0.1cm 0.0cm 0.0cm 0.1cm, clip, width=0.35\textwidth]{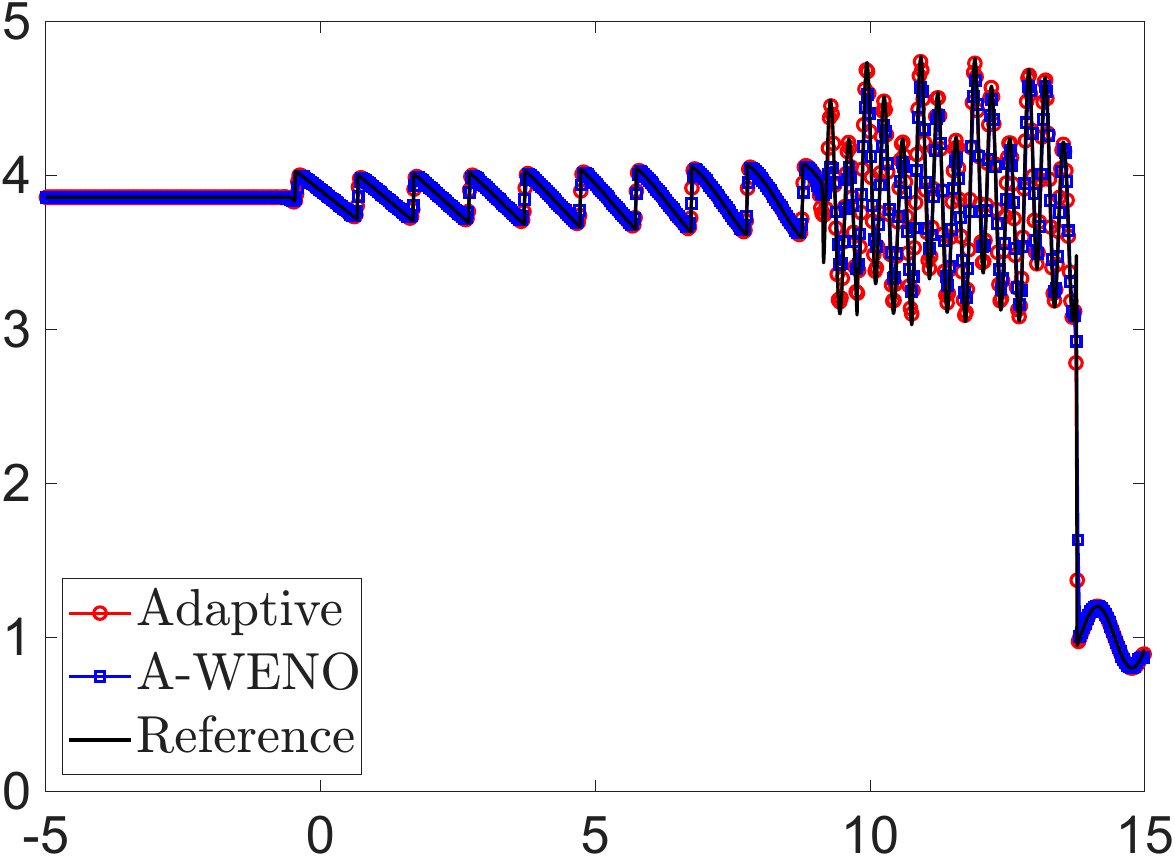}\hspace*{0.8cm}
	    \includegraphics[trim=0.1cm 0.1cm 0.1cm 0.0cm, clip, width=0.36\textwidth]{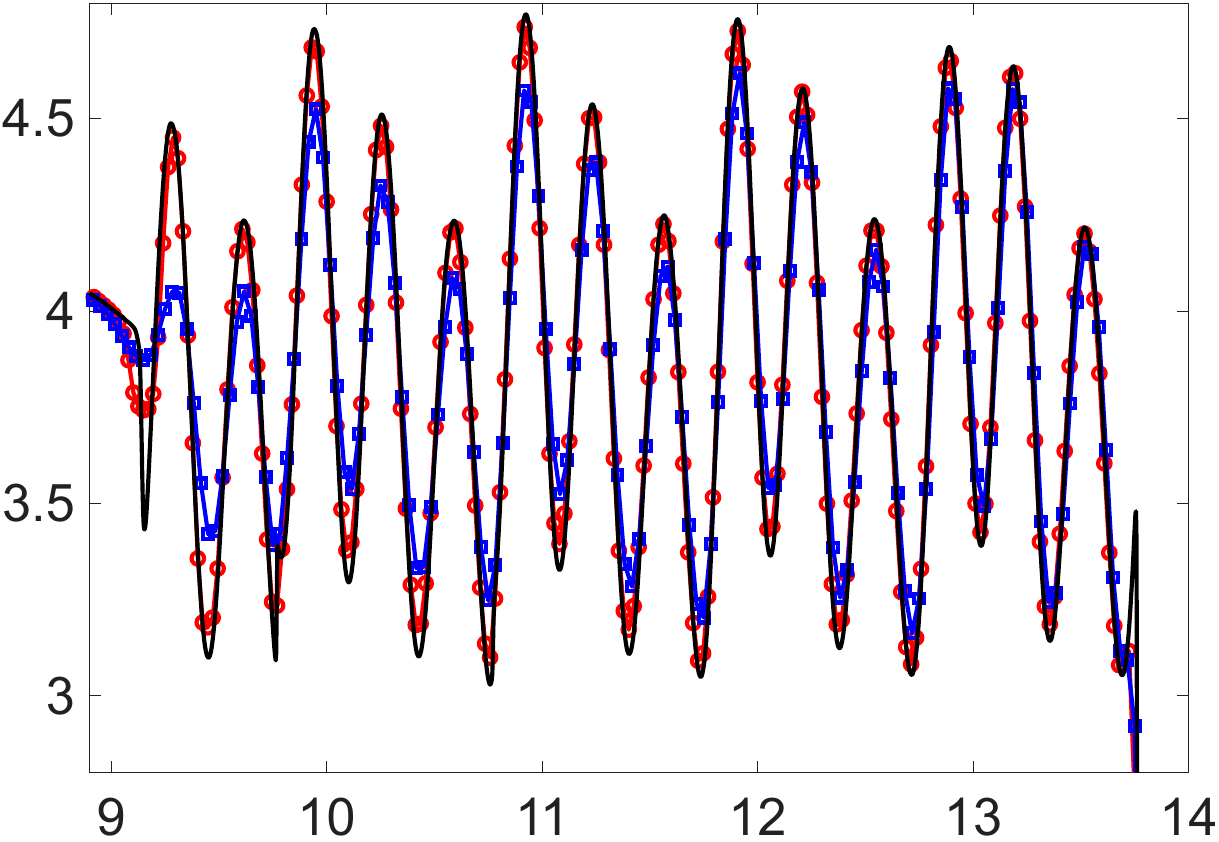}}
\caption{\sf Example 1: Density $\rho$ computed by the adaptive and A-WENO scheme with $\dx=2/87$ and $\dx=1/30$, respectively (left) and
zoom at $x\in[8.9,14]$ (right).\label{fig52}}
\end{figure}

We compute the numerical solutions by the adaptive (with the adaption coefficients $\kappa_{\rho u}=10^{-3}$ and $\kappa_p=10^{-5}$) and
A-WENO schemes until the final time $t=5$ on a uniform mesh with $\dx=1/30$. The numerical results are shown in Figure \ref{fig51} along
with the reference solution obtained by the A-WENO scheme on a much finer uniform mesh with $\dx=1/400$. As can be observed, the adaptive
scheme produces results that are only slightly better than those obtained by the A-WENO scheme. However, the CPU time consumed by the
adaptive scheme is substantially smaller than that consumed by the A-WENO scheme, which performs Ai-WENO-Z interpolations throughout the
entire computational domain. Therefore, to assess the efficiency of the proposed adaptive strategy, we rerun the adaptive scheme on a finer
mesh with $\dx=2/87$, while keeping the CPU time identical to that of the A-WENO scheme on the mesh with $\dx=1/30$. The obtained results
are shown in Figure \ref{fig52}, where one can see that for the same computational cost, the adaptive scheme clearly outperforms the A-WENO
scheme.

\smallskip
\paragraph{Example 2 (Shock-Entropy Wave Interaction Problem)}
In this example, we consider the shock--entropy wave interaction problem taken from \cite{TITAREV2004238}. The initial conditions,
\begin{equation*}
(\rho(x,0),u(x,0),p(x,0))=\begin{cases}(1.51695,0.523346,1.805),&x<-4.5,\\(1+0.1\sin(20x),0,1),&x>-4.5,\end{cases}
\end{equation*}
are prescribed in the computational domain $[-5,5]$ subject to the free boundary conditions. This test illustrates the interaction between a
Mach 1.1 forward-facing shock and small-scale density fluctuations. As the shock propagates downstream, the perturbations are convected and
amplified ahead of it, leading to the formation of high-frequency entropy waves.

\begin{figure}[ht!]
\centerline{\includegraphics[trim=0.0cm 0.0cm 0.1cm 0.0cm, clip, width=0.35\textwidth]{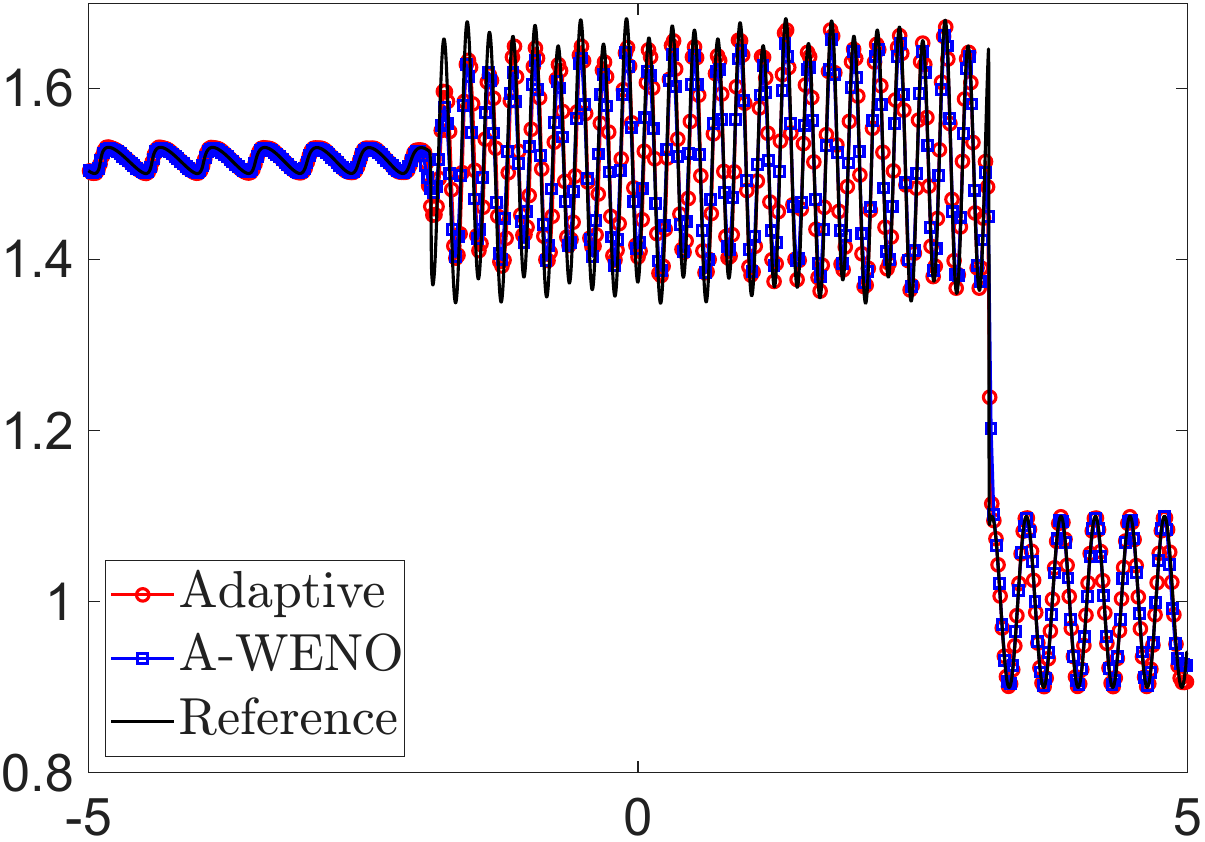}\hspace*{0.8cm}
	    \includegraphics[trim=0.0cm 0.0cm 0.1cm 0.0cm, clip, width=0.36\textwidth]{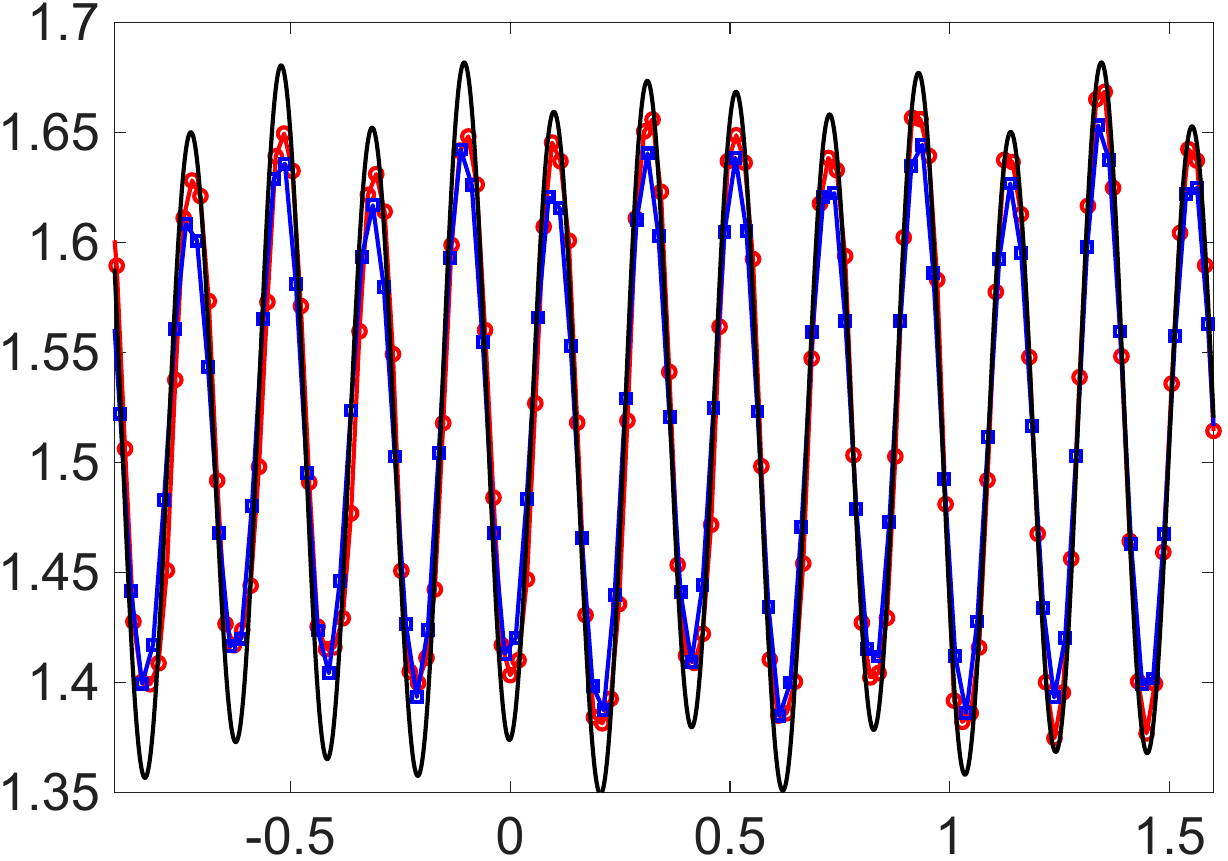}}
\caption{\sf Example 2: Density $\rho$ computed by the adaptive and A-WENO scheme with $\dx=2/105$ and $\dx=1/40$, respectively (left) and
zoom at $x\in[-0.9,1.6]$ (right).\label{fig53}}
\end{figure}

We compute the numerical solutions by the adaptive (with the adaption coefficients $\kappa_{\rho u}=5\cdot10^{-3}$ and $\kappa_p=10^{-3}$)
and A-WENO schemes until the final time $t=5$. As in the previous example, we have observed that the adaptive and A-WENO schemes produce
very similar computed solutions when implemented on the same meshes. However, in order to perform a fair comparison between the investigated
methods, we conduct the simulations on meshes, which would require the same CPU times for both methods. The obtained numerical results are
reported in Figure \ref{fig53} along with a reference solution computed using the A-WENO scheme on a much finer uniform mesh with
$\dx=1/800$. In particular, the adaptive scheme is applied on a uniform mesh with $\dx=2/105$, whereas the A-WENO scheme is computed on a
coarser mesh with $\dx=1/40$. As illustrated in the right panel, which zooms in on the interval $[-0.9,1.6]$ containing high-frequency wave
structures, the adaptive scheme achieves noticeably higher resolution compared to the A-WENO one.  

\smallskip
\paragraph{Example 3 (Blast Wave Problem)}
In the final 1-D example, we consider the blast wave problem proposed in \cite{Woodward1984TheNS}. The initial conditions,
\begin{equation*}
(\rho(x,0),u(x,0),p(x,0))= \begin{cases}(1,0,1000),&x<0.1,\\(1,0,0.01),& 0.1<x<0.9,\\(1,0,100),&x>0.9,\end{cases}
\end{equation*}
are imposed in the computational domain $[0,1]$ subject to the solid wall boundary conditions.

\begin{figure}[ht!]
\centerline{\includegraphics[trim=0.1cm 0.1cm 0.1cm 0.0cm, clip, width=0.35\textwidth]{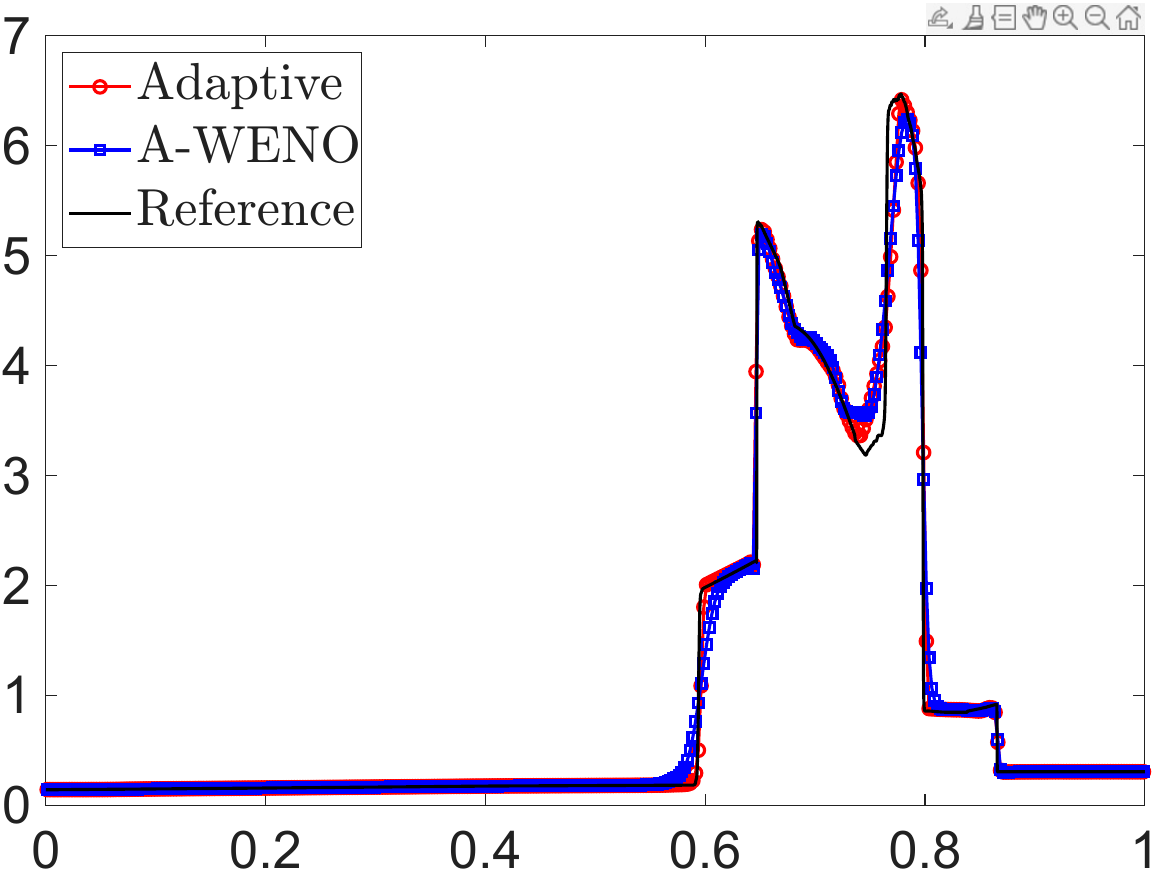}\hspace*{0.8cm}
            \includegraphics[trim=0.1cm 0.1cm 0.1cm 0.0cm, clip, width=0.37\textwidth]{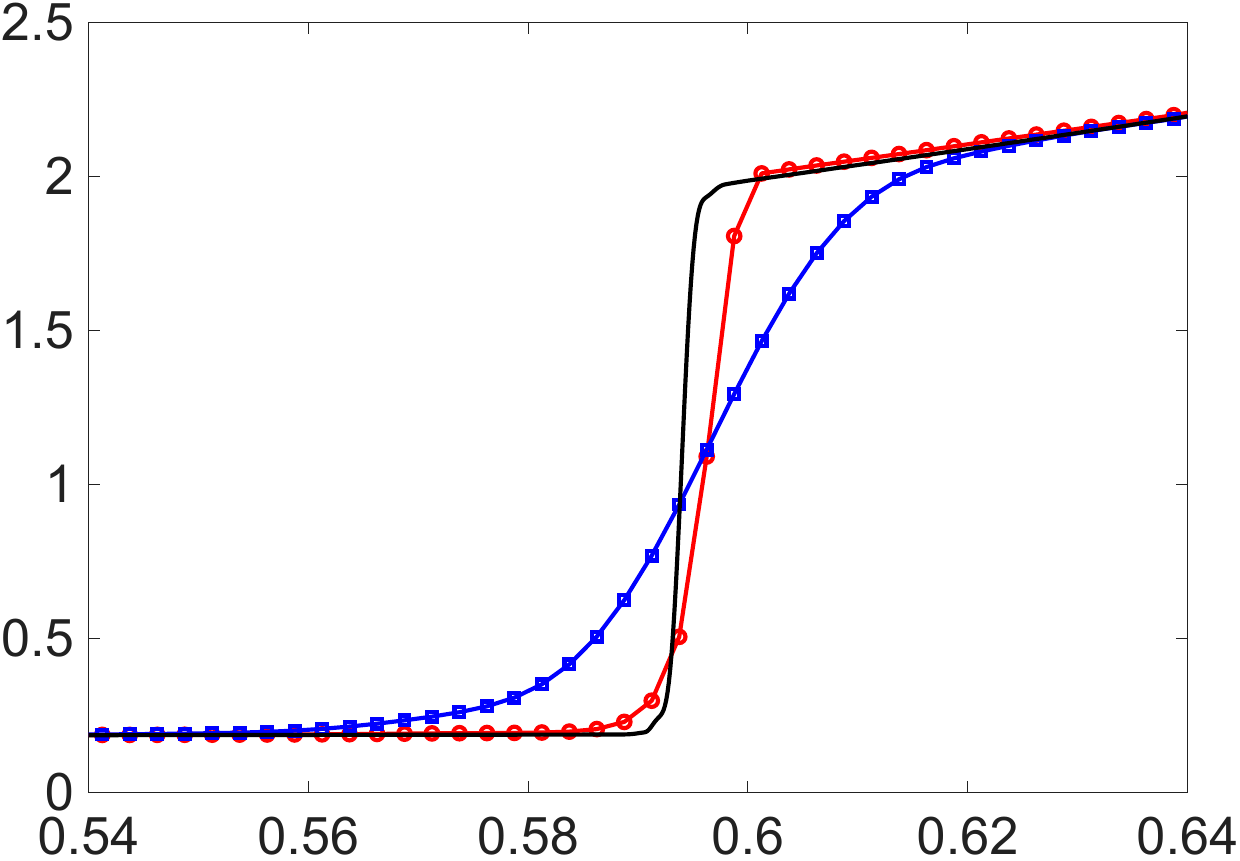}}
\caption{\sf Example 3: Density $\rho$ (left) and zoom at $x\in[0.56,0.62]$ (right).\label{fig54}}
\end{figure}

We compute the numerical solutions by the adaptive (with the adaption coefficients $\kappa_{\rho u}=10^{-4}$ and $\kappa_p=5\cdot10^{-2}$)
and A-WENO schemes until the final time $t=0.038$ on a uniform mesh with $\dx=1/400$. The obtained results are shown in Figure \ref{fig54}
along with the reference solution computed by the A-WENO scheme on a much finer mesh with $\dx=1/8000$. As one can see, the adaptive scheme
resolves the contact wave near $x\in[0.59,0.60]$ significantly sharper than the A-WENO scheme. This improvement is attributed to SI, which
correctly identifies Region RC. Consequently, the overcompressive SBM reconstruction together with the LDCU numerical fluxes is activated
near the contact waves, leading to their enhanced resolution.

\subsection{Two-Dimensional Examples}\label{sec52}
The 2-D numerical examples include two 2-D Riemann problems with strong discontinuities as well as the implosion problem and the
Rayleigh-Taylor (RT) instability problem, both featuring strong shocks and complex wave interactions.

In all of the 2-D examples presented below, we report plots of the different regions of the computational domain identified by the SIs at
the final time. Specifically, the areas classified as Regions S, RC, and RNC based on the SIs in the $x$- and $y$-directions are displayed
using three different colors: blue, green, and red, respectively.

\smallskip
\paragraph{Example 4 (Riemann Problems)}
In this example, we consider two Riemann problems from \cite{KTrp} with the following two sets of initial data:
\begin{equation*}
\mbox{Configuration 3:}\quad(\rho,u,v,p)\left|_{(x,y,0)}\right.=\begin{cases}(1.5,0,0,1.5),&x>1,~y>1,\\(0.5323,1.206,0,0.3),&x<1,~y>1,\\
(0.138,1.206,1.206,0.029),&x<1,~y<1,\\(0.5323,0,1.206,0.3),&x>1,~y<1,\end{cases}
\end{equation*}
prescribed in the computational domain $[0,1.2]\times[0,1.2]$;
\begin{equation*}
\mbox{Configuration 12:}\quad(\rho,u,v,p)\left|_{(x,y,0)}\right.=\begin{cases}(0.5313,0,0,0.4),&x>0.5,~y>0.5,\\
(1,0.7276,0,1),&x<0.5,~y>0.5,\\(0.8,0,0,1),&x<0.5,~y<0.5,\\(1,0,0.7276,1),&x>0.5,~y<0.5,\end{cases}
\end{equation*}
prescribed in the computational domain $[0,0.6]\times[0,0.6]$. For both configurations, the free boundary conditions are imposed.

\begin{figure}[ht!]
\centerline{\includegraphics[trim=1.0cm 2.8cm 0.8cm 2.3cm, clip, width=0.35\textwidth]{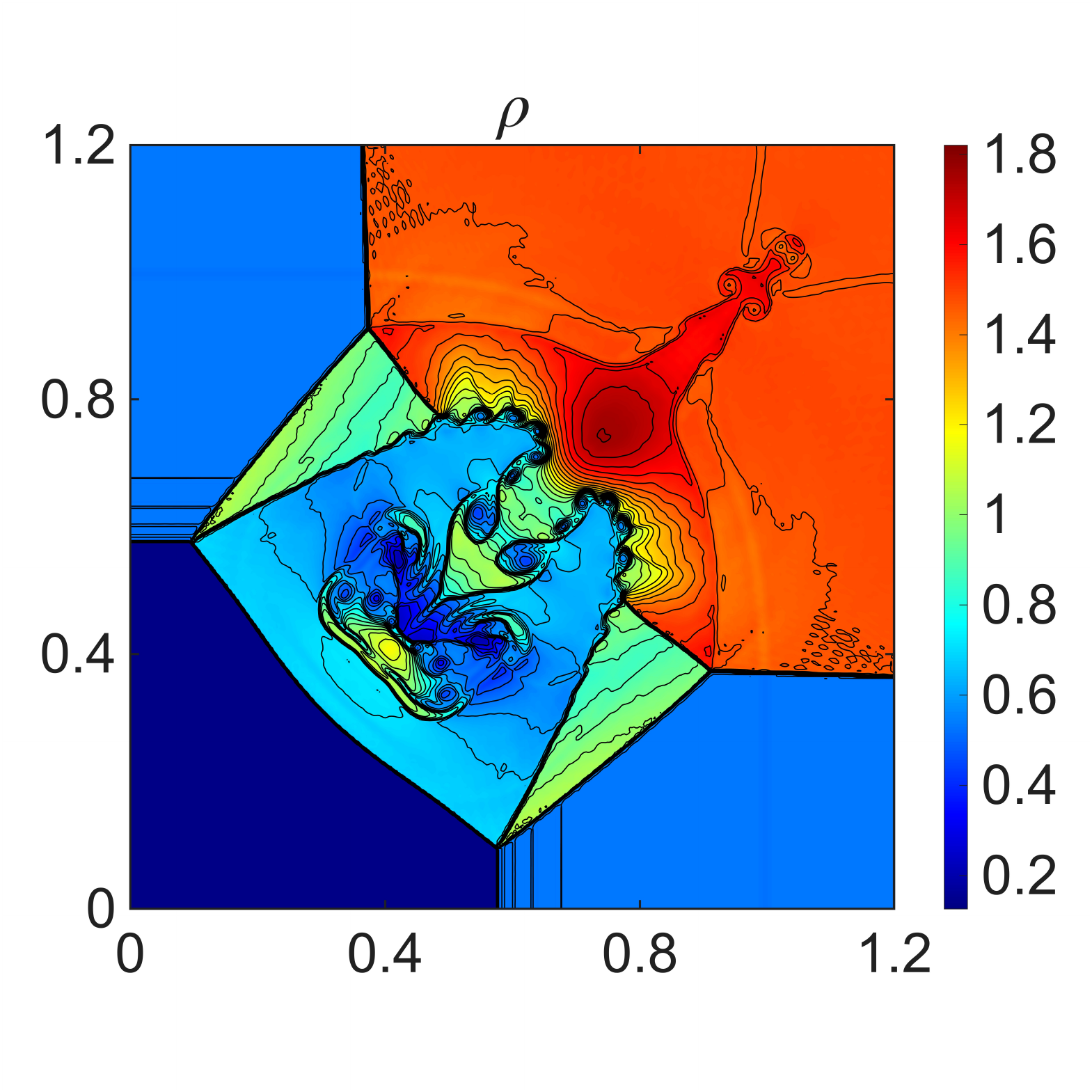}
\hspace*{0.5cm}
            \includegraphics[trim=1.0cm 2.8cm 0.8cm 2.3cm, clip, width=0.35\textwidth]{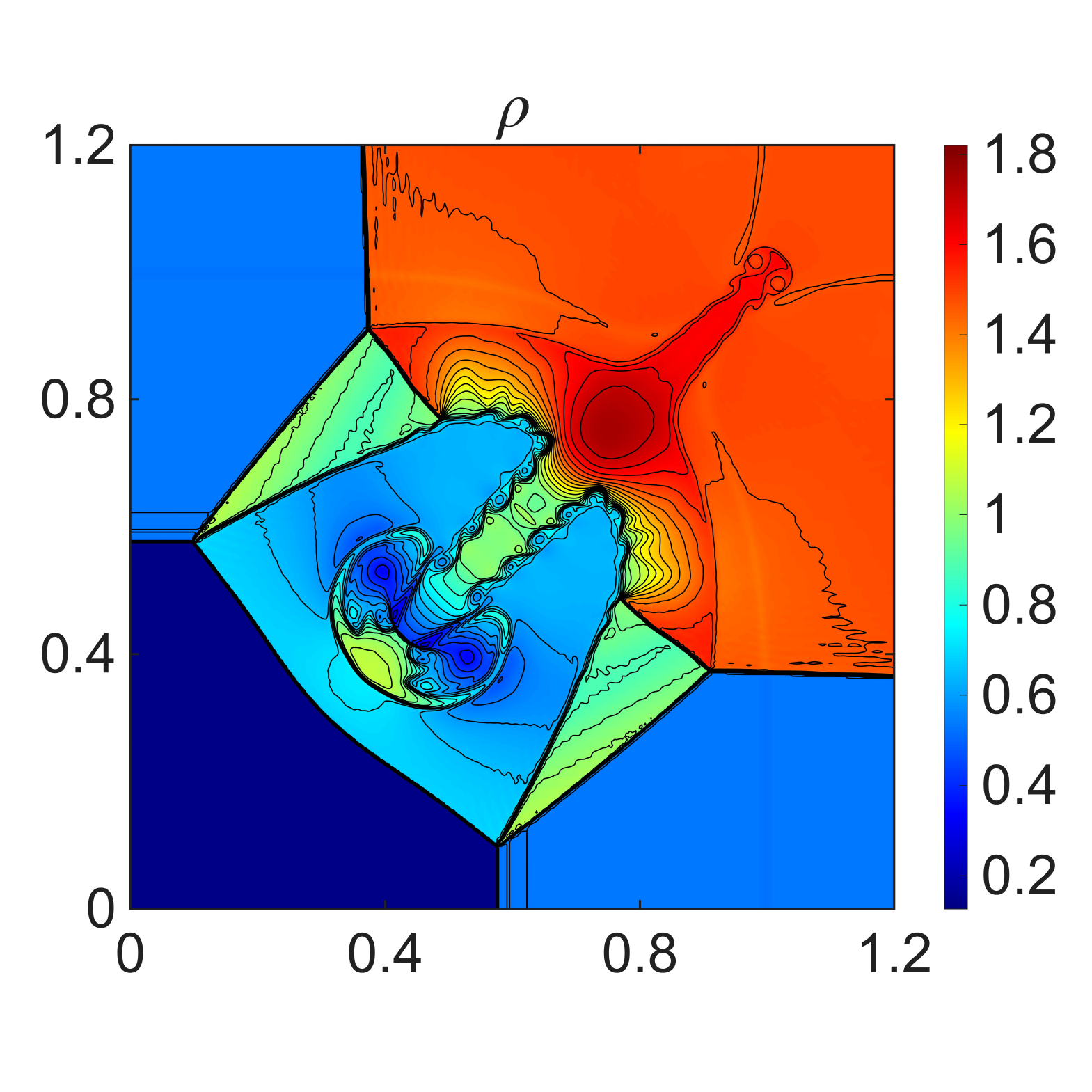}}
\vskip10pt
\hspace*{-0.2cm}\centerline{\includegraphics[trim=1.3cm 1.9cm 1.4cm 1.0cm, clip, width=0.31\textwidth]{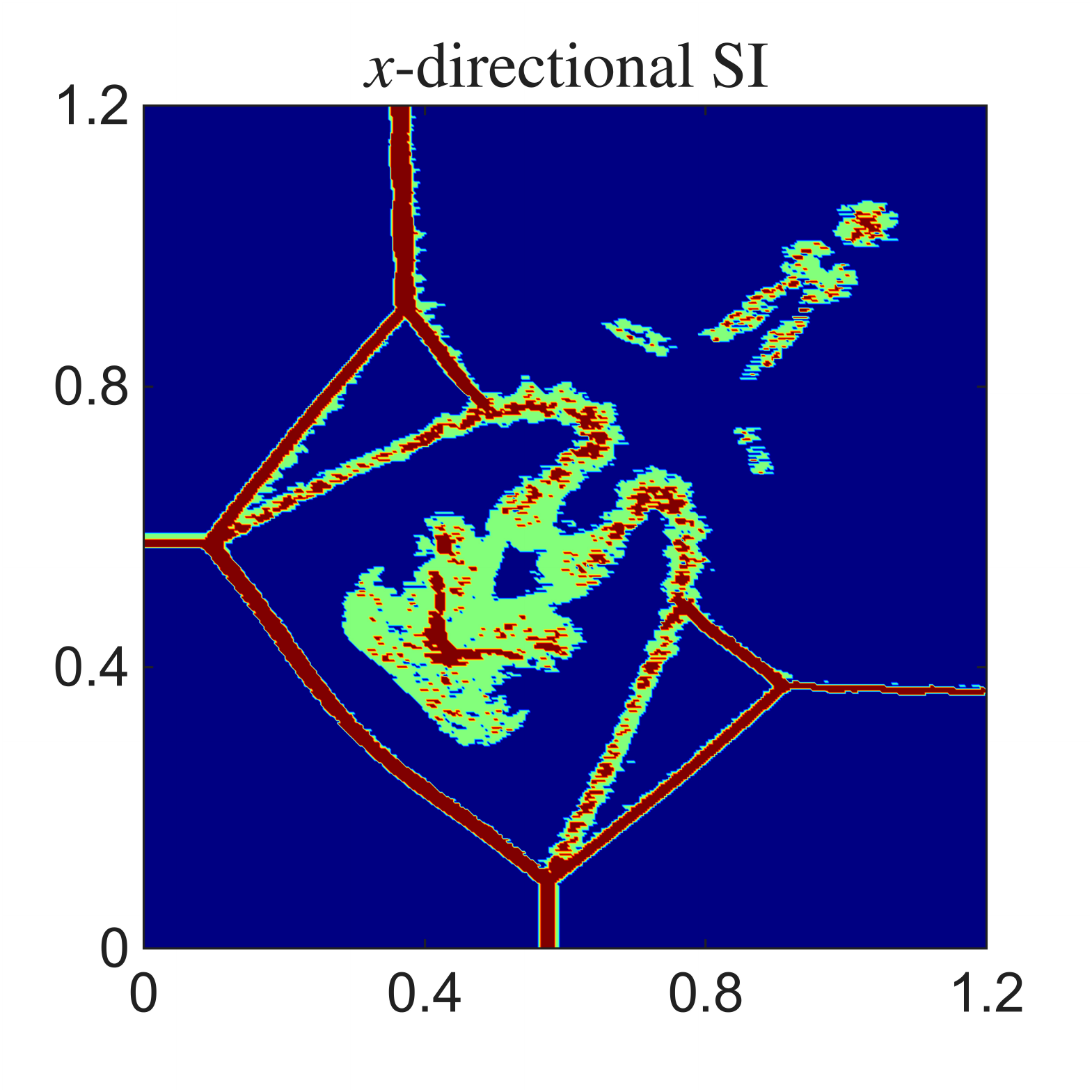}
\hspace*{1.0cm}
	    \includegraphics[trim=1.0cm 1.9cm 1.4cm 1.0cm, clip, width=0.31\textwidth]{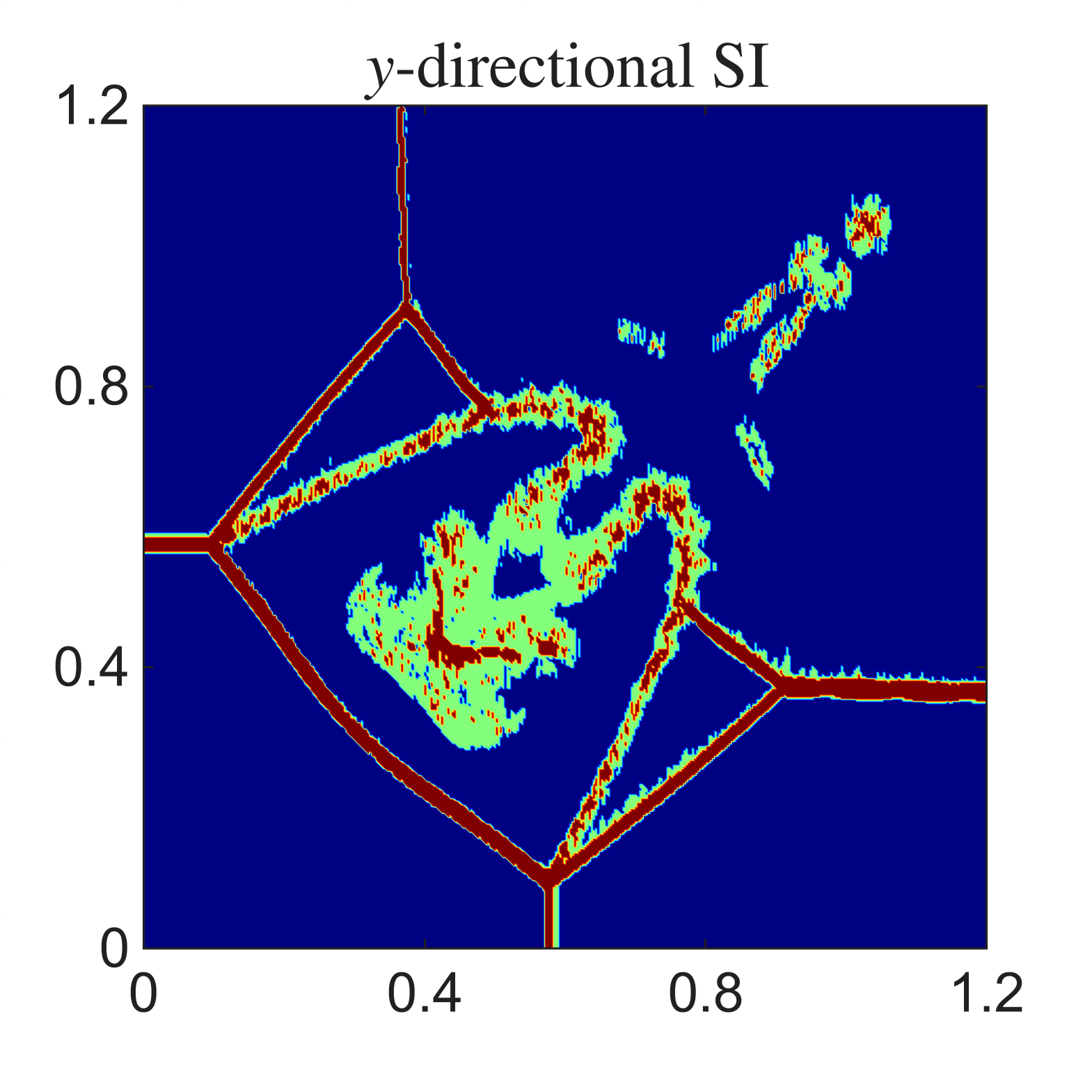}}
\caption{\sf Example 4 (Configuration 3): Density $\rho$ computed by the adaptive (top left) and A-WENO (top right) schemes on a uniform
mesh with $\dx=\dy=3/1000$, along with SIs in the $x$-direction (bottom left) and $y$-direction (bottom right).\label{fig55}}
\end{figure}

\begin{figure}[ht!]
\centerline{\includegraphics[trim=1.0cm 2.8cm 0.8cm 2.3cm, clip, width=0.35\textwidth]{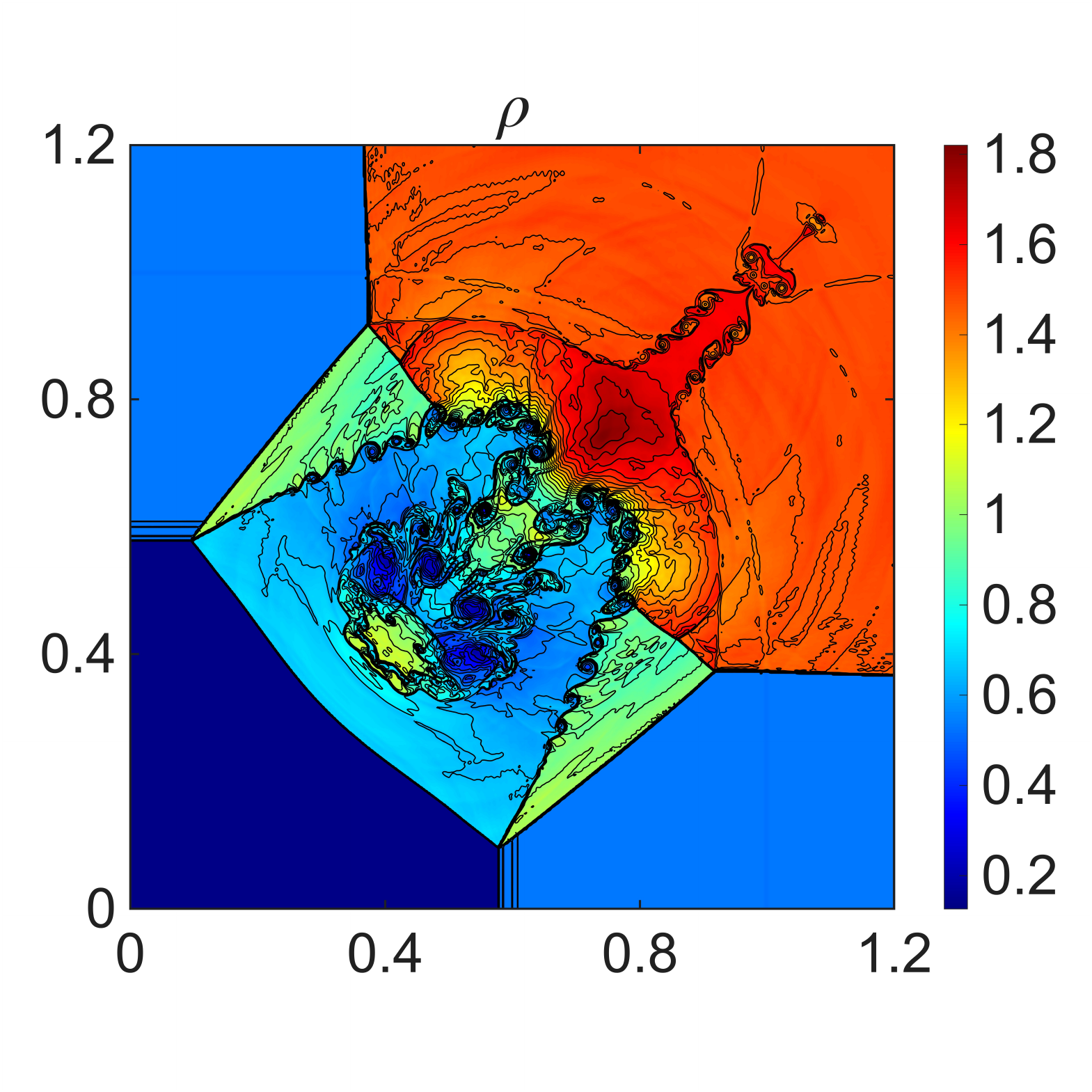}
\hspace*{0.5cm}
            \includegraphics[trim=1.0cm 2.8cm 0.8cm 2.3cm, clip, width=0.35\textwidth]{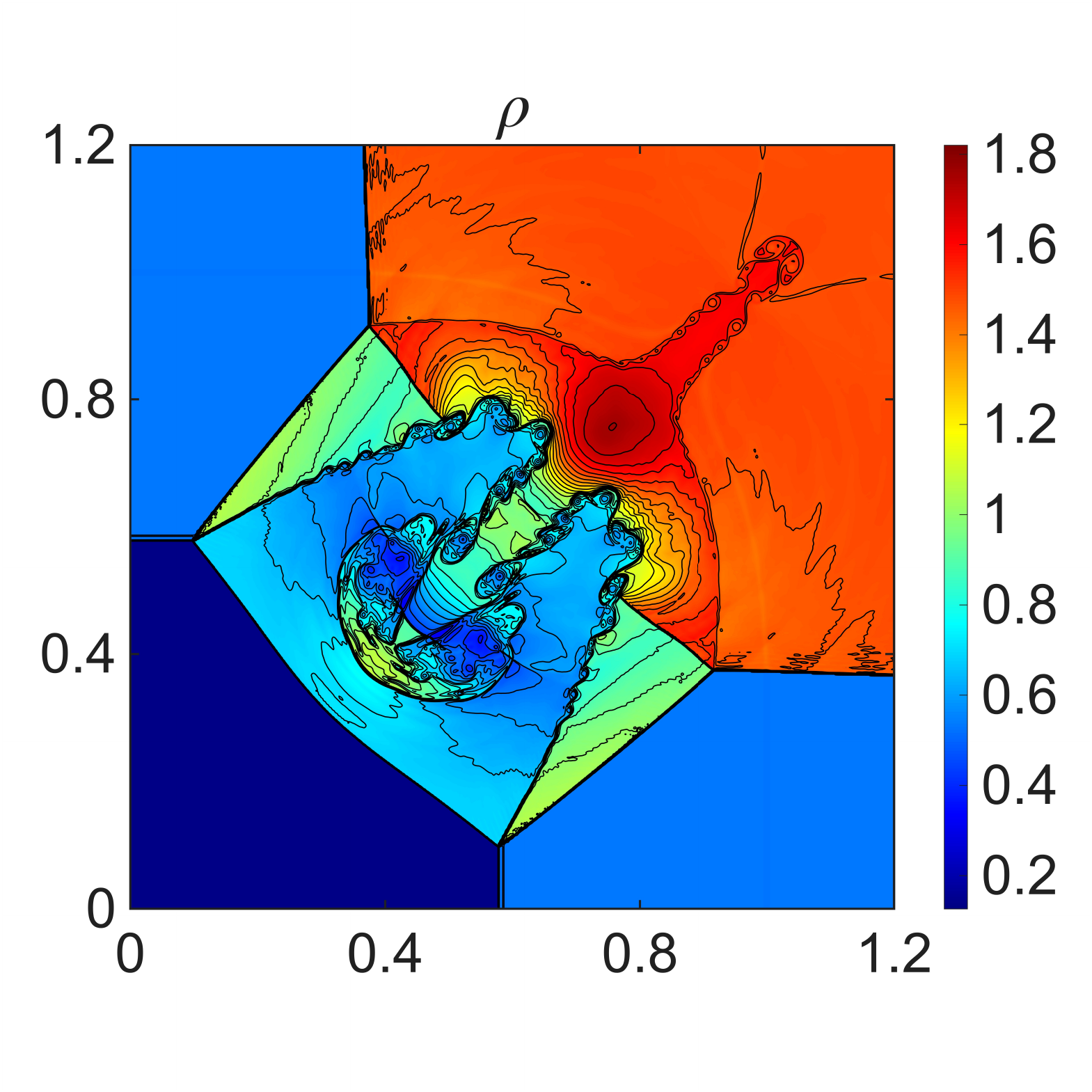}}
\vskip10pt
\hspace*{-0.2cm}\centerline{\includegraphics[trim=1.0cm 1.9cm 1.4cm 1.0cm, clip, width=0.31\textwidth]{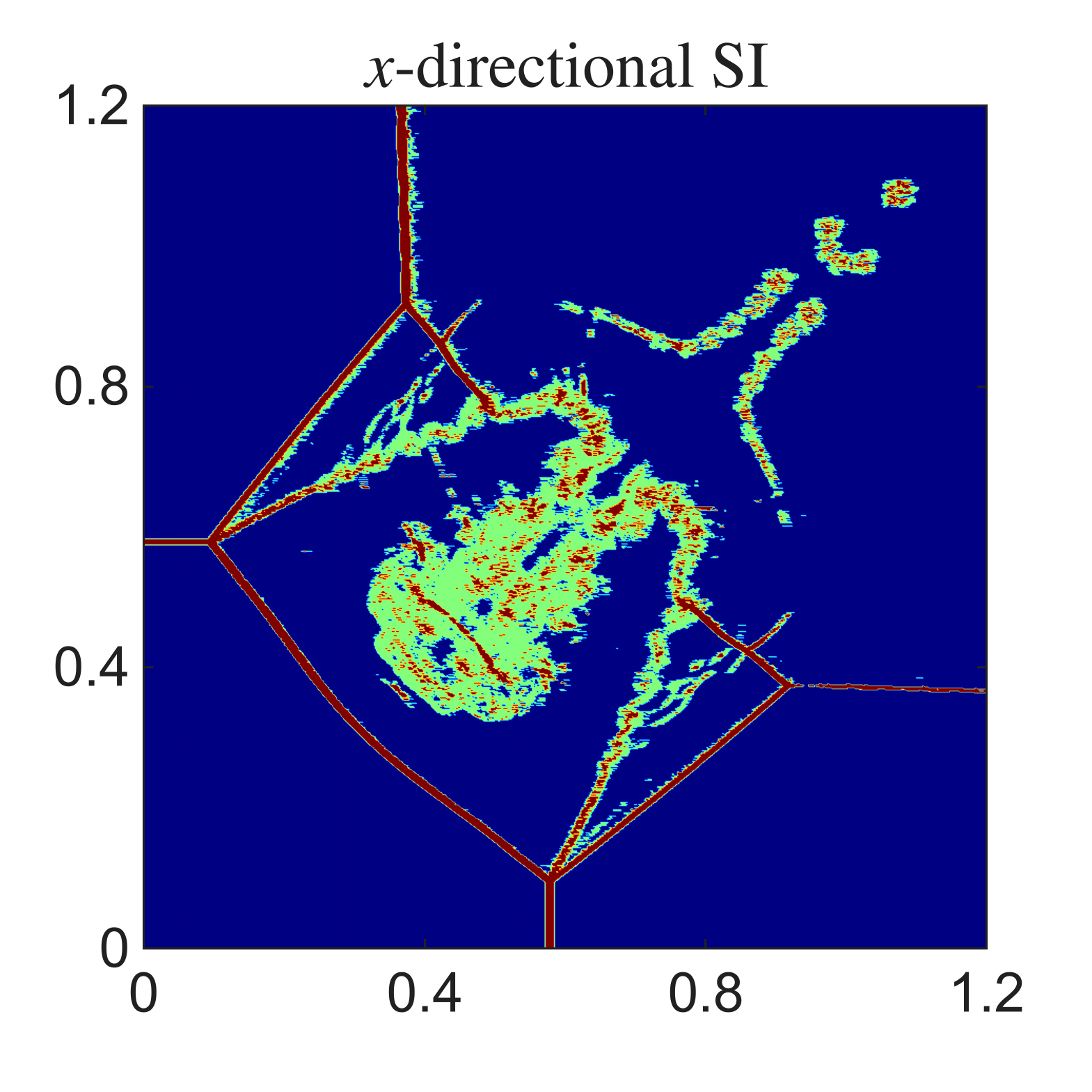}
\hspace*{1.0cm}
            \includegraphics[trim=1.0cm 1.9cm 1.4cm 1.0cm, clip, width=0.31\textwidth]{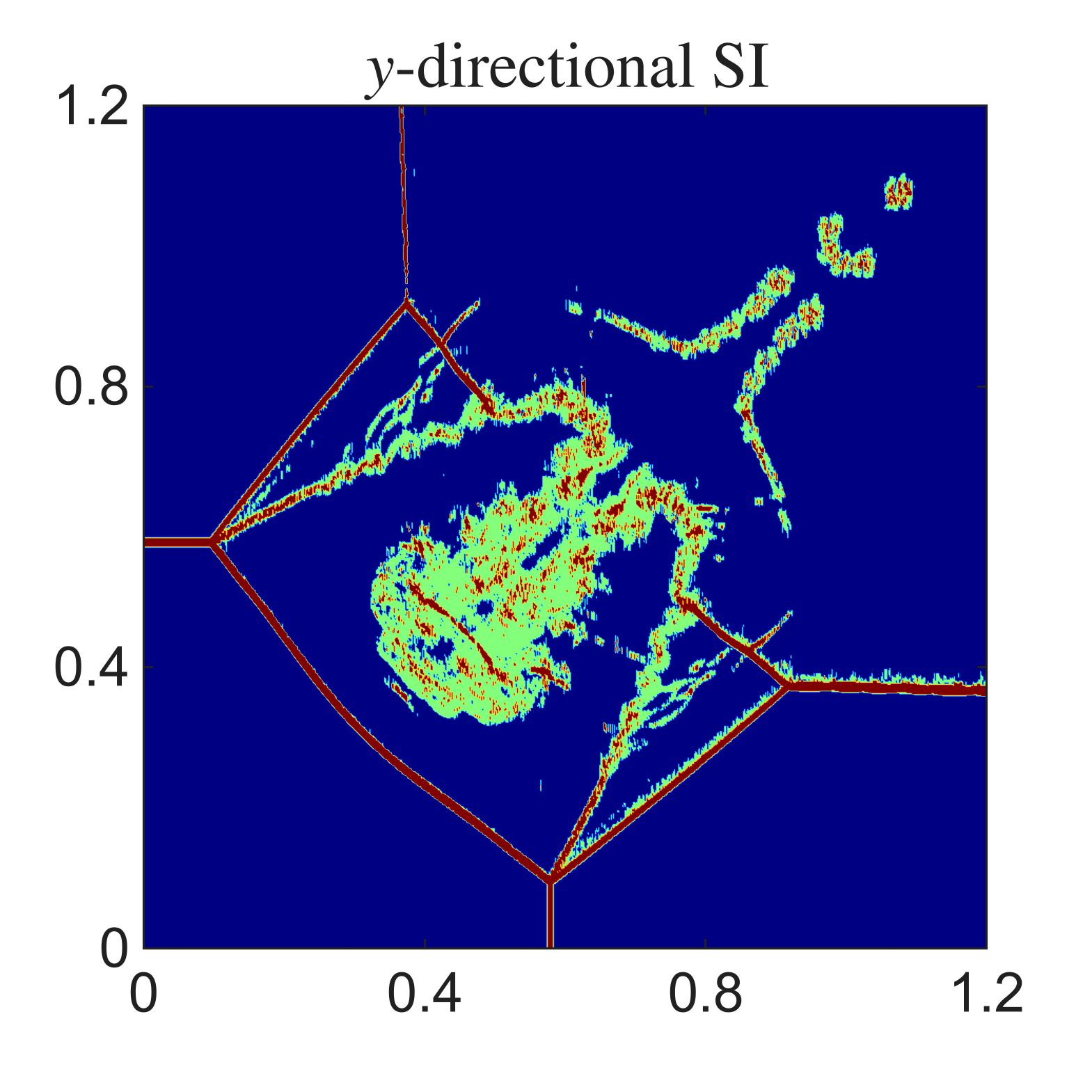}}
\caption{\sf Example 4 (Configuration 3): Same as in Figure \ref{fig55} but for the mesh with $\dx=\dy=3/2000$.\label{fig56}}
\end{figure}

We begin with Configuration 3 and compute the numerical solutions by the adaptive (with the adaption coefficients
$\kappa_{\rho u}=\kappa_{\rho v}=10^{-2}$ and $\kappa_p=5\cdot10^{-2}$) and A-WENO schemes until the final time $t=1$ on a uniform mesh with
$\dx=\dy=3/1000$. The numerical results are presented in Figure \ref{fig55}, where one can see that both methods accurately capture the
complex interaction of four shocks and the resulting mushroom-shaped vortical structures. Notably, the proposed adaptive scheme resolves the
fine-scale Kelvin-Helmholtz instabilities with lower numerical dissipation. To demonstrate the robustness of the proposed adaptive scheme in
terms of the independence on the choice of the adaption coefficients on the mesh size, we refine the mesh to $\dx=\dy=3/2000$ and conduct
the same simulations. The obtained results are displayed in Figure \ref{fig56}. One can observe that the adaptive scheme still clearly
outperforms the A-WENO scheme in resolving complex wave structures. On both meshes, Regions S, RC, and RNC are adequately
indicated, as can be seen in the bottom rows in Figures \ref{fig55} and \ref{fig56}.

\begin{figure}[ht!]
\centerline{\includegraphics[trim=0.9cm 2.8cm 0.8cm 2.3cm, clip, width=0.34\textwidth]{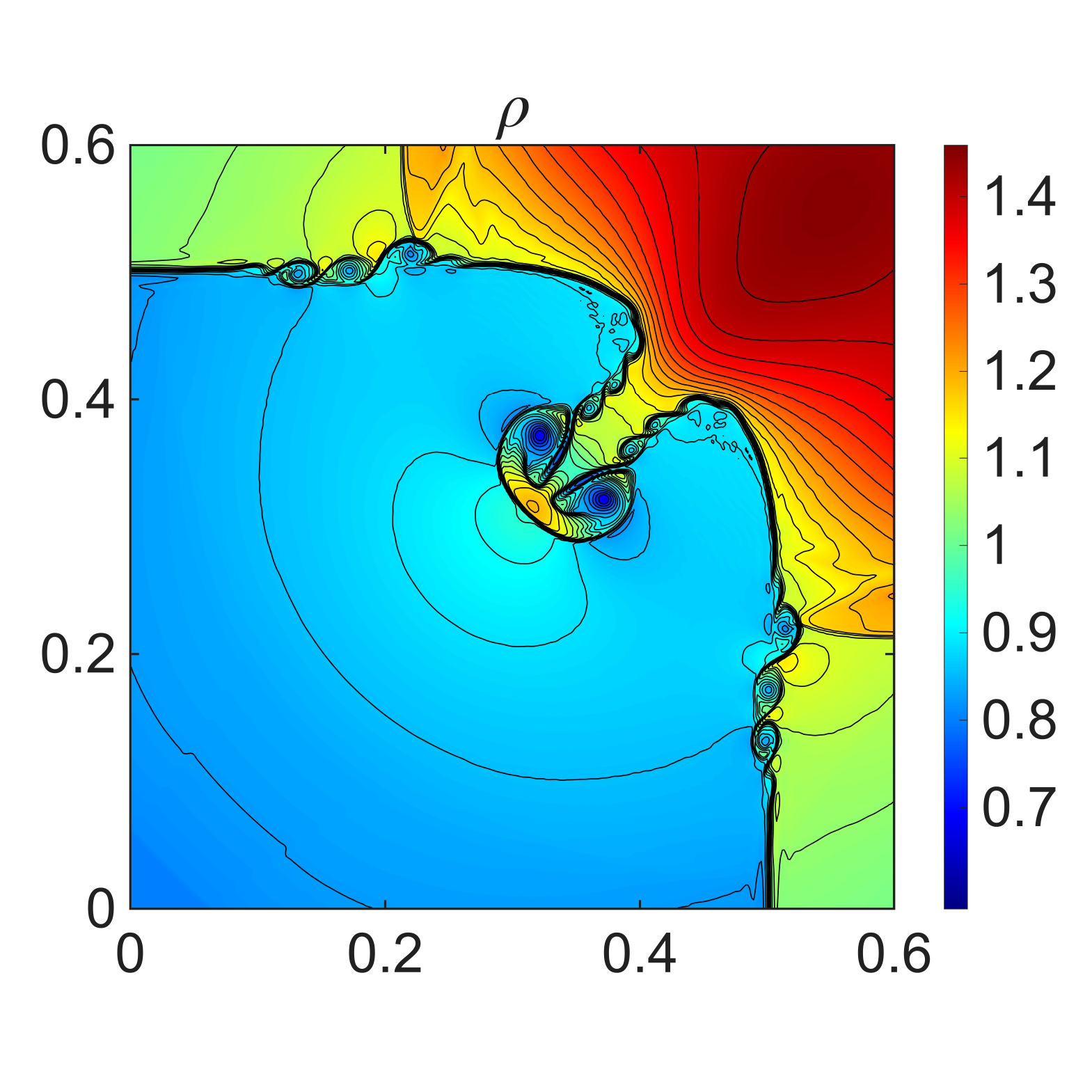}
\hspace*{0.5cm}
            \includegraphics[trim=0.9cm 2.8cm 0.8cm 2.3cm, clip, width=0.34\textwidth]{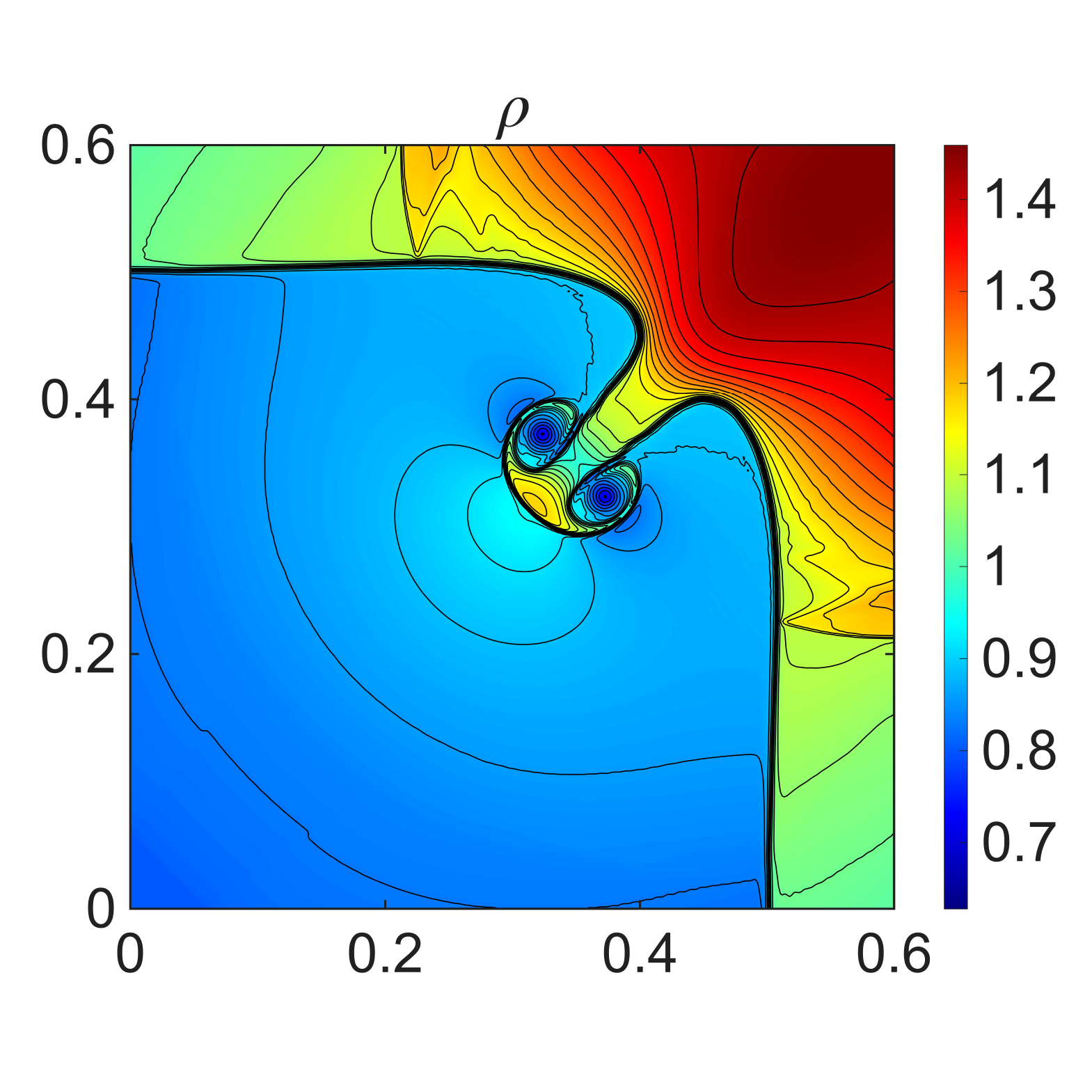}}
\vskip10pt   
\hspace*{-0.2cm}\centerline{\includegraphics[trim=0.9cm 1.9cm 1.0cm 1.0cm, clip, width=0.305\textwidth]{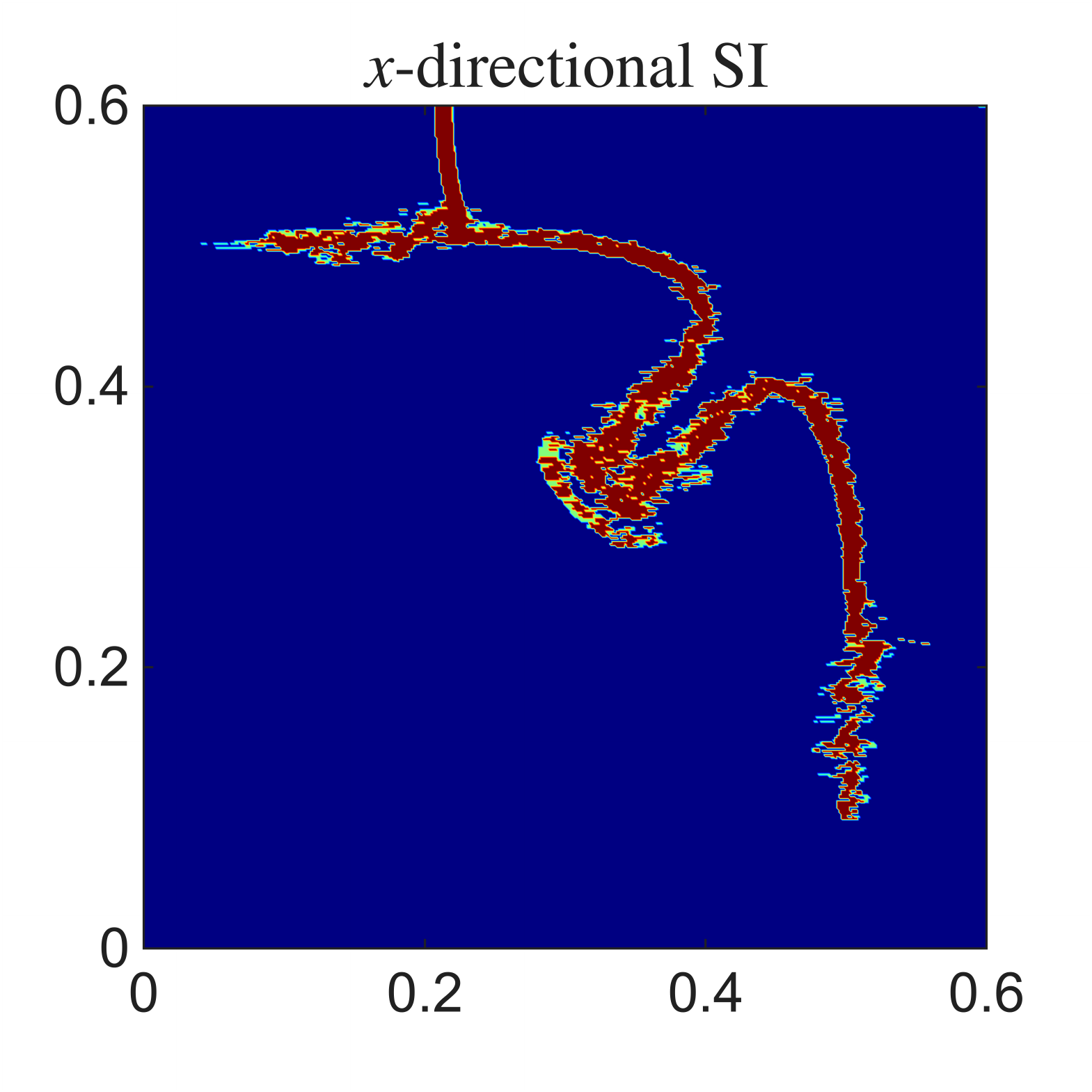}
\hspace*{1.0cm}
            \includegraphics[trim=0.9cm 1.9cm 1.0cm 1.0cm, clip, width=0.305\textwidth]{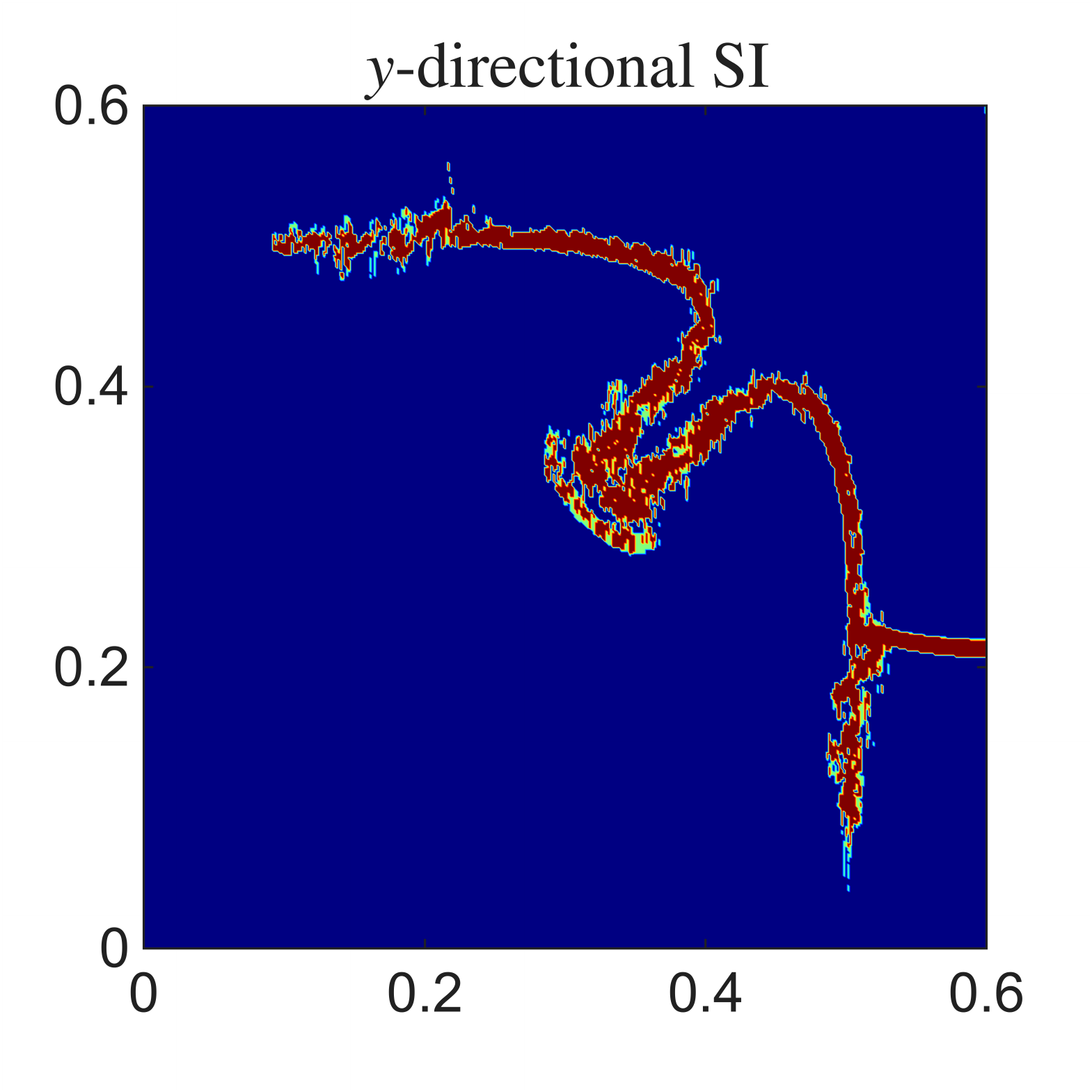}}
\caption{\sf Example 4 (Configuration 12): Density $\rho$ computed by the adaptive (top left) and A-WENO (top right) schemes on a uniform
mesh with $\dx=\dy=3/2000$, along with SIs in the $x$-direction (bottom left) and $y$-direction (bottom right).\label{fig57}}
\end{figure}

\begin{figure}[ht!]
\centerline{\includegraphics[trim=0.9cm 2.8cm 0.8cm 2.3cm, clip, width=0.34\textwidth]{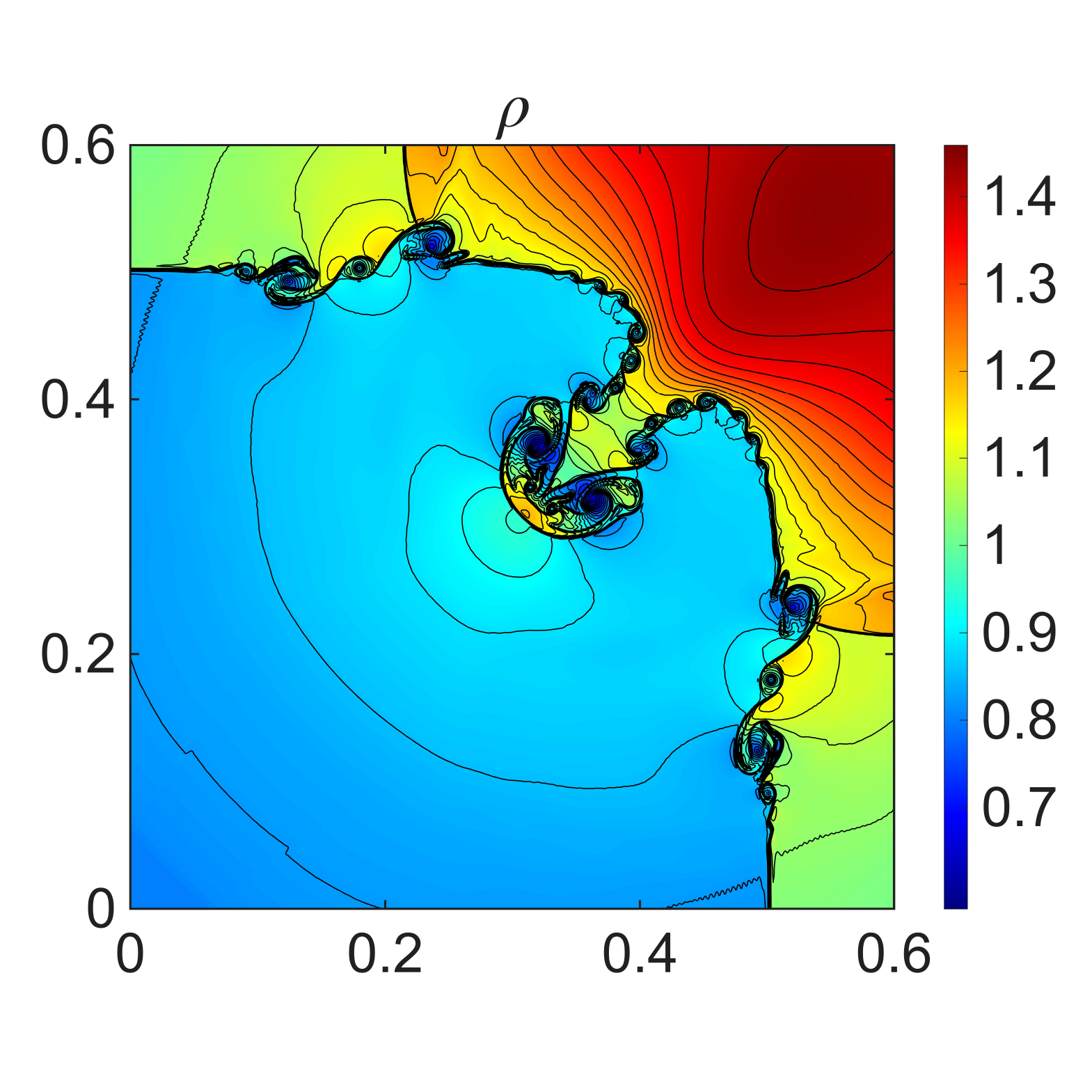}
\hspace*{0.5cm}
	    \includegraphics[trim=0.9cm 2.8cm 0.8cm 2.3cm, clip, width=0.34\textwidth]{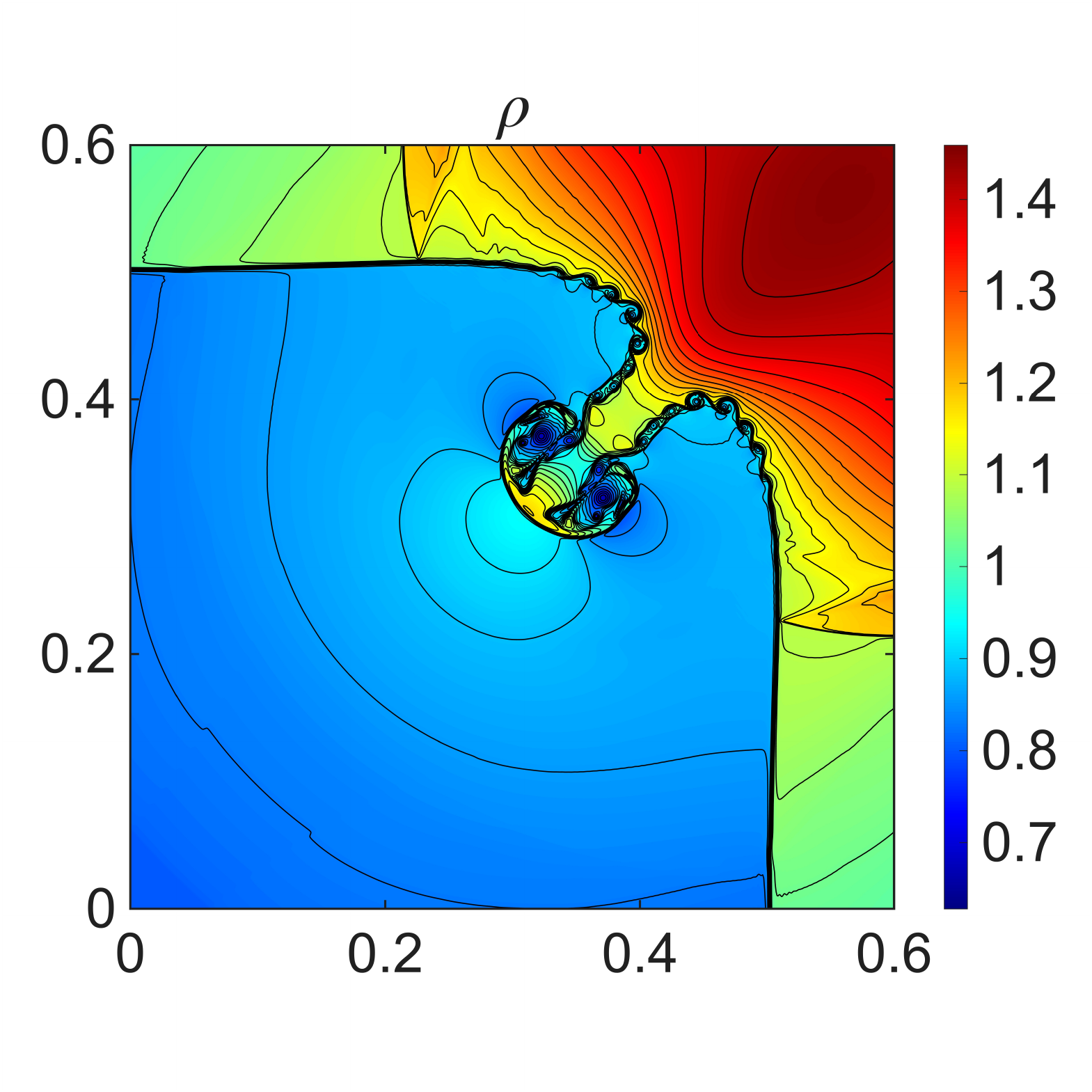}}
\vskip10pt   
\hspace*{-0.2cm}\centerline{\includegraphics[trim=0.9cm 1.9cm 1.0cm 1.0cm, clip, width=0.305\textwidth]{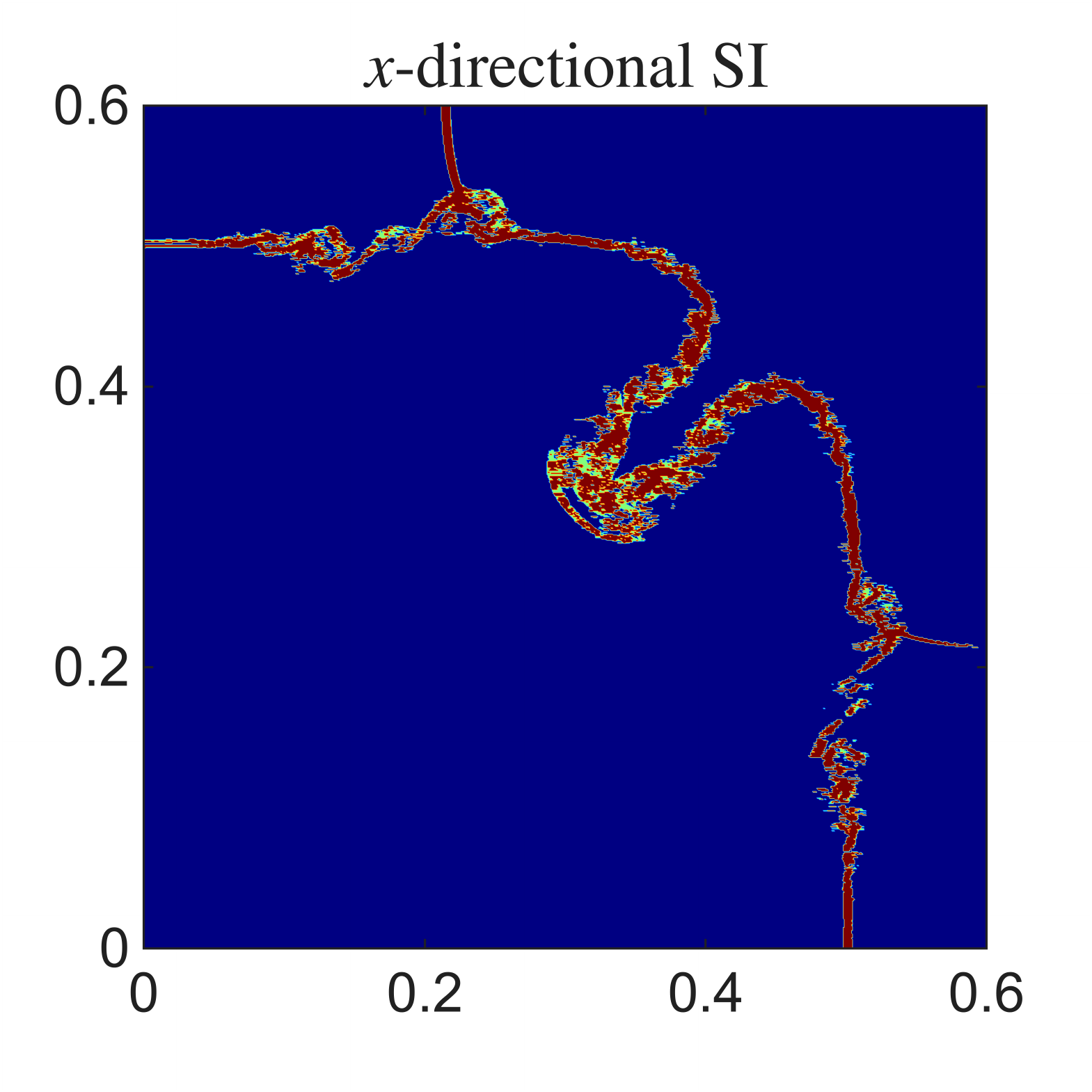}
\hspace*{1.0cm}
	    \includegraphics[trim=0.9cm 1.9cm 1.0cm 1.0cm, clip, width=0.305\textwidth]{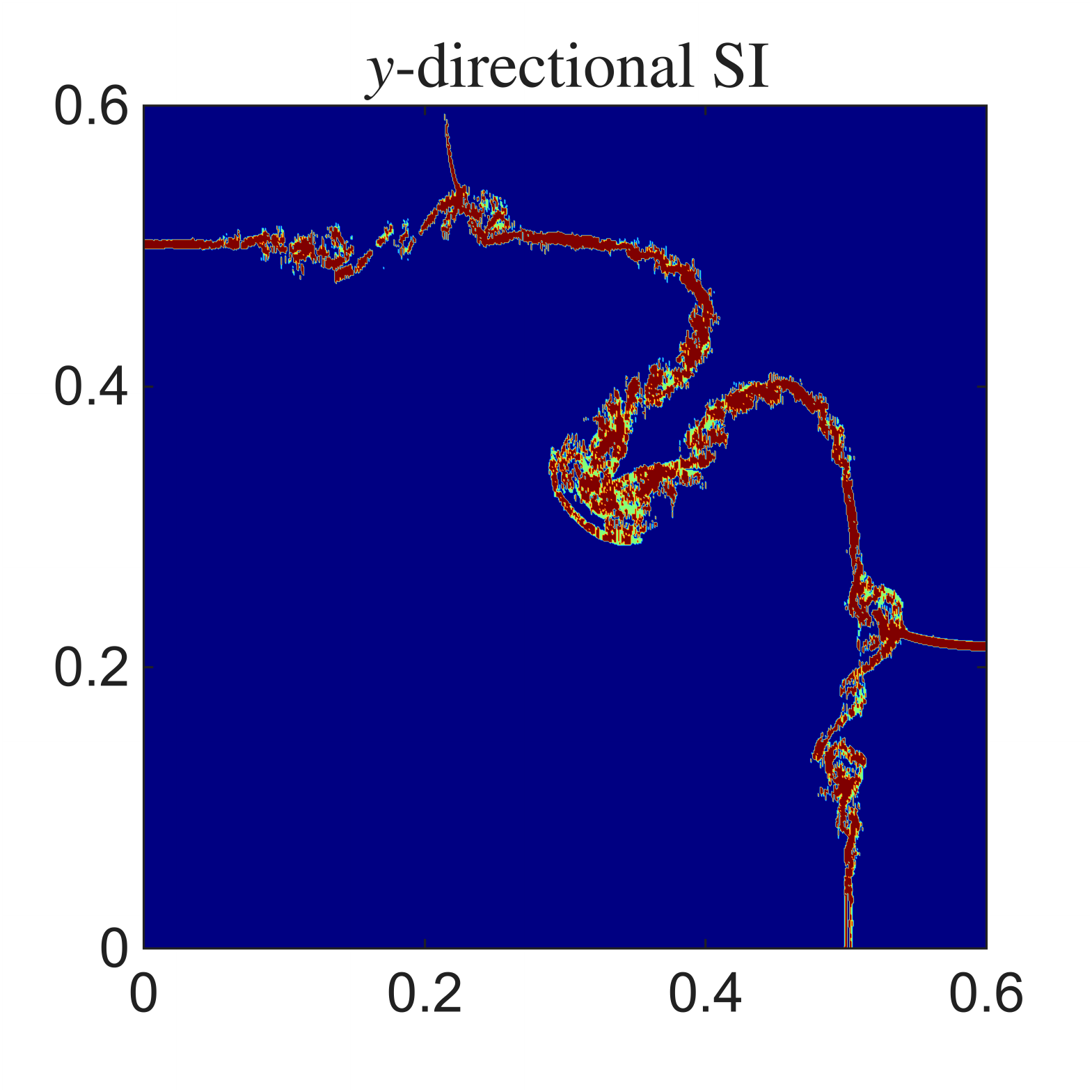}}
\caption{\sf Example 4 (Configuration 12): Same as in Figure \ref{fig57} but for the mesh with $\dx=\dy=3/4000$.\label{fig58}}
\end{figure}

Next, we compute the numerical solutions for Configuration 12 using both the adaptive (with the adaption coefficients
$\kappa_{\rho u}=\kappa_{\rho v}=0.9$ and $\kappa_p=1$) and A-WENO schemes until the final time $t=0.5$. The numerical results, obtained on
two uniform meshes with $\dx=\dy=3/2000$ and $3/4000$, are presented in Figures \ref{fig57} and \ref{fig58}, respectively. It is evident
that the adaptive scheme yields sharper resolution and finer-scale wave structures than its A-WENO counterpart, and that Regions S, RC, and
RNC are adequately indicated.

\smallskip
\paragraph{Example 5 (Implosion problem)}
In this section, we consider the implosion problem taken from \cite{liska2003comparison}. The initial conditions,
\begin{equation*}
(\rho(x,y,0),u(x,y,0),v(x,y,0),p(x,y,0))=
\left\{\begin{aligned}
&(0.125,0,0,0.14),&&x+y<0.15,\\
&(1,0,0,1),&&\mbox{otherwise},
\end{aligned}\right.
\end{equation*}
are prescribed in the computational domain $[0,0.3]\times[0,0.3]$ subject to the solid wall boundary conditions. As time evolves, the
surrounding high-density fluid converges toward the origin, compressing the low-density core and giving rise to a sequence of complex wave
interactions, including multiple reflections and the eventual formation of a prominent jet-like structure that propagates along the symmetry
line $y=x$.

\begin{figure}[ht!]
\centerline{\includegraphics[trim=1.0cm 2.9cm 0.9cm 2.3cm, clip, width=0.34\textwidth]{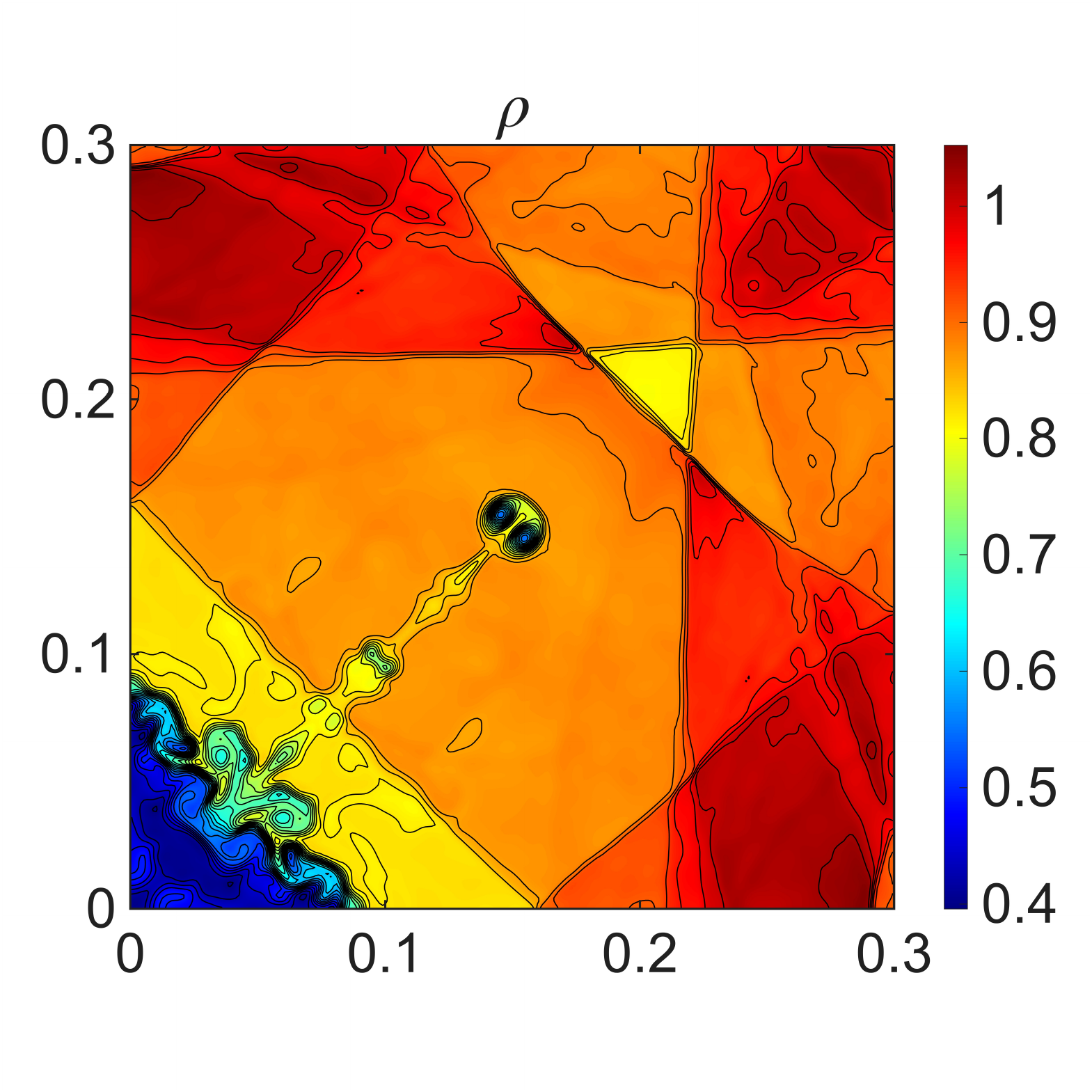}
\hspace*{0.5cm}
            \includegraphics[trim=1.0cm 2.9cm 0.9cm 2.3cm, clip, width=0.34\textwidth]{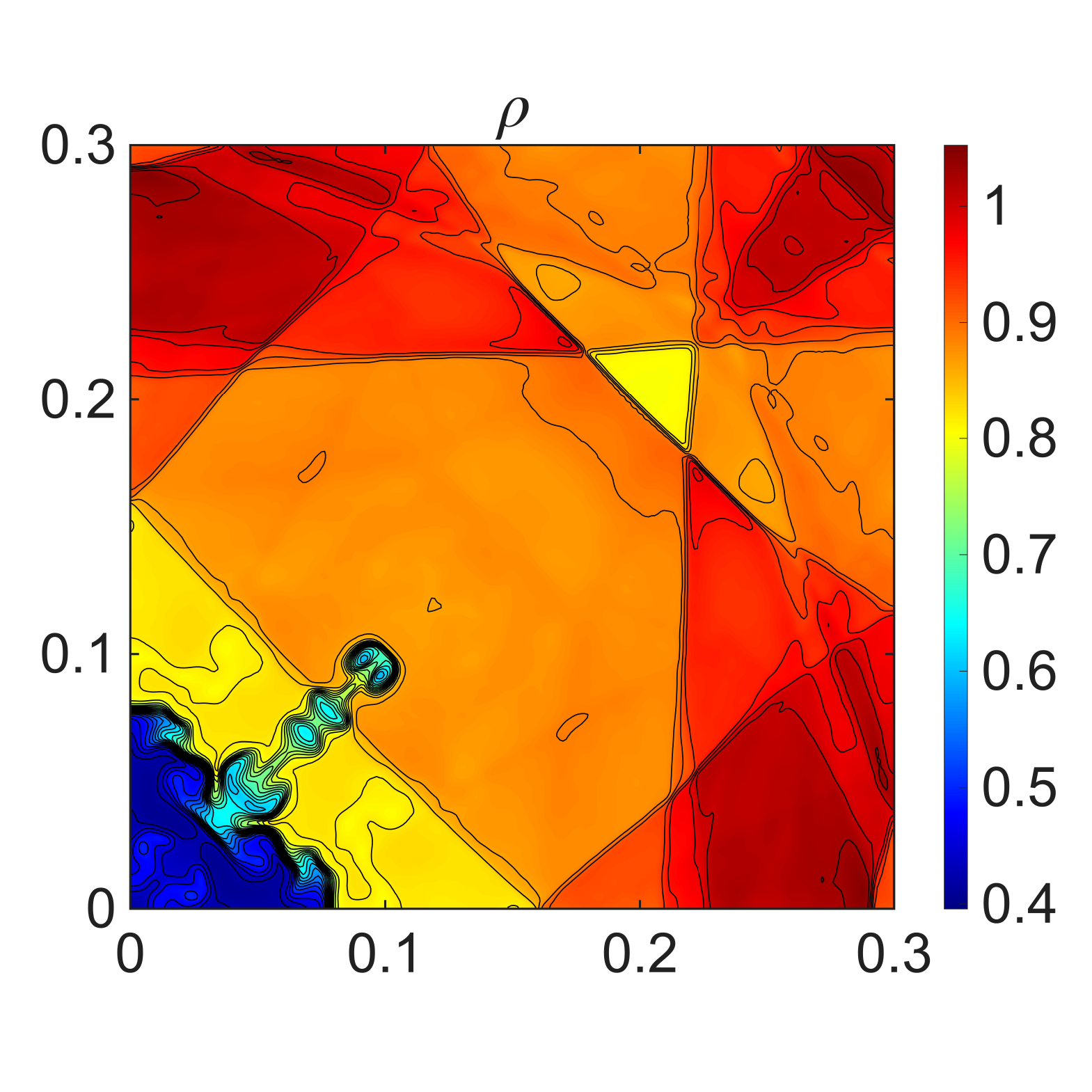}}
\vskip10pt
\hspace*{-0.3cm}\centerline{\includegraphics[trim=0.9cm 1.9cm 1.7cm 1.0cm, clip, width=0.30\textwidth]{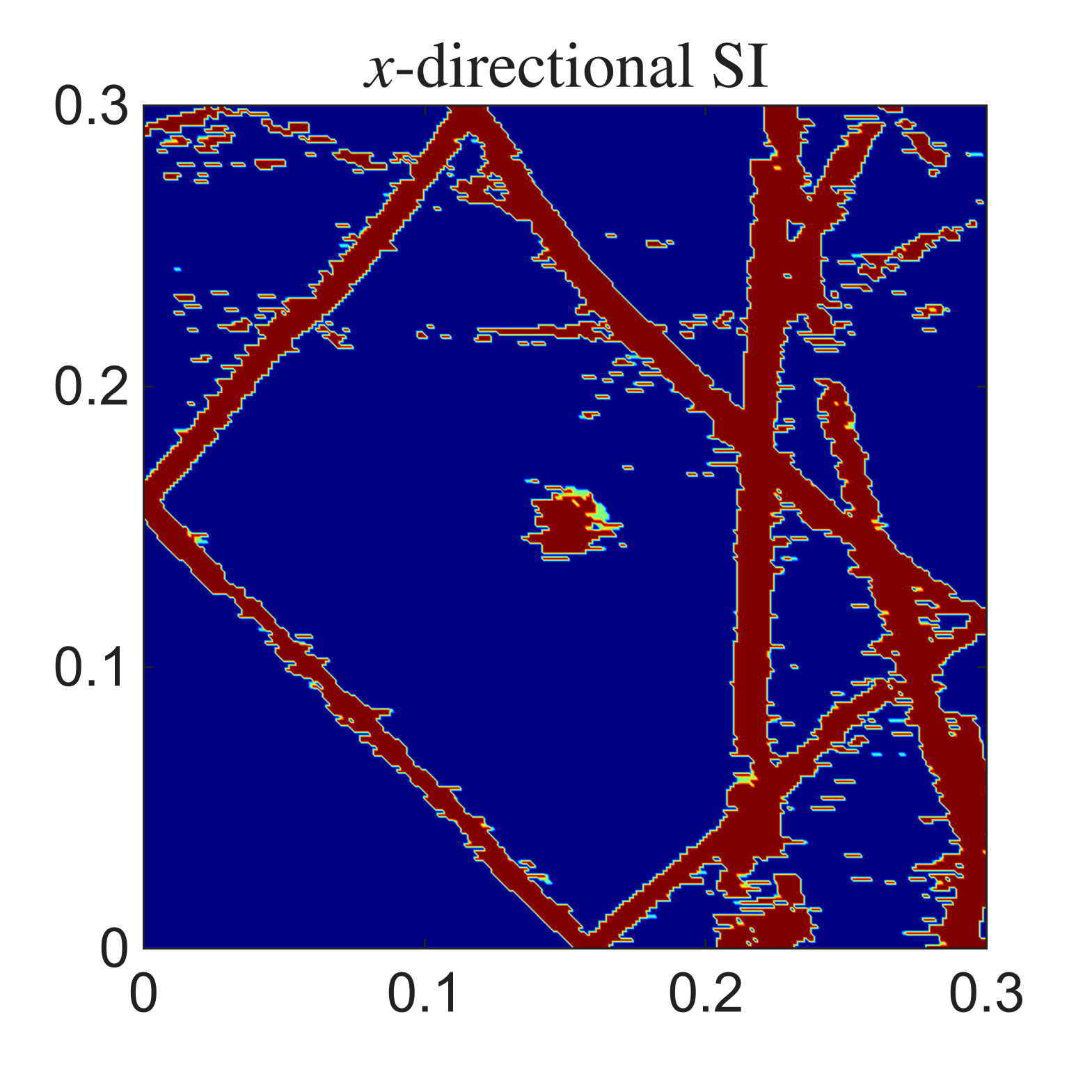}\hspace*{1.1cm}
        \includegraphics[trim=0.9cm 1.9cm 1.7cm 1.0cm, clip, width=0.30\textwidth]{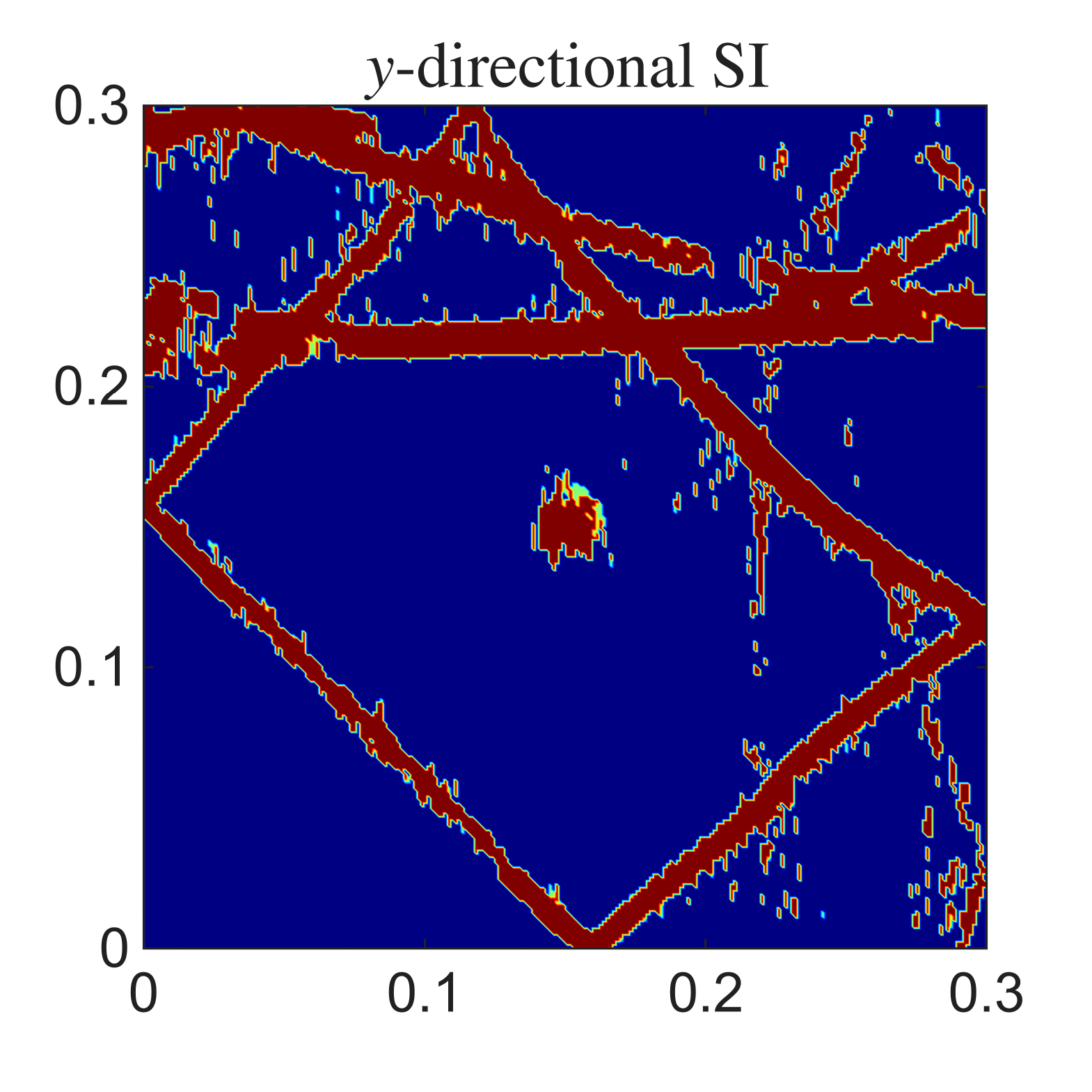}}
\caption{\sf Example 5: Density $\rho$ computed by the adaptive (top left) and A-WENO schemes (top right), along with SIs in the 
$x$-direction (bottom left) and the $y$-direction (bottom right). All with $\dx = 3/2500$.\label{fig59}}
\end{figure}

\begin{figure}[ht!]
\centerline{\includegraphics[trim=1.0cm 2.9cm 0.9cm 2.3cm, clip, width=0.34\textwidth]{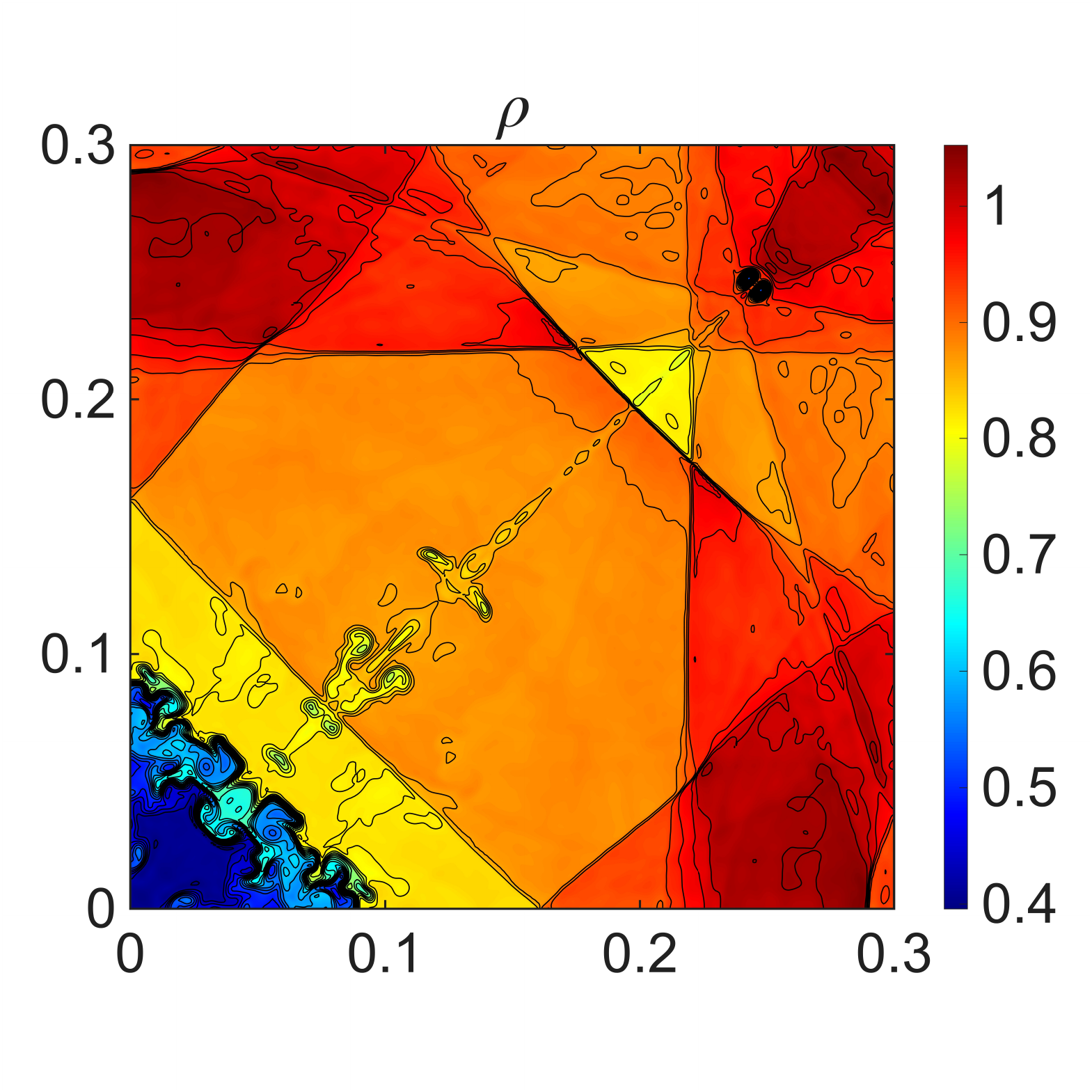}
\hspace*{0.5cm}
	    \includegraphics[trim=1.0cm 2.9cm 0.9cm 2.3cm, clip, width=0.34\textwidth]{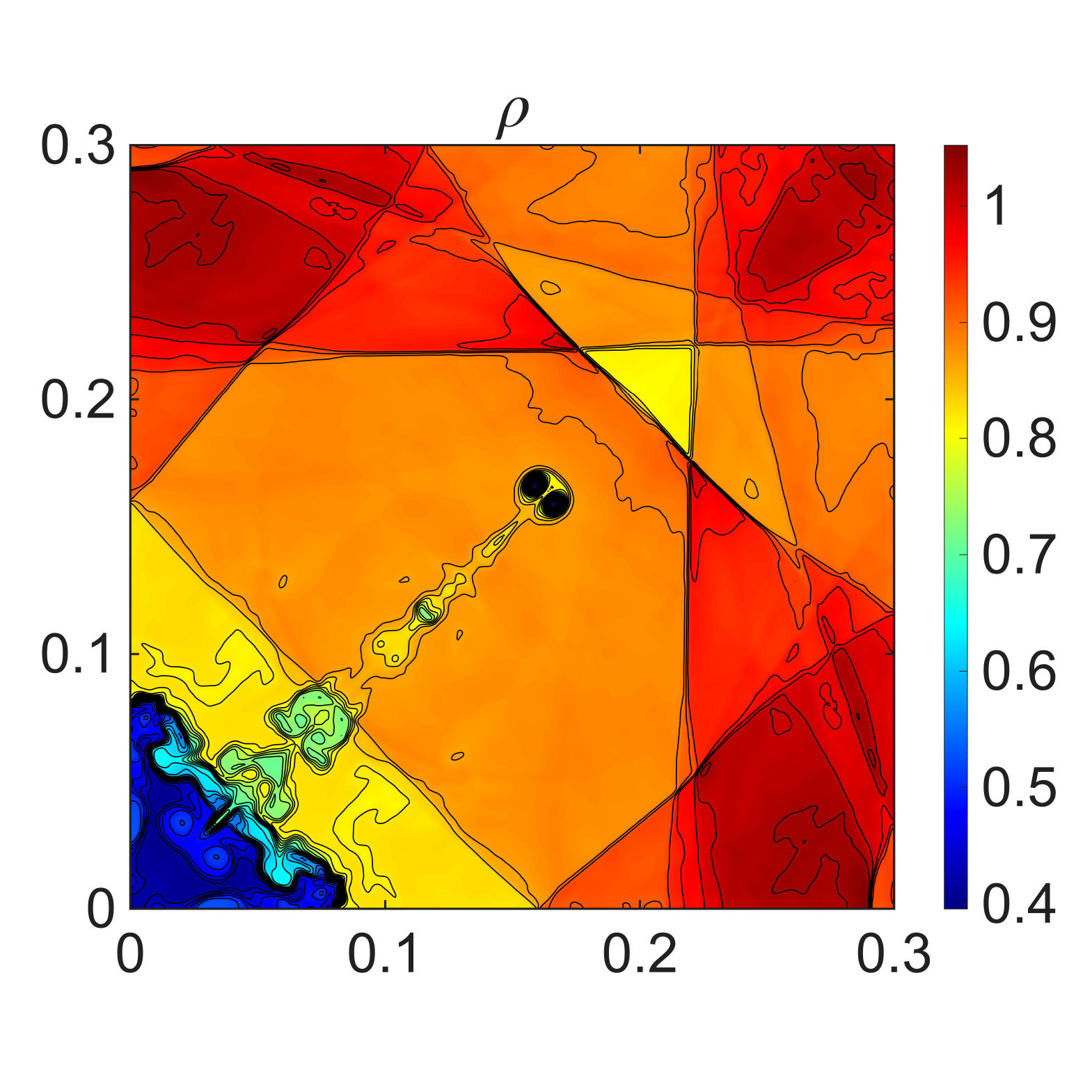}}
\vskip10pt
\hspace*{-0.3cm}\centerline{\includegraphics[trim=0.9cm 1.9cm 1.7cm 1.0cm, clip, width=0.30\textwidth]{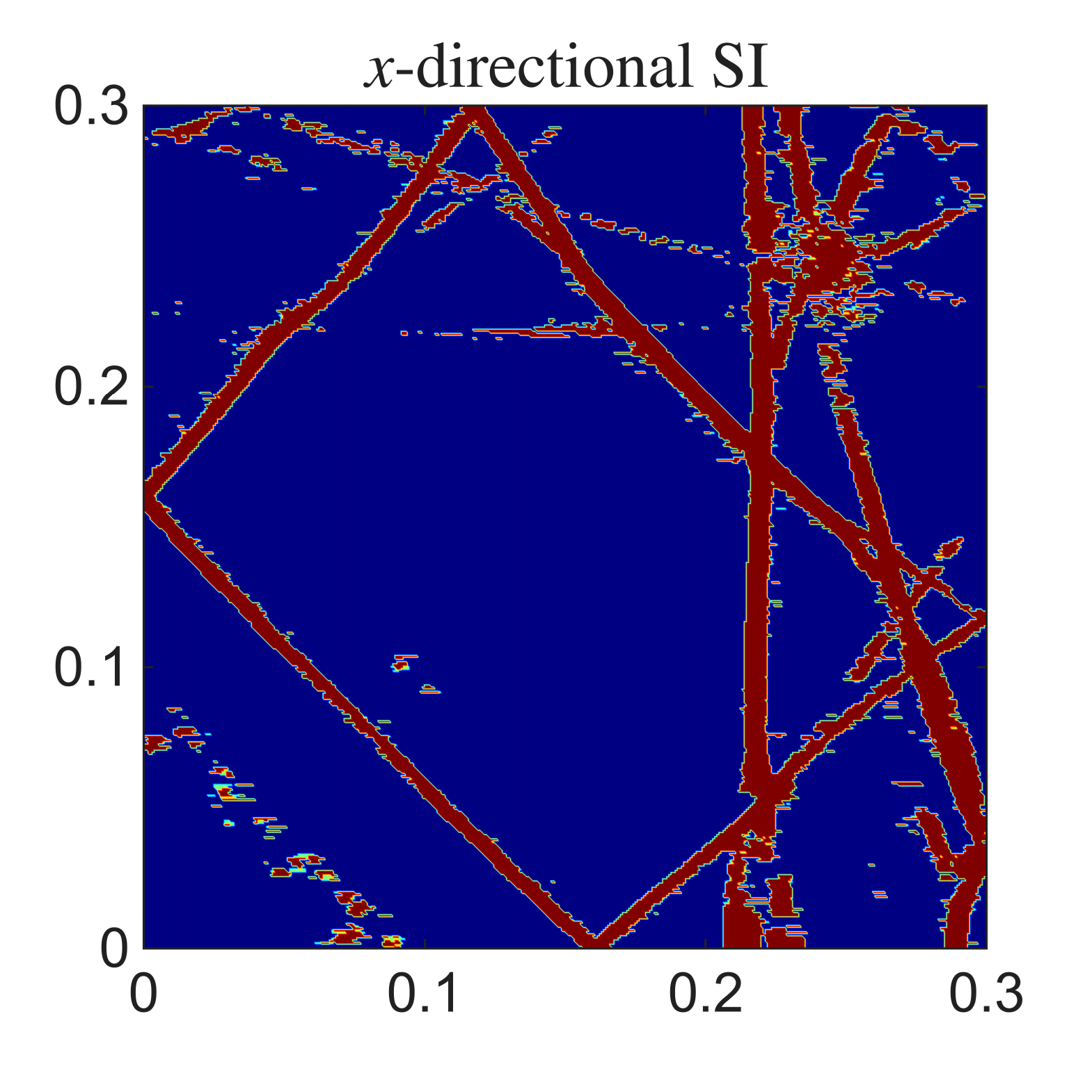}\hspace*{1.1cm}
	\includegraphics[trim=0.9cm 1.9cm 1.7cm 1.0cm, clip, width=0.30\textwidth]{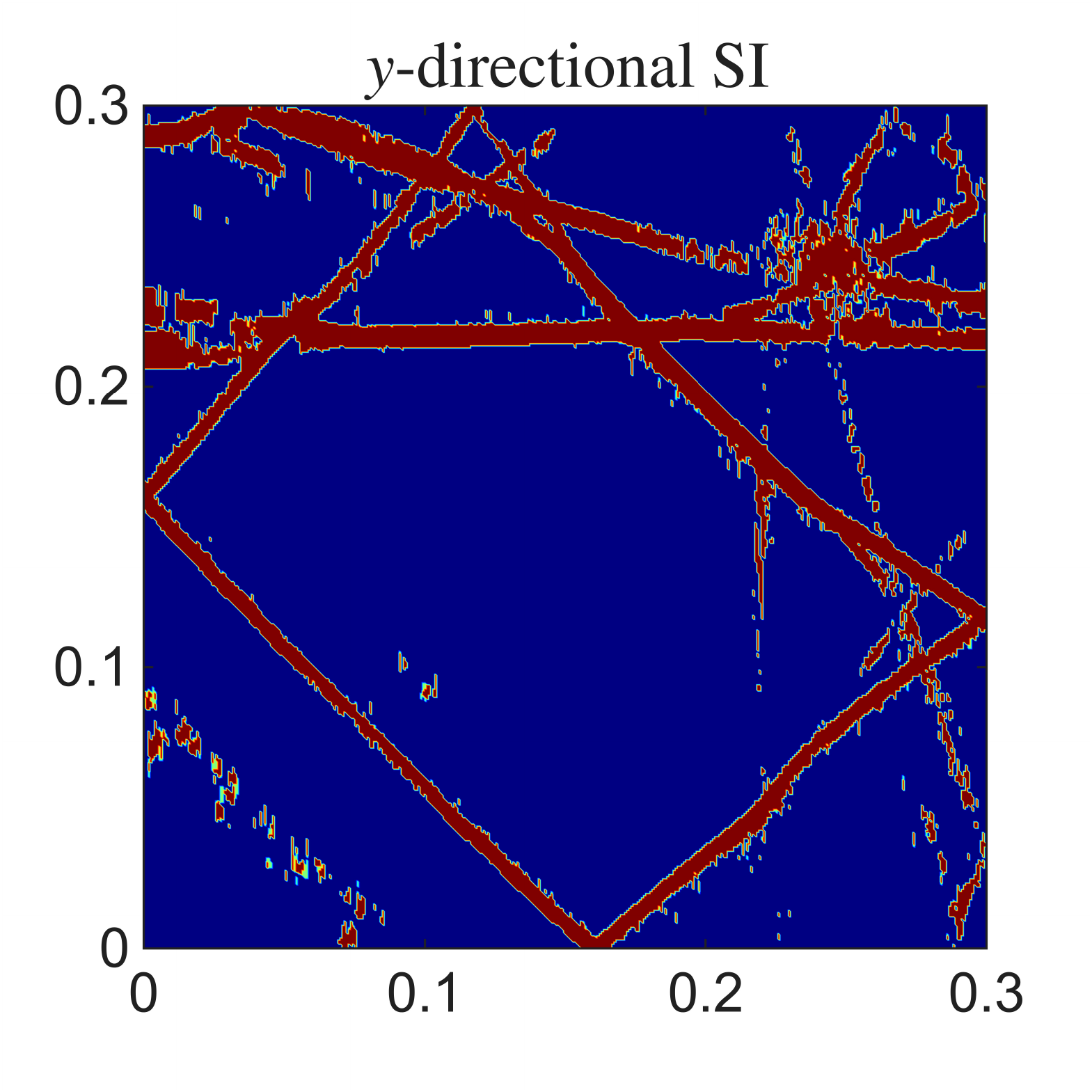}}
\caption{\sf Example 5: Same as in Figure \ref{fig59} but with $\dx = 3/4000$.\label{fig510}}
\end{figure}

We compute the numerical solutions by the adaptive (with the adaption coefficients $\kappa_{\rho u}=\kappa_{\rho v}=5\cdot10^{-2}$ and
$\kappa_p=2\cdot10^{-2}$) and A-WENO schemes until the final time $t=2.5$ on two uniform meshes with $\dx=\dy=3/2500$ and $3/4000$, and plot
the obtained solutions in Figures \ref{fig59} and \ref{fig510}, respectively. As one can observe, the jet in the adaptive solution reaches
further along the diagonal with sharper interfaces compared to the result obtained by the A-WENO scheme, indicating that the adaptive
strategy effectively reduces numerical dissipation. These results demonstrate that the proposed adaptive scheme achieves superior resolution
in capturing fine-scale flow structures and high-gradient phenomena without introducing spurious oscillations. In addition, one can see that
Regions S, RC, and RNC are adequately indicated.

\smallskip
\paragraph{Example 6 (Rayleigh-Taylor Instability)}
In the final example, we consider the RT instability, which arises at the interface between heavier and lighter fluids under the action of
gravity. To incorporate gravitational effects, appropriate source terms are added to the vertical momentum and total energy equations. As a
result, the governing equations become the 2-D compressible Euler equations with gravity, written as
$$
\U_t+\F(\U)_x+\G(\U)_y=\bm S,
$$
where $\U$, $\F$ and $\G$ are given by \eref{1.3} and $\bm S=(0,0,\rho,\rho v)^\top$.

The initial data,
\begin{equation*}
(\rho,u,v,p)\left|_{(x,y,0)}\right.=\begin{cases}(2,0,-0.025c\cos(8\pi x),2y+1),&y<0.5,\\
(1,0,-0.025c\cos(8\pi x),y+1.5),&y>0.5,\end{cases}
\end{equation*}
with $c$ being the speed of sound, are prescribed in the computational domain $[0,0.25]\times[0,1]$ subject to the solid wall boundary 
conditions at the left and right boundaries and the following Dirichlet boundary conditions at the top and bottom boundaries:
$$
(\rho,u,v,p)\left|_{(x,1,t)}\right.=(1,0,0,2.5),\quad(\rho,u,v,p)\left|_{(x,0,t)}\right.=(2,0,0,1).
$$

\begin{figure}[ht!]
\centerline{\includegraphics[trim=0.4cm 1.6cm 0.9cm 0.9cm, clip, width=0.22\textwidth]{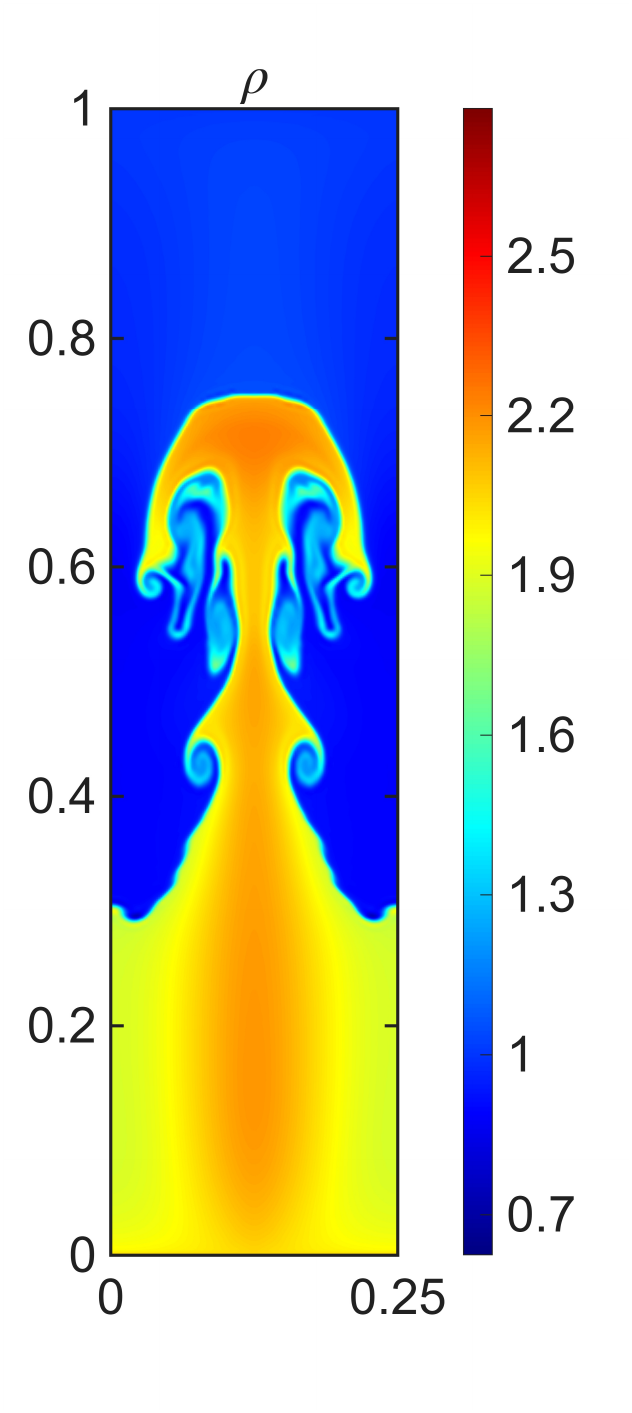}
\hspace*{0.6cm}
          \includegraphics[trim=0.4cm 1.6cm 0.9cm 0.9cm, clip, width=0.22\textwidth]{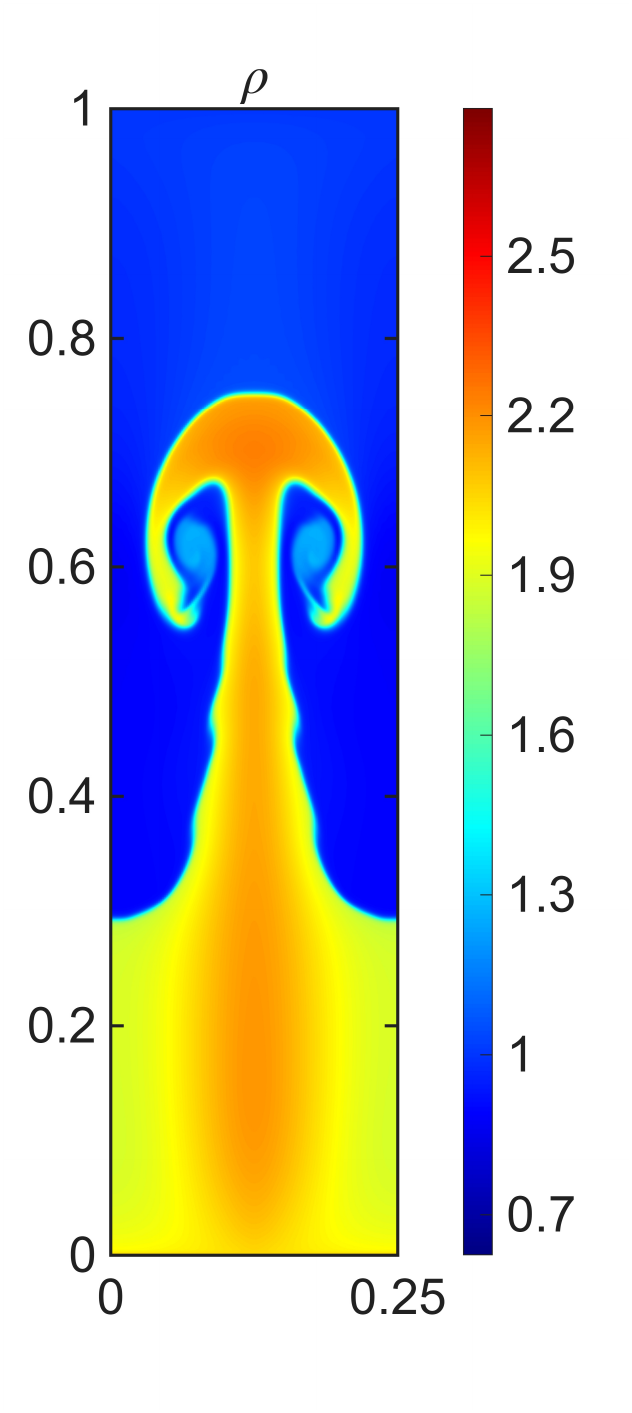}
\hspace*{0.2cm}
          \includegraphics[trim=0.4cm 1.6cm 0.9cm 0.9cm, clip, width=0.22\textwidth]{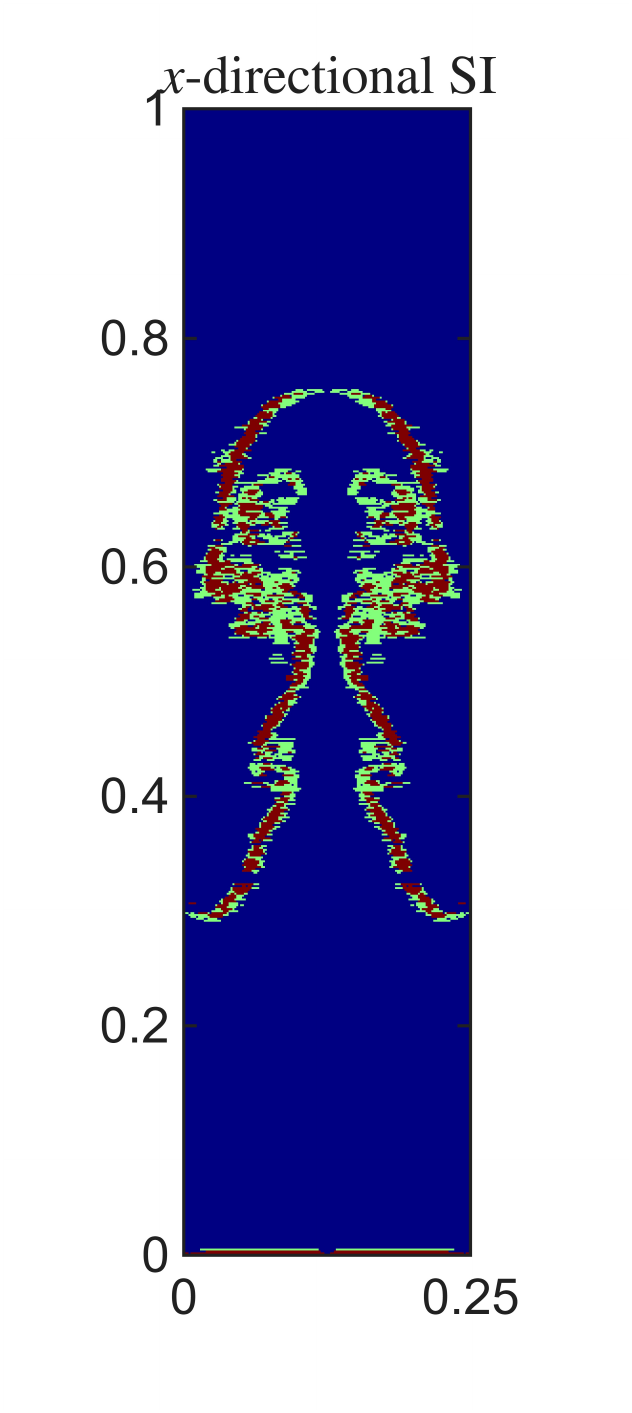}
\hspace*{0.0cm}
	  \includegraphics[trim=0.4cm 1.6cm 0.9cm 0.9cm, clip, width=0.22\textwidth]{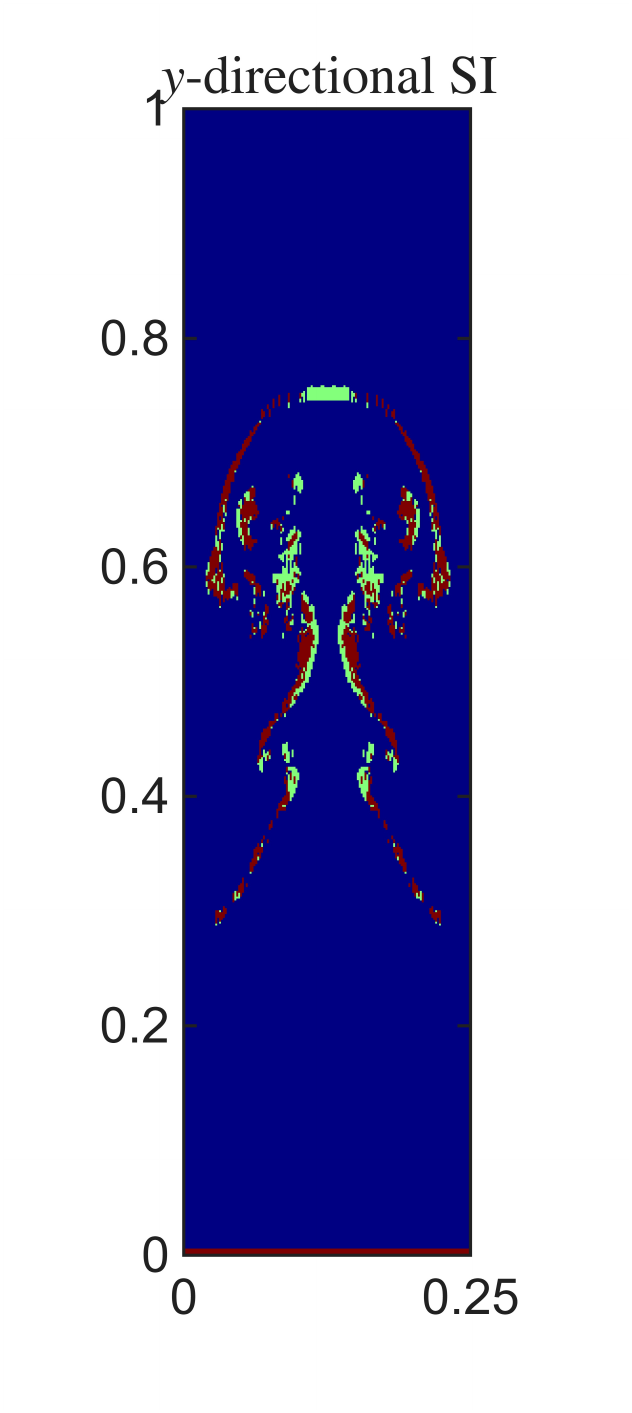}}
\vskip10pt
\centerline{\includegraphics[trim=0.4cm 1.6cm 0.9cm 0.9cm, clip, width=0.22\textwidth]{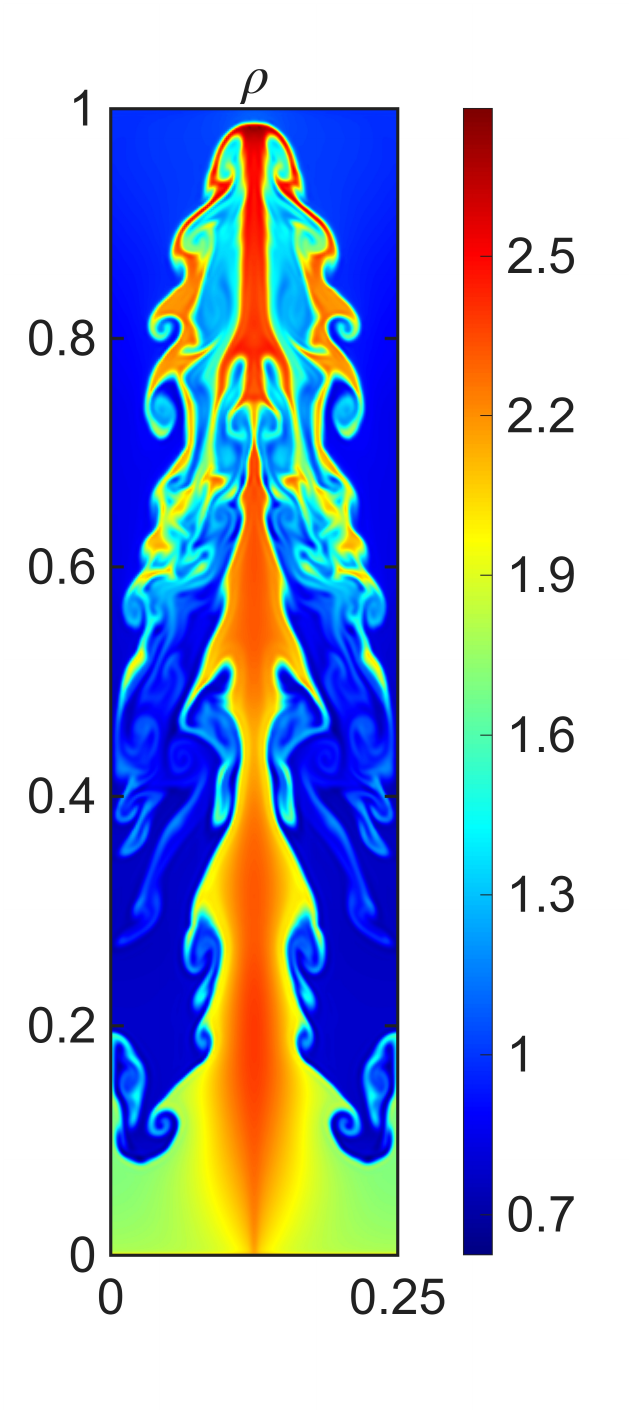}
\hspace*{0.6cm}
          \includegraphics[trim=0.4cm 1.6cm 0.9cm 0.9cm, clip, width=0.22\textwidth]{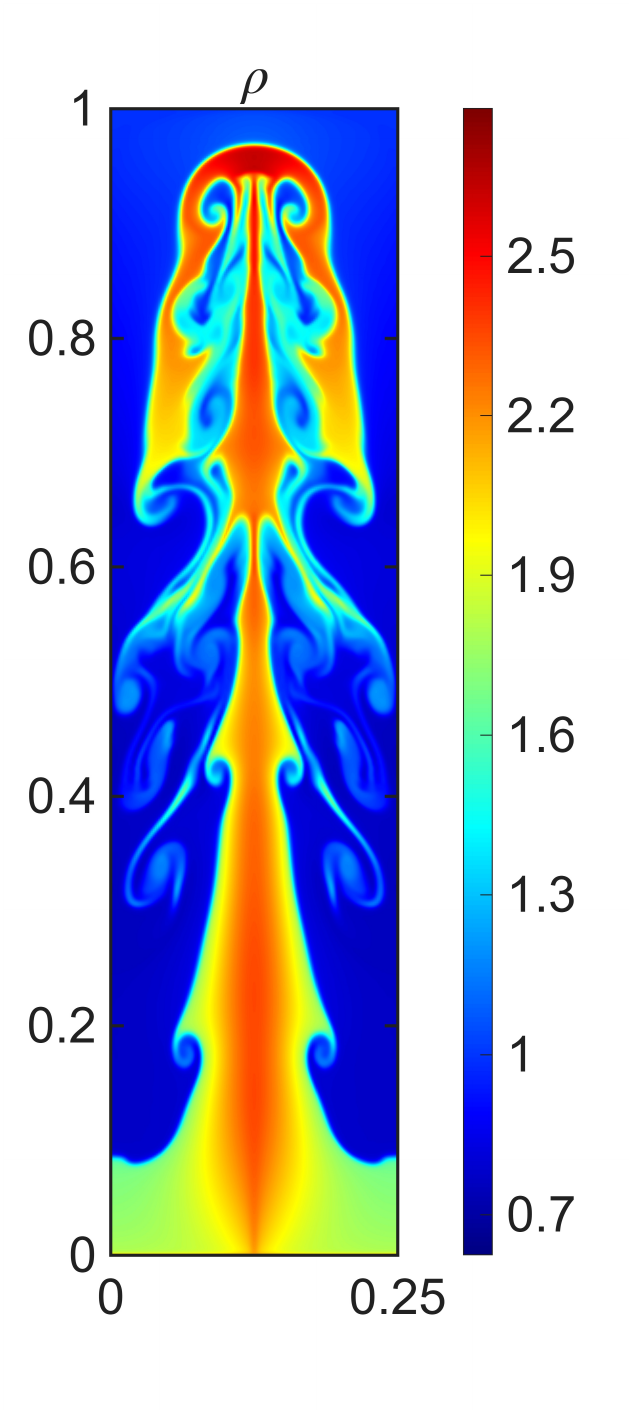}
\hspace*{0.2cm}
	  \includegraphics[trim=0.4cm 1.6cm 0.9cm 0.9cm, clip, width=0.22\textwidth]{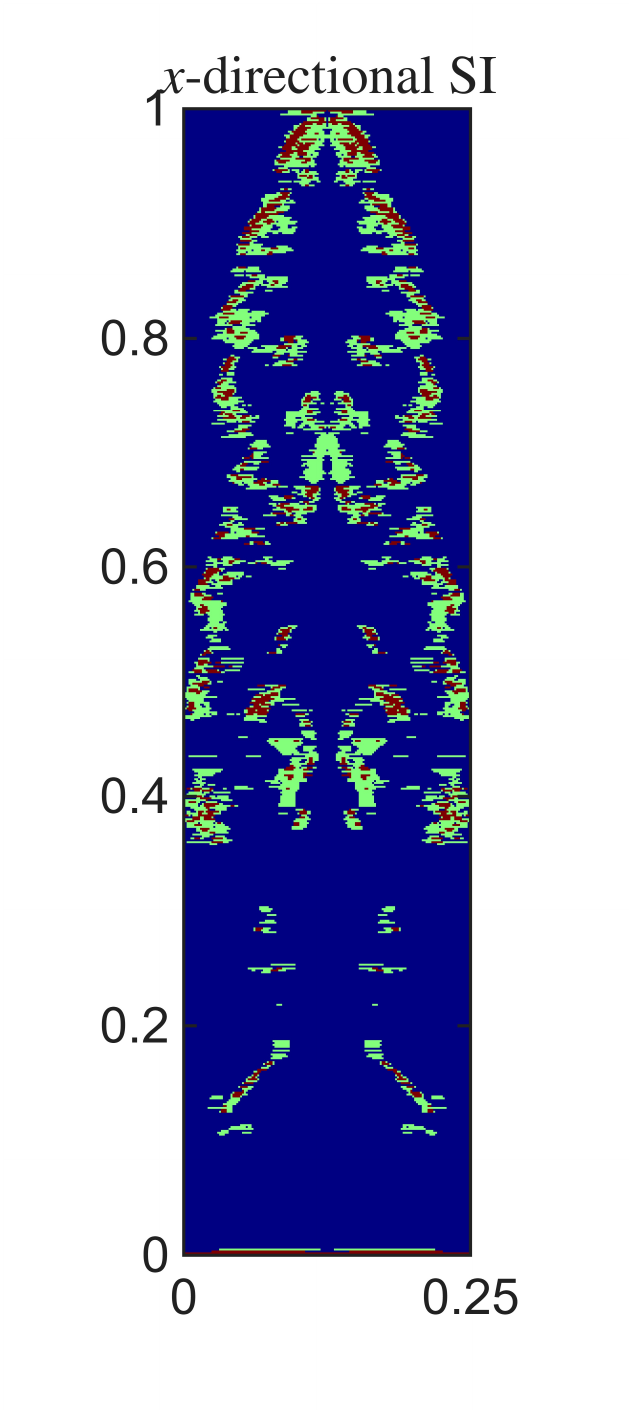}
	  \includegraphics[trim=0.4cm 1.6cm 0.9cm 0.9cm, clip, width=0.22\textwidth]{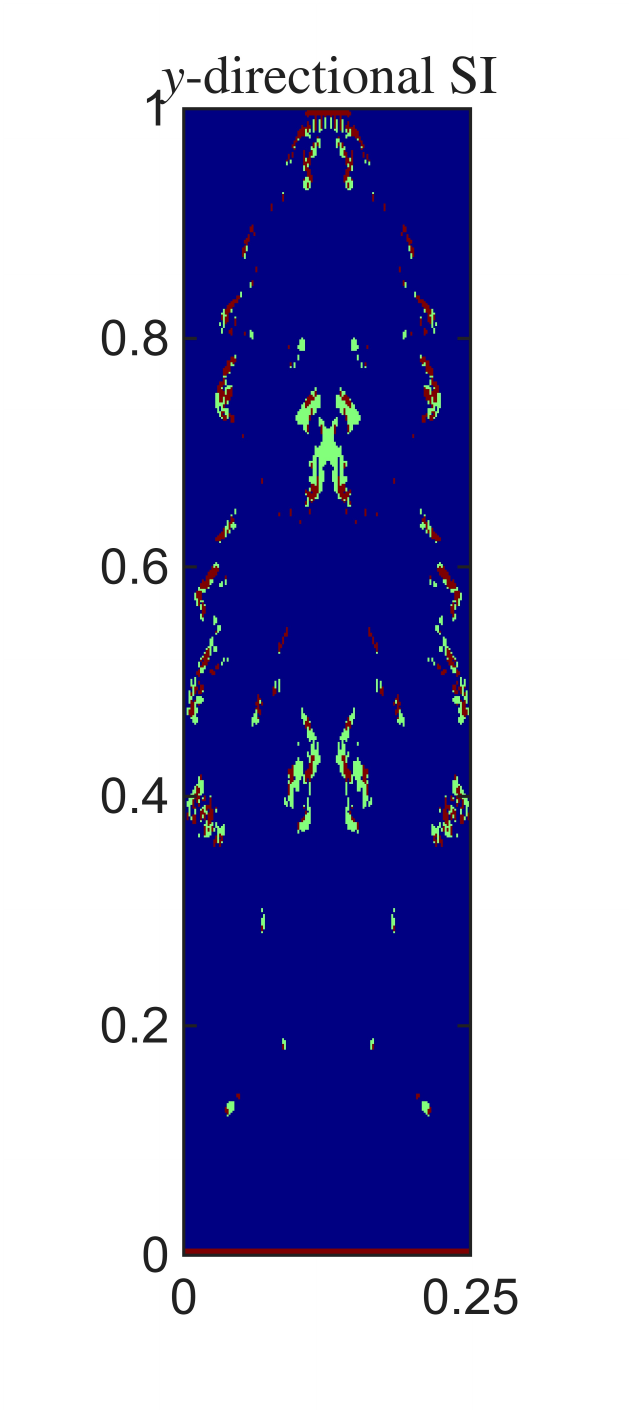}}
\caption{\sf Example 6: Density $\rho$ computed by the adaptive (left column) and A-WENO (second column) schemes, along with SIs in the
$x$-direction (third column) and the $y$-direction (right column) at times $t=1.95$ (upper row) and at=$2.95$ (bottom row).\label{fig511}}
\end{figure}

We compute the numerical solutions by the adaptive (with the adaption coefficients $\kappa_{\rho u}=\kappa_{\rho v}=1.2$ and $\kappa_p=1$)
and A-WENO schemes until the final time $t=2.95$ on a uniform mesh with $\dx=\dy=1/600$. The density computed at times $t=1.95$ and $2.95$
is shown in Figure \ref{fig511}, in which one can observe that the adaptive scheme resolves finer features, including Kelvin–Helmholtz
roll-ups in the falling fluid, which are significantly smeared or entirely absent in the A-WENO results. Moreover, the SIs shown in the
right columns demonstrate that the adaptive scheme accurately detects regions with strong gradients, thereby enabling sharper resolution of
complex interfaces and vortical structures.

\section{Conclusions}\label{sec6}
In this paper, we have developed new adaptive high-order numerical methods for hyperbolic systems of conservation laws within a dual
formulation (DF) framework. The key idea is to simultaneously evolve, at each time step, solutions of both the conservative and
nonconservative formulations of the same governing system and exploiting the discrepancy between the resulting numerical solutions to
develop reliable adaptive criteria. In particular, the smoothness indicator (SI) introduced in this work relies on the difference between
the numerical solutions computed by the conservative and nonconservative schemes. While the two solutions remain close up to truncation
errors in smooth regions, their difference becomes ${\cal O}(1)$ in nonsmooth regions. For the Euler equations of gas dynamics, the use of
momentum- and pressure-based SIs makes it possible to distinguish between the neighborhoods of contact discontinuities and other ``rough''
parts of the computed solutions. This information is then used to adaptively select appropriate numerical discretizations across the
computational domain. In the vicinities of contact discontinuities, we employ the low-dissipation central-upwind numerical flux and a
second-order piecewise linear reconstruction with the slopes computed using an overcompressive SBM limiter. Elsewhere, we use an alternative
weighted essentially non-oscillatory (A-WENO) framework with the central-upwind finite-volume numerical fluxes and either unlimited (in
smooth regions) or Ai-WENO-Z (in the nonsmooth regions away from contact discontinuities) fifth-order interpolation.

Extensive one- and two-dimensional numerical experiments confirm that the proposed adaptive approach achieves a favorable balance between
accuracy and efficiency. Compared with classical non-adaptive high-order schemes, the new method provides an improved resolution of complex
solution structures at a significantly reduced computational cost. Overall, the results demonstrate that the DF-based adaptive paradigm
offers a flexible and effective framework for the numerical simulation of problems involving intricate interactions between smooth and
nonsmooth solution features. The reliance on primitive formulations of the governing equations turns out to be particularly useful in many
contexts. Our DF-approach offers a flexible way to exploit it, and further applications in various directions, including compressible
multifluid flows and development of asymptotic-preserving schemes for a variety of hyperbolic systems of conservation laws, are currently
under investigation by the authors.


\bibliographystyle{siam}
\bibliography{Reference_UV_adaptation}

\end{document}